\documentclass[final]{siamltex}
\usepackage{graphicx}
\usepackage{amsfonts}
\usepackage{amsmath}

\usepackage[square,comma,numbers,sort&compress]{natbib}

\usepackage{hyperref}

\def\N{{\mathbb N}}
\def\R{{\mathbb R}}

\def\Z{{\mathbb Z}}

\def\ve{\varepsilon}

\def\pa{{\partial\Omega}}

\def\paD{{\partial D}}

\def\ce{{\rm ce}}
\def\se{{\rm se}}
\def\Mc{{\rm Mc}}
\def\Ms{{\rm Ms}}

\DeclareMathOperator{\esssup}{ess\,sup}

\pdfoutput=1

\title{Geometrical structure of Laplacian eigenfunctions}

\author{
D.~S. Grebenkov\footnotemark[1]\ \footnotemark[3] \ \footnotemark[4] \ \footnotemark[5]
\and B.-T. Nguyen\footnotemark[1]}

\begin{document}
\maketitle

\renewcommand{\thefootnote}{\fnsymbol{footnote}}
\footnotetext[1]{Laboratoire de Physique de la Mati\`ere Condens\'ee, CNRS -- Ecole Polytechnique, 91128 Palaiseau, France}
\footnotetext[3]{Laboratoire Poncelet, CNRS -- Independent University of Moscow, Bolshoy Vlasyevskiy Pereulok 11, 119002 Moscow, Russia}
\footnotetext[4]{Chebyshev Laboratory, Saint Petersburg State University, 14th line of Vasil'evskiy Ostrov 29, Saint Petersburg, Russia}
\footnotetext[5]{Corresponding author: denis.grebenkov@polytechnique.edu}
\renewcommand{\thefootnote}{\arabic{footnote}}

\begin{abstract}
We summarize the properties of eigenvalues and eigenfunctions of the
Laplace operator in bounded Euclidean domains with Dirichlet, Neumann
or Robin boundary condition.  We keep the presentation at a level
accessible to scientists from various disciplines ranging from
mathematics to physics and computer sciences.  The main focus is put
onto multiple intricate relations between the shape of a domain and
the geometrical structure of eigenfunctions.
\end{abstract}




\pagestyle{myheadings}
\thispagestyle{plain}


\hfill {\it Dedicated to Professor Bernard Sapoval for his 75th birthday}

\section{Introduction}

This review focuses on the classical eigenvalue problem for the
Laplace operator $\Delta = \partial^2/\partial x_1^2 + ... +
\partial^2/\partial x_d^2$ in an open bounded connected domain
$\Omega\subset \R^d$ ($d = 2,3,...$ being the space dimension),
\begin{equation}
\label{eq:eigen}
- \Delta u_m(x) = \lambda_m u_m(x) \quad (x\in\Omega) ,
\end{equation}
with Dirichlet, Neumann or Robin boundary condition on a piecewise
smooth boundary $\pa$:
\begin{equation}
\label{eq:BC}
\begin{split}
u_m(x) & = 0  \quad (x\in \pa)  \quad \mathrm{(Dirichlet)}, \\
\frac{\partial}{\partial n} u_m(x) & = 0  \quad (x\in \pa) \quad  \mathrm{(Neumann)}  , \\
\frac{\partial}{\partial n} u_m(x) + h u_m(x) & = 0 \quad (x\in \pa) \quad  \mathrm{(Robin)}, \\
\end{split}
\end{equation}
where $\partial/\partial n$ is the normal derivative pointed outwards
the domain, and $h$ is a posi\-tive constant.  The spectrum of the
Laplace operator is known to be discrete, the eigenvalues $\lambda_m$
are nonnegative and ordered in an ascending order by the index $m =
1,2,3,...$,
\begin{equation}
(0 \leq) \lambda_1 < \lambda_2 \leq \lambda_3 \leq \ldots \nearrow \infty 
\end{equation}  
(with possible multiplicities), while the eigenfunctions $\{u_m(x)\}$
form a complete basis in the functional space $L_2(\Omega)$ of
measurable and square-integrable functions on $\Omega$
\cite{Courant,Reed}.  By definition, the function $0$ satisfying
Eqs. (\ref{eq:eigen}, \ref{eq:BC}) is excluded from the set of
eigenfunctions.  Since the eigenfunctions are defined up to a
multiplicative factor, it is sometimes convenient to normalize them to
get the unit $L_2$-norm:
\begin{equation}
\|u_m\|_{2} \equiv \|u_m\|_{L_2(\Omega)} \equiv \left(\int\limits_\Omega dx~ |u_m(x)|^2\right)^{1/2} = 1
\end{equation}
(note that there is still ambiguity up to the multiplication by
$e^{i\alpha}$, with $\alpha\in\R$).

Laplacian eigenfunctions appear as vibration modes in acoustics, as
electron wave functions in quantum waveguides, as natural basis for
constructing heat kernels in the theory of diffusion, etc.  For
instance, vibration modes of a thin membrane (a drum) with a fixed
boundary are given by Dirichlet Laplacian eigenfunctions $u_m$, with
the drum frequencies proportional to $\sqrt{\lambda_m}$
\cite{Rayleigh}.  A particular eigenmode can be excited at the
corresponding frequency \cite{Sapoval91,Sapoval93,Sapoval97}.  In the
theory of diffusion, an interpretation of eigenfunctions is less
explicit.  The first eigenfunction represents the long-time asymptotic
spatial distribution of particles diffusing in a bounded domain (see
below).  A conjectural probabilistic representation of higher-order
eigenfunctions through a Fleming-Viot type model was developed by
Burdzy {\it et al.} \cite{Burdzy96,Burdzy00}.

The eigenvalue problem (\ref{eq:eigen}, \ref{eq:BC}) is archetypical
in the theory of elliptic operators, while the properties of the
underlying eigenfunctions have been thoroughly investigated in various
mathematical and physical disciplines, including spectral theory,
probability and stochastic processes, dynamical systems and quantum
billiards, condensed matter physics and quantum mechanics, theory of
acoustical, optical and quantum waveguides, computer sciences, etc.
Many books and reviews were dedicated to different aspects of
Laplacian eigenvalues, eigenfunctions and their applications (see,
e.g.,
\cite{Payne67,Polya,Bandle,Hile80,Kuttler84,Chavel84,Davies2,Edmunds,Jakobson01,Ashbaugh00,Hale05,Ashbaugh07,Arendt09,Benguria11}).
The diversity of notions and methods developed by mathematicians,
physicists and computer scientists often makes the progress in one
discipline almost unknown or hardly accessible to scientists from the
other disciplines.  One of the goals of the review is to bring
together various facts about Laplacian eigenvalues and eigenfunctions
and to present them at a level accessible to scientists from various
disciplines.  For this purpose, many technical details and
generalities are omitted in favor to simple illustrations.  While the
presentation is focused on the Laplace operator in bounded Euclidean
domains with piecewise smooth boundaries, a number of extensions are
relatively straightforward.  For instance, the Laplace operator can be
extended to a second order elliptic operator with appropriate
coefficients, the piecewise smoothness of a boundary can often be
relaxed \cite{Lions,Grisvard}, while Euclidean domains can be replaced
by Riemannian manifolds or weighted graphs \cite{Jakobson01}.  The
main emphasis is put onto the geometrical structure of Laplacian
eigenfunctions and on their relation to the shape of a domain.
Although the bibliography counts more than five hundred citations, it
is far from being complete, and readers are invited to refer to other
reviews and books for further details and references.

The review is organized as follows.  We start by recalling in
Sec. \ref{sec:general} general properties of the Laplace operator.
Explicit representations of eigenvalues and eigenfunctions in simple
domains are summarized in Sec. \ref{sec:simple}.  In
Sec. \ref{sec:eigenvalues} we review the properties of eigenvalues and
their relation to the shape of a domain: Weyl's asymptotic law,
isoperimetric inequalities and the related shape optimization
problems, and Kac's inverse spectral problem.  Although eigenfunctions
are not involved at this step, valuable information can be learned
about the domain from the eigenvalues alone.  The next step consists
in the analysis of nodal lines/surfaces or nodal domains in
Sec. \ref{sec:nodal}.  The nodal lines tell us how the zeros of
eigenfunctions are spatially distributed, while their amplitudes are
still ignored.  In Sec. \ref{sec:estimates}, several estimates for the
amplitudes of eigenfunctions are summarized.  Most of these results
were obtained during the last twenty years.

Section \ref{sec:localization} is devoted to the property of
eigenfunctions known as localization.  We start by recalling the
notion of localization in quantum mechanics: the strong localization
by a potential (Sec. \ref{sec:potential}), Anderson localization
(Sec. \ref{sec:Anderson}) and trapped modes in infinite waveguides
(Sec. \ref{sec:trapping}).  In all three cases, the eigenvalue problem
is different from Eqs. (\ref{eq:eigen}, \ref{eq:BC}), due to either
the presence of a potential, or the unboundness of a domain.
Nevertheless, these cases are instructive, as similar effects may be
observed for the eigenvalue problem (\ref{eq:eigen}, \ref{eq:BC}).  In
particular, we discuss in Sec. \ref{sec:expon} an exponentially
decaying upper bound for the norm of eigenfunctions in domains with
branches of variable cross-sectional profiles.  Section
\ref{sec:dumbbell} reviews the properties of low-frequency
eigenfunctions in ``dumbbell'' domains, in which two (or many)
subdomains are connected by narrow channels.  This situation is
convenient for a rigorous analysis as the width of channels plays the
role of a small parameter \cite{Sanchez}.  A number of asymptotic
results for eigenvalues and eigenfunctions were derived, for
Dirichlet, Neumann and Robin boundary conditions.  A harder case of
irregular or fractal domains is discussed in Sec. \ref{sec:irregular}.
Here, it is difficult to identify a relevant small parameter to get
rigorous estimates.  In spite of numerous numerical examples of
localized eigenfunctions (both for Dirichlet and Neumann boundary
conditions), a comprehensive theory of localization is still missing.
Section \ref{sec:high-freq} is devoted to high-frequency localization
and the related scarring problems in quantum billiards.  We start by
illustrating the classical whispering gallery, bouncing ball and
focusing modes in circular and elliptical domains.  We also provide
examples for the case without localization.  A brief overview of
quantum billiards is presented.  In the last
Sec. \ref{sec:conclusion}, we mention some issues that could not be
included into the review, e.g., numerical methods for computation of
eigenfunctions or their numerous applications.

\section{Basic properties}
\label{sec:general}

We start by recalling basic properties of the Laplacian eigenvalues
and eigenfunctions (see \cite{Courant,Reed,Birman} or other standard
textbooks).

(i) The eigenfunctions are infinitely differentiable inside the domain
$\Omega$.  For any open subset $V\subset \Omega$, the restriction of
$u_m$ on $V$ cannot be strictly $0$ \cite{Kuttler84}.

(ii) Multiplying Eq. (\ref{eq:eigen}) by $u_m$, integrating over
$\Omega$ and using the Green's formula yield
\begin{equation}
\label{eq:lambda_u}
\lambda_m = \frac{\int\limits_\Omega dx~ |\nabla u_m|^2 - \int\limits_{\pa} dx ~u_m \frac{\partial u_m}{\partial n}}
{\int\limits_\Omega dx ~u_m^2} = \frac{\| \nabla u_m\|_{L_2(\Omega)}^2 + h \| u_m\|_{L_2(\pa)}^2}{\| u_m\|_{L_2(\Omega)}^2} ,
\end{equation}
where $\nabla$ stands for the gradient operator, and we used Robin
boundary condition (\ref{eq:BC}) in the last equality; for Dirichlet
or Neumann boundary conditions, the boundary integral (second term)
vanishes.  This formula ensures that all eigenvalues are nonnegative.

(iii) Similar expression appears in the variational formulation of the
eigenvalue problem, known as the minimax principle \cite{Courant}
\begin{equation}
\label{eq:minimax}
\lambda_m = \min \max \frac{\| \nabla v\|_{L_2(\Omega)}^2 + h \| v\|_{L_2(\pa)}^2}{\| v\|_{L_2(\Omega)}^2} ,
\end{equation}
where the maximum is over all linear combinations of the form
\begin{equation*}
v = a_1 \phi_1 + ... + a_m \phi_m ,
\end{equation*}
and the minimum is over all choices of $m$ linearly independent
continuous and piecewise-differentiable functions $\phi_1$, ...,
$\phi_m$ (said to be in the Sobolev space $H^1(\Omega)$)
\cite{Courant,Henrot}.  Note that the minimum is reached exactly on
the eigenfunction $u_m$.  For Dirichlet eigenvalue problem, there is a
supplementary condition $v = 0$ on the boundary $\pa$ so that the
second term in Eq. (\ref{eq:minimax}) is canceled.  For Neumann
eigenvalue problem, $h = 0$ and the second term vanishes again.

(iv) The minimax principle implies the monotonous increase of the
eigenvalues $\lambda_m$ with $h$, namely if $h < h'$, then
$\lambda_m(h) \leq \lambda_m(h')$.  In particular, any eigenvalue
$\lambda_m(h)$ of the Robin problem lies between the corresponding
Neumann and Dirichlet eigenvalues.

\begin{figure}
\begin{center}
\includegraphics[width=50mm]{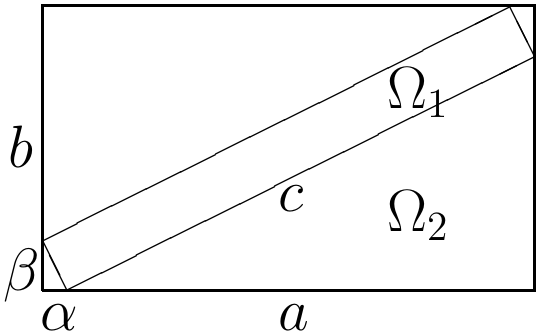}
\end{center}
\caption{
A counter-example for the property of domain monotonicity for Neumann
boundary condition.  Although a smaller rectangle $\Omega_1$ is
inscribed into a larger rectangle $\Omega_2$ (i.e., $\Omega_1\subset
\Omega_2$), the second eigenvalue $\lambda_2(\Omega_1) = \pi^2/c^2$
is smaller than the second eigenvalue $\lambda_2(\Omega_2) =
\pi^2/a^2$ (if $a > b$) when $c = \sqrt{(a-\alpha)^2 + (b-\beta)^2} >
a$ (courtesy by N. Saito; see also \cite{Henrot}, Sec. 1.3.2). }
\label{fig:N_counter}
\end{figure}

(v) For Dirichlet boundary condition, the minimax principle implies
the property of domain monotonicity: eigenvalues monotonously decrease
when the domain enlarges, i.e., $\lambda_m(\Omega_1) \geq
\lambda_m(\Omega_2)$ if $\Omega_1 \subset \Omega_2$.  This property
does not hold for Neumann or Robin boundary conditions, as illustrated
by a simple counter-example on Fig. \ref{fig:N_counter}.

(vi) The eigenvalues are invariant under translations and rotations of
the domain.  This is a key property for an efficient image recognition
and analysis \cite{Reuter06,Saito08,Saito09}.  When a domain is
expanded by factor $\alpha$, all the eigenvalues are rescaled by
$1/\alpha^2$.

(vii) The first eigenfunction $u_1$ does not change the sign and can
be chosen positive.  Because of the orthogonality of eigenfunctions,
$u_1$ is in fact the only eigenfunction not changing its sign.

(viii) The first eigenvalue $\lambda_1$ is simple and strictly
positive for Dirichlet and Robin boundary conditions; for Neumann
boundary condition, $\lambda_1 = 0$ and $u_1$ is a constant.

(ix) The completeness of eigenfunctions in $L_2(\Omega)$ can be
expressed as
\begin{equation}
\sum\limits_m  u_m(x) u_m^*(y) = \delta(x-y)  \qquad (x,y\in\Omega),
\end{equation}
where asterisk denotes the complex conjugate, $\delta(x)$ is the Dirac
distribution, and the eigenfunctions are $L_2$-normalized.
Multiplying this relation by a function $f\in L_2(\Omega)$ and
integrating over $\Omega$ yields the decomposition of $f(x)$ over
$u_m(x)$:
\begin{equation*}
f(x) = \sum\limits_m u_m(x) \int\limits_\Omega dy ~f(y) ~ u_m^*(y) .
\end{equation*}

(x) The Green function $G(x,y)$ for the Laplace operator which
satisfies
\begin{equation}
-\Delta G(x,y) = \delta(x-y)  \qquad (x,y\in \Omega)
\end{equation}
(with an appropriate boundary condition), admits the spectral
decomposition over the $L_2$-normalized eigenfunctions
\begin{equation}
\label{eq:G_decomp}
G(x,y) = \sum\limits_m \lambda_m^{-1} u_m(x) u_m^*(y) .
\end{equation}
(for Neumann boundary condition, $\lambda_1 = 0$ has to be excluded;
in that case, the Green function is defined up to an additive
constant).

Similarly, the heat kernel (or diffusion propagator) $G_t(x,y)$
satisfying
\begin{equation}
\label{eq:heat}
\begin{split}
\frac{\partial}{\partial t} G_t(x,y) - \Delta G_t(x,y) & = 0  \qquad (x,y\in \Omega), \\
G_{t=0}(x,y) & = \delta(x-y)  \\
\end{split}
\end{equation}
(with an appropriate boundary condition), admits the spectral
decomposition
\begin{equation}
\label{eq:Gt_decomp}
G_t(x,y) = \sum\limits_m e^{-\lambda_m t} u_m(x) u_m^*(y) .
\end{equation} 
The Green function and heat kernel allow one to solve the standard
boundary value and Cauchy problems for the Laplace and heat equations,
respectively \cite{Crank,Carslaw}.  The decompositions
(\ref{eq:G_decomp}, \ref{eq:Gt_decomp}) are the major tool for getting
explicit solutions in simple domains for which the eigenfunctions are
known explicitly (see Sec. \ref{sec:simple}).  This representation is
also important for the theory of diffusion due to the probabilistic
interpretation of $G_t(x,y)dx$ as the conditional probability for
Brownian motion started at $y$ to arrive in the $dx$ vicinity of $x$
after a time $t$ \cite{Feller,Bass,Bass2,Garnett,Redner,Hughes,Weiss,Port,Borodin}.
Setting Dirichlet, Neumann or Robin boundary conditions, one can
respectively describe perfect absorptions, perfect reflections and
partial absorption/reflection on the boundary \cite{Grebenkov06}.

For Dirichlet boundary condition, if $\Omega \subset \Omega'$, then $0
\leq G_t^{(\Omega)}(x,y) \leq G_t^{(\Omega')}(x,y)$
\cite{vandenBerg89}.  In particular, taking $\Omega' = \R^d$, one gets
\begin{equation}
0 \leq G_t(x,y) \leq (4\pi t)^{-d/2} \exp\left(- \frac{|x-y|^2}{4t}\right) ,
\end{equation}
where the Gaussian heat kernel for free space is written on the
right-hand side.  The above domain monotonicity for heat kernels may
not hold for Neumann boundary condition \cite{Bass93}.

(xi) For Dirichlet boundary condition, the eigenvalues vary
continuously under a ``continuous'' perturbation of the domain
\cite{Courant}.  For Neumann boundary condition, the situation is much
more delicate.  The continuity still holds when a bounded domain with
a smooth boundary is deformed by a ``continuously differentiable
transformation'', while in general this statement is false, with an
explicit counter-example provided in \cite{Courant}.  Note that the
continuity of the spectrum is important for numerical computations of
the eigenvalues by finite element or other methods, in which an
irregular boundary is replaced by a suitable polygonal or piecewise
smooth approximation.  The underlying assumption that the eigenvalues
are very little affected by such domain perturbations, holds in great
generality for Dirichlet boundary condition, but is much less evident
for Neumann boundary condition \cite{Burenkov02}.  The spectral
stability of elliptic operators under domain perturbations has been
thoroughly investigated
\cite{Hale05,Henry,Burenkov02,Burenkov07,Burenkov08,Burenkov08b,Burenkov12}.
It is also worth stressing that the spectrum of the Laplace operator
in a bounded domain with Neumann boundary condition on an irregular
boundary may not be discrete, with explicit counter-examples provided
in \cite{Hempel91}.

\section{Eigenbasis for simple domains}
\label{sec:simple}

We list the examples of ``simple'' domains, in which symmetries allow
for variable separations and thus explicit representations of
eigenfunctions in terms of elementary or special functions.

\subsection{Intervals, rectangles, parallelepipeds}
\label{sec:u_rectangle}

For rectangle-like domains $\Omega = [0,\ell_1]\times ...\times
[0,\ell_d]\subset \R^d$ (with the sizes $\ell_i > 0$), the natural
variable separation yields
\begin{equation}
\begin{split}
u_{n_1,...,n_d}(x_1,...,x_d) & = u_{n_1}^{(1)}(x_1) \ldots u_{n_d}^{(d)}(x_d) , \qquad
\lambda_{n_1,...,n_d} = \lambda_{n_1}^{(1)} + ... + \lambda_{n_d}^{(d)} , \\
\end{split}
\end{equation}
where the multiple index $n_1...n_d$ is used instead of $m$, and
$u_{n_i}^{(i)}(x_i)$ and $\lambda_{n_i}^{(i)}$ ($i=1,...,d$)
correspond to the one-dimensional problem on the interval
$[0,\ell_i]$.  Depending on the boundary condition, $u_n^{(i)}(x)$ are
sines (Dirichlet), cosines (Neumann) or their combinations (Robin):
\begin{equation}
\begin{split}
u_n^{(i)}(x) & = \sin(\pi (n+1) x/\ell_i) , \quad \lambda_n^{(i)} = \pi^2 (n+1)^2/\ell_i^2, 
 \hskip 10mm  \rm{(Dirichlet)} , \\
u_n^{(i)}(x) & = \cos(\pi n x/\ell_i) , \hskip 12mm   \lambda_n^{(i)} = \pi^2 n^2/\ell_i^2, 
 \hskip 19mm \rm{(Neumann)} , \\
u_n^{(i)}(x) & = \sin(\alpha_n x/\ell_i) + \frac{\alpha_n}{h \ell_i} \cos(\alpha_n x/\ell_i) , \quad \lambda_n^{(i)} = \alpha_n^2/\ell_i^2, 
 \hskip 3mm \rm{(Robin)} , \\
\end{split}
\end{equation}
where $n = 0,1,2,...$ and the coefficients $\alpha_n$ depend on the
parameter $h$ and satisfy the equation obtained by imposing the Robin
boundary condition in Eq. (\ref{eq:BC}) at $x = \ell_i$:
\begin{equation}
\frac{2\alpha_n}{h\ell_i} \cos\alpha_n + \left(1 - \frac{\alpha_n^2}{h^2\ell_i^2}\right)\sin\alpha_n = 0 .
\end{equation}
According to the property (iv) of Sec. \ref{sec:general}, this
equation has the unique solution $\alpha_n$ on each interval
$[n\pi,(n+1)\pi]$ ($n = 0,1,2,...$), that makes its numerical
computation by bisection (or other) method easy and fast.  All the
eigenvalues $\lambda_n^{(i)}$ are simple (not degenerate), while
\begin{equation}
\| u_n^{(i)}(x) \|_{L_2((0,\ell_i))} = \left(\frac{\alpha_n^2 + 2h\ell_i + h^2\ell_i^2}{2h^2}\right)^{1/2} .
\end{equation}

In turn, the eigenvalues $\lambda_{n_1,...,n_d}$ can be degenerate if
there exists a rational ratio $(\ell_i/\ell_j)^2$ (with $i\ne j$).
For instance, the first Dirichlet eigenvalues of the unit square are
$2\pi^2$, $5\pi^2$, $5\pi^2$, $8\pi^2$, .... , with the twice
degenerate second eigenvalue.  An eigenfunction associated to a
degenerate eigenvalue is a linear combination of the corresponding
functions.  For the above example $u(x_1,x_2) = c_1
\sin(\pi x_1)\sin(2\pi x_2) + c_2 \sin(2\pi x_1)\sin(\pi x_2)$ with
any $c_1$ and $c_2$ such that $c_1^2+c_2^2 \ne 0$.

\subsection{Disk, sector and circular annulus}

The rotation symmetry of a circular annulus, $\Omega = \{ x\in\R^2 ~:~
R_0 < |x| < R\}$, allows one to write the Laplace operator in polar
coordinates,
\begin{equation}
\begin{cases}  x_1 = r \cos \varphi , \cr  x_2 = r \sin \varphi ,  \end{cases}  \qquad
\Delta = \frac{\partial^2}{\partial r^2} + \frac{1}{r} \frac{\partial}{\partial r} + \frac{1}{r^2} \frac{\partial^2}{\partial\varphi^2} ,
\end{equation}
that leads to variable separation and an explicit representation of
eigenfunctions
\begin{equation}
\label{eq:u_annulus}
u_{nkl}(r,\varphi) = \bigl[J_n(\alpha_{nk} r/R) + c_{nk} Y_n(\alpha_{nk} r/R)\bigr] 
\times \begin{cases} \cos(n\varphi), \quad l = 1,  \cr  \sin(n\varphi), \quad l = 2 ~ (n\ne 0),\end{cases}
\end{equation}
where $J_n(z)$ and $Y_n(z)$ are the Bessel functions of the first and
second kind \cite{Abramowitz,Watson,Bowman}, and the coefficients
$\alpha_{nk}$ and $c_{nk}$ are set by the boundary conditions at $r =
R$ and $r = R_0$:
\begin{equation}
\label{eq:BC_disk}
\begin{split}
0 & = \frac{\alpha_{nk}}{R} \biggl[J'_n(\alpha_{nk}) + c_{nk} Y'_n(\alpha_{nk})\biggr] 
+ h \biggl[J_n(\alpha_{nk}) + c_{nk} Y_n(\alpha_{nk})\biggr], \\
0 & = - \frac{\alpha_{nk}}{R} \biggl[J'_n\bigl(\alpha_{nk} \frac{R_0}{R}\bigr) + c_{nk} Y'_n\bigl(\alpha_{nk} \frac{R_0}{R}\bigr)\biggr] 
+ h \biggl[J_n\bigl(\alpha_{nk} \frac{R_0}{R}\bigr) + c_{nk} Y_n\bigl(\alpha_{nk} \frac{R_0}{R}\bigr)\biggr]. \\
\end{split}
\end{equation}
For each $n = 0,1,2,...$, the system of these equations has infinitely
many solutions $\alpha_{nk}$ which are enumerated by the index $k =
1,2,3,...$ \cite{Watson}.  The eigenfunctions are enumerated by the
triple index $nkl$, with $n = 0,1,2,...$ counting the order of Bessel
functions, $k = 1,2,3,...$ counting solutions of
Eqs. (\ref{eq:BC_disk}), and $l = 1,2$.  Since $u_{0k2}(r,\varphi)$
are trivially zero (as $\sin (n\varphi) = 0$ for $n = 0$), they are
excluded.  The eigenvalues $\lambda_{nk} = \alpha_{nk}^2/R^2$, which
are independent of the last index $l$, are simple for $n = 0$ and
twice degenerate for $n > 0$.  In the latter case, an eigenfunction is
any nontrivial linear combination of $u_{nk1}$ and $u_{nk2}$.  The
squared $L_2$-norm of the eigenfunction is
\begin{equation}
\begin{split}
\| u_{nkl}(r,\varphi) \|_2^2 = \frac{\pi (2-\delta_{n,0}) R^2}{2\alpha_{nk}^2} \biggl[ & \biggl(\alpha_{nk}^2 + h^2 R^2 - n^2\biggr) v_{nk}^2(R) \\
& - \biggl((\alpha_{nk}^2 + h^2 R^2)\frac{R_0^2}{R^2} - n^2\biggr) v_{nk}^2(R_0)\biggr] , \\
\end{split}
\end{equation}
where $v_{nk}(r) = J_n(\alpha_{nk} r/R) + c_{nk} Y_n(\alpha_{nk}
r/R)$.

For the special case of a disk ($R_0 = 0$), all the coefficients
$c_{nk}$ in front of the Bessel functions $Y_n(z)$ (divergent at $0$)
are set to $0$:
\begin{equation}
\label{eq:u_disk}
u_{nkl}(r,\varphi) = J_n(\alpha_{nk} r/R) 
\times \begin{cases} \cos(n\varphi), \quad l = 1,  \cr  \sin(n\varphi), \quad l = 2 ~ (n\ne 0),\end{cases}
\end{equation}
where $\alpha_{nk}$ are either the positive roots $j_{nk}$ of the
Bessel function $J_n(z)$ (Dirichlet), or the positive roots
$\tilde{j}_{nk}$ of its derivative $J'_n(z)$ (Neumann), or the
positive roots of their linear combination $J'_n(z) + h J_n(z)$
(Robin).  The asymptotic behavior of zeros of Bessel functions was
thoroughly investigated.  For fixed $k$ and large $n$, the Olver's
expansion holds $j_{nk} \simeq n + \delta_k n^{1/3} + O(n^{-1/3})$
(with known coefficients $\delta_k$) \cite{Olver51,Olver52,Elbert01},
while for fixed $n$ and large $k$, the McMahon's expansion holds:
$j_{nk} \simeq \pi(k + n/2 - 1/4) + O(k^{-1})$ \cite{Watson}.  Similar
asymptotic relations are applicable for Neumann and Robin boundary
conditions.

For a circular sector of radius $R$ and of angle $\pi \beta$, the
eigenfunctions are
\begin{equation}
u_{nk}(r,\varphi) = J_{n/\beta}(\alpha_{nk} r/R) \times \begin{cases} \sin(n \varphi/\beta)  \quad \textrm{(Dirichlet)}  \cr
\cos(n\varphi/\beta) \quad \textrm{(Neumann)}  \end{cases}  \quad (r < R,~ 0 < \varphi < \pi \beta)
\end{equation}
i.e., they are expressed in terms of Bessel functions of fractional
order, and $\alpha_{nk}$ are the positive roots of $J_{n/\beta}(z)$
(Dirichlet) or $J'_{n/\beta}(z)$ (Neumann).  The Robin boundary
condition and a sector of a circular annulus can be treated similarly.

\subsection{Sphere and spherical shell}

The rotation symmetry of a spherical shell in three dimensions,
$\Omega = \{ x\in\R^3 ~:~ R_0 < |x| < R\}$, allows one to write the
Laplace operator in spherical coordinates,
\begin{equation}
\begin{cases} x_1 = r \sin\theta \cos\varphi , \cr    x_2 = r \sin\theta \sin\varphi , \cr    x_3 = r \cos\theta  , \cr  \end{cases}  \qquad
\Delta = \frac{\partial^2}{\partial r^2} + \frac{2}{r} \frac{\partial}{\partial r} + \frac{1}{r^2} 
\left(\frac{1}{\sin\theta} \frac{\partial}{\partial\theta} \sin\theta \frac{\partial}{\partial \theta} +
\frac{\partial^2}{\partial\varphi^2}\right) ,
\end{equation}
that leads to the variable separation and an explicit representation
of eigenfunctions
\begin{equation}
u_{nkl}(r,\theta,\varphi) = \bigl[j_n(\alpha_{nk} r/R) + c_{nk} y_n(\alpha_{nk} r/R)\bigr] P_n^l(\cos\theta) e^{il\varphi},
\end{equation}
where $j_n(z)$ and $y_n(z)$ are the spherical Bessel functions of the
first and second kind, 
\begin{equation}
\label{eq:Bessels}
j_n(z) = \sqrt{\frac{\pi}{2z}} J_{n+1/2}(z),  \qquad   y_n(z) = \sqrt{\frac{\pi}{2z}} Y_{n+1/2}(z),
\end{equation}
$P_n^l(z)$ are associated Legendre polynomials (note that the
angular part, $P_n^l(\cos\theta) e^{il\varphi}$, is also called
spherical harmonic and denoted as $Y_{nl}(\theta,\varphi)$, up to a
normalization factor).  The coefficients $\alpha_{nk}$ and $c_{nk}$
are set by the boundary conditions at $r = R$ and $r = R_0$ similar to
Eq. (\ref{eq:BC_disk}).  The eigenfunctions are enumerated by the
triple index $nkl$, with $n = 0,1,2,...$ counting the order of
spherical Bessel functions, $k = 1,2,3,...$ counting zeros, and $l =
-n,-n+1,...,n$.  The eigenvalues $\lambda_{nk} =
\alpha_{nk}^2/R^2$, which are independent of the last index $l$, have
the degeneracy $2n+1$.  The squared $L_2$-norm of the eigenfunction is
\begin{equation}
\begin{split}
\| u_{nkl}(r,\theta,\varphi) \|_2^2 & = \frac{2\pi R^3}{(2n+1)\alpha_{nk}^2} \biggl[ \biggl(\alpha_{nk}^2 + h^2R^2 - hR - n(n+1)\biggr) v_{nk}^2(R) \\
& - \biggl(\alpha_{nk}^2 (R_0/R)^3 + h^2 R_0^2 -  h R_0 - n(n+1) R_0/R\biggr) v_{nk}^2(R_0)\biggr] , \\
\end{split}
\end{equation}
where $v_{nk}(r) = j_n(\alpha_{nk} r/R) + c_{nk} y_n(\alpha_{nk}
r/R)$.

In the special case of a sphere ($R_0 = 0$), one has $c_{nk} = 0$ and
the equations are simplified.  For balls and spherical shells in
higher dimensions ($d > 3$), the radial dependence of eigenfunctions
is expressed through a linear combination of so-called ultra-spherical
Bessel functions $r^{1-d/2} J_{\frac{d}{2}-1+n}(\alpha_{nk}r/R)$ and
$r^{1-d/2} Y_{\frac{d}{2}-1+n}(\alpha_{nk}r/R)$.

\subsection{Ellipse and elliptical annulus}
\label{sec:ellipse}

\begin{figure}
\begin{center}
\includegraphics[width=60mm]{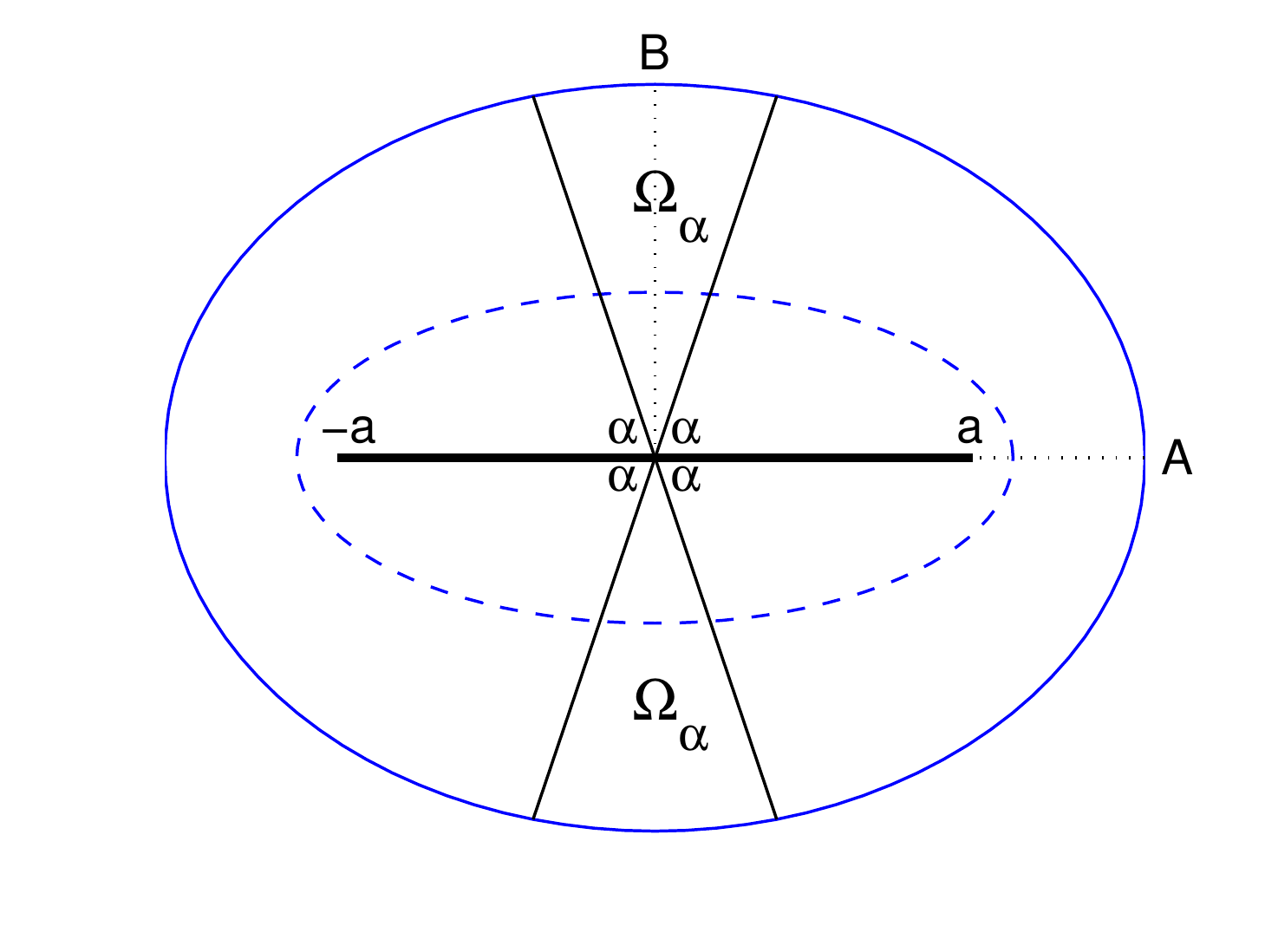}
\end{center}
\caption{
Two ellipses of radii $R = 0.5$ (dashed line) and $R = 1$ (solid
line), with the focal distance $a = 1$.  The major and minor
semi-axes, $A = a \cosh R$ and $B = a \sinh R$, are shown by black
dotted lines.  The horizontal thick segment connects the foci. }
\label{fig:ellipse}
\end{figure}

In elliptic coordinates, the Laplace operator reads as
\begin{equation}
\label{eq:elliptic}
\begin{cases}  x_1 = a \cosh r \cos \theta , \cr  x_2 = a \sinh r \sin \theta , \end{cases}    \hskip 10mm  
\Delta = \frac{1}{a^2(\sinh^2 r + \sin^2\theta)} \left(\frac{\partial^2}{\partial r^2} + \frac{\partial^2}{\partial\theta^2}\right),
\end{equation} 
where $a > 0$ is the prescribed distance between the origin and the
foci, $r \geq 0$ and $0\le \theta < 2\pi$ are the radial and angular
coordinates (Fig. \ref{fig:ellipse}).  An ellipse is a curve of
constant $r = R$ so that its points $(x_1,x_2)$ satisfy $x_1^2/A^2 +
x_2^2/B^2 = 1$, where $R$ is the ``radius'' of the ellipse and $A = a
\cosh R$ and $B = a \sinh R$ are the major and minor semi-axes.  Note
that the eccentricity $e = a/A = 1/\cosh R$ is strictly positive.  A
filled ellipse (i.e., the interior of an given ellipse) can be
characterized in elliptic coordinates as $0 \leq r < R$ and $0\le
\theta < 2\pi$.  Similarly, an elliptical annulus (i.e., the interior
between two ellipses with the same foci) is characterized by $R_0 < r
< R$ and $0\le \theta < 2\pi$.

In the elliptic coordinates, the variables can be separated,
$u(r,\theta) = f(r) g(\theta)$, from which Eq. (\ref{eq:eigen}) reads
as
\begin{equation*}
\left(\frac{1}{f(r)} \frac{d^2f}{dr^2} + \frac{\lambda a^2}{2} \cosh(2r)\right) =
- \left(\frac{1}{g(\theta)} \frac{d^2g}{d\theta^2} - \frac{\lambda a^2}{2} \cos(2\theta)\right) 
\end{equation*}
so that both sides are equal to a constant (denoted $c$).  As a
consequence, the angular and radial parts, $g(\theta)$ and $f(r)$, are
solutions of the Mathieu equation and the modified Mathieu equation,
respectively \cite{McLachlan,Zhang,Chen94}
\begin{equation*}
g''(\theta) + \left(c - 2q \cos 2\theta\right) g(\theta) = 0, \qquad
f''(r) - \left(c - 2q \cosh 2r\right) f(r) = 0,
\end{equation*}
where $q = \lambda a^2/4$ and the parameter $c$ is called the
characteristic value of Mathieu functions.  Periodic solutions of the
Mathieu equation are possible for specific values of $c$.  They are
denoted as $\ce_n(\theta,q)$ and $\se_{n+1}(\theta,q)$ (with $n =
0,1,2,...$) and called the angular Mathieu functions of the first and
second kind.  Each function $\ce_n(\theta,q)$ and
$\se_{n+1}(\theta,q)$ corresponds to its own characteristic value $c$
(the relation being implicit, see \cite{McLachlan}).

For the radial part, there are two linearly independent solutions for
each characteristic value $c$: two modified Mathieu functions
$\Mc_n^{(1)}(r,q)$ and $\Mc_n^{(2)}(r,q)$ correspond to the same $c$
as $\ce_n(\theta,q)$, and two modified Mathieu functions
$\Ms_{n+1}^{(1)}(r,q)$ and $\Ms_{n+1}^{(2)}(r,q)$ correspond to the
same $c$ as $\se_{n+1}(\theta,q)$.  As a consequence, there are four
families of eigenfunctions (distinguished by the index $l=1,2,3,4$) in
an elliptical domain
\begin{eqnarray*}
u_{nk1}(r,\theta) &=& \ce_n(\theta,q_{nk1}) \Mc_n^{(1)}(r,q_{nk1}) , \\
u_{nk2}(r,\theta) &=& \ce_n(\theta,q_{nk2}) \Mc_n^{(2)}(r,q_{nk2}) , \\
u_{nk3}(r,\theta) &=& \se_{n+1}(\theta,q_{nk3}) \Ms_{n+1}^{(1)}(r,q_{nk3}) , \\ 
u_{nk4}(r,\theta) &=& \se_{n+1}(\theta,q_{nk4}) \Ms_{n+1}^{(2)}(r,q_{nk4}) ,  
\end{eqnarray*}
where the parameters $q_{nkl}$ are determined by the boundary
condition.  For instance, for a filled ellipse of radius $R$ with
Dirichlet boundary condition, there are four individual equations for
the parameter $q$ for each $n = 0,1,2,...$
\begin{equation*}
\Mc_n^{(1)}(R,q_{nk1}) = 0 ,  \hskip 2mm   \Mc_n^{(2)}(R,q_{nk2}) = 0 , \hskip 2mm
\Ms_{n+1}^{(1)}(R,q_{nk3}) = 0 ,  \hskip 2mm   \Ms_{n+1}^{(2)}(R,q_{nk4}) = 0 ,
\end{equation*} 
each of them having infinitely many positive solutions $q_{nkl}$
enumerated by $k=1,2,\dots$ \cite{McLachlan,Abramowitz}.
Finally, the associated eigenvalues are $\lambda_{nkl} = 4q_{nkl}/a^2$.

\subsection{Equilateral triangle}

Lam\'e discovered the Dirichlet eigenvalues and eigenfunctions of the
equilateral triangle $\Omega = \{ (x_1,x_2)\in\R^2~: 0 < x_1 < 1,~ 0 <
x_2 < x_1\sqrt{3},~ x_2 < \sqrt{3}(1-x_1)\}$ by using reflections and
the related symmetries \cite{Lame}:
\begin{equation}
\lambda_{mn} = \frac{16\pi^2}{27}\bigl(m^2 + n^2 - mn\bigr) \qquad (m,n\in \Z),
\end{equation}
where $3$ divides $m+n$, $m\ne 2n$, and $n\ne 2m$, and the associate
eigenfunction is
\begin{equation}
u_{mn}(x_1,x_2) = \sum\limits_{(m',n')} \pm \exp\left[\frac{2\pi i}{3}\biggl(m'x_1 + (2n'-m')\frac{x_2}{\sqrt{3}}\biggr)\right],
\end{equation}
where $(m',n')$ runs over $(-n,m-n)$, $(-n,-m)$, $(n-m,-m)$,
$(n-m,n)$, $(m,n)$ and $(m,m-n)$ with the $\pm$ sign alternating
(see also \cite{Mathews} for basic introduction, as well as
\cite{Makai70,Lee78}).  Pinsky showed that this set of eigenfunctions
is complete in $L_2(\Omega)$ \cite{Pinsky80,Pinsky85}.  Note that the
conditions $m\ne 2n$ and $n\ne 2m$ should be satisfied for all 6 pairs
in the sum that yields one additional condition: $m\ne -n$.  The
following relations hold: $u_{-m,-n} = u_{m,n}^*$, $u_{n,m} = -
u_{m,n}^*$ and $u_{m,0} = u_{m,m}$.  All symmetric eigenfunctions are
enumerated by the index $(m,0)$.  The eigenvalue $\lambda_{mn}$
corresponds to a symmetric eigenfunction if and only if $m$ is a
multiple of $3$ \cite{Pinsky80}.

The eigenfunctions for Neumann boundary condition are
\begin{equation}
u_{mn}(x_1,x_2) = \sum\limits_{(m',n')} \exp\left[\frac{2\pi i}{3}\biggl(m'x_1 + (2n'-m')\frac{x_2}{\sqrt{3}}\biggr)\right] ,
\end{equation}
where the only condition is that $m+n$ are multiples of $3$ (and no
sign change).  Further references and extensions (e.g., to Robin
boundary conditions) can be found in a series of works by McCartin
\cite{McCartin03,McCartin02,McCartin05,McCartin07,McCartin11}.  
McCartin also developed a classification of all polygonal domains
possessing a complete set of trigonometric eigenfunctions of the
Laplace operator under either Dirichlet or Neumann boundary conditions
\cite{McCartin08}.

\section{Eigenvalues}
\label{sec:eigenvalues}

\subsection{Weyl's law}

The Weyl's law is one of the first connections between the spectral
properties of the Laplace operator and the geometrical structure of a
bounded domain $\Omega$.  In 1911, Hermann Weyl derived the asymptotic
behavior of the Laplacian eigenvalues \cite{Weyl1911,Weyl1912}:
\begin{equation}
\lambda_m \propto \frac{4\pi^2}{(\omega_d \mu_d(\Omega))^{2/d}}~ m^{2/d}  \qquad (m\to\infty),
\end{equation}
where $\mu_d(\Omega)$ is the Lebesgue measure of $\Omega$ (its area in
2D and volume in 3D), and
\begin{equation}
\label{eq:omega_d}
\omega_d = \frac{\pi^{d/2}}{\Gamma(d/2+1)}
\end{equation}
is the volume of the unit ball in $d$ dimensions ($\Gamma(z)$ being
the Gamma function).  As a consequence,
plotting eigenvalues versus $m^{2/d}$ allows one to extract the area
in 2D or the volume in 3D.  This result can equivalently be written
for the counting function $N(\lambda) = \#\{m~:~ \lambda_m <
\lambda\}$ (i.e., the number of eigenvalues smaller than $\lambda$):
\begin{equation}
N(\lambda) \propto \frac{\omega_d \mu_d(\Omega)}{(2\pi)^d} ~
\lambda^{d/2} \qquad (\lambda \to \infty).
\end{equation}

Weyl also conjectured the second asymptotic term which in two and
three dimensions reads as
\begin{equation}
N(\lambda) \propto \begin{cases} \displaystyle \frac{\mu_2(\Omega)}{4\pi} ~ \lambda \mp \frac{\mu_1(\pa)}{4\pi} \sqrt{\lambda}  \qquad (d=2)  \cr
\displaystyle \frac{\mu_3(\Omega)}{6\pi^2} ~ \lambda^{3/2} \mp \frac{\mu_2(\pa)}{16\pi} \lambda  \qquad (d=3) \end{cases}   \qquad (\lambda \to \infty),
\end{equation}
where $\mu_2(\Omega)$ and $\mu_1(\pa)$ are the area and
perimeter of $\Omega$ in 2D, $\mu_3(\Omega)$ and $\mu_2(\pa)$ are the
volume and surface area of $\Omega$ in 3D, and sign ``--'' (resp.\
``+'') refers to the Dirichlet (resp.\ Neumann) boundary condition.
The correction terms which yield information about the boundary of the
domain, were justified, under certain conditions on $\Omega$ (e.g.,
convexity) only in 1980 \cite{Ivrii80,Melrose80} (see \cite{Arendt09}
for a historical review and further details).

Alternatively, one can study the heat trace (or partition function)
\begin{equation}
Z(t) \equiv \int\limits_\Omega dx ~ G_t(x,x) = \sum\limits_{m=1}^\infty e^{-\lambda_m t} = 
\int\limits_0^\infty e^{-\lambda t} dN(\lambda) 
\end{equation}
(here $G_t(x,y)$ is the heat kernel, cf. Eq. (\ref{eq:heat})), for
which the following asymptotic expansion holds
\cite{Minakshisundaram53,McKean67,Protter87,Branson90,Davies,Davies2,Gilkey}
\begin{equation}
Z(t) = (4\pi t)^{-d/2} \left(\sum\limits_{k=0}^K c_k t^{k/2} + o(t^{(K+1)/2})\right) \qquad (t\to 0),
\end{equation}
where the coefficients $c_k$ are again related to the geometrical
characteristics of the domain:
\begin{equation}
c_0 = \mu_d(\Omega), \quad c_1 = - \frac{\sqrt{\pi}}{2} \mu_{d-1}(\pa), \quad ...
\end{equation}
(see \cite{Reuter06} for further discussion).  Some estimates for the
trace of the Dirichlet Laplacian were given by Davies \cite{Davies84}
(see also \cite{vandenBerg99b} for the asymptotic behavior of the heat
content).

A number of extensions have been proposed.  Berry conjectured that,
for irregular boundaries, for which the Lebesgue measure in the
correction term is infinite, the correction term should be
$\lambda^{H/2}$ instead of $\lambda^{(d-1)/2}$, where $H$ is the
Hausdorff dimension of the boundary \cite{Berry79,Berry80}.  However,
Brossard and Carmona constructed a counter-example to this conjecture
and suggested a modified version, in which the Hausdorff dimension was
replaced by Minkowski dimension \cite{Brossard86}.  The modified
Weyl-Berry conjecture was discussed at length by Lapidus who proved it
for $d = 1$ \cite{Lapidus91,Lapidus93} (see these references for
further discussion).  For dimension $d$ higher than $1$, this
conjecture was disproved by Lapidus and Pomerance \cite{Lapidus96}.
The correction term to the Weyl's formula for domains with rough
boundary (in particular, for Lipschitz class) was studied by Netrusov
and Safarov \cite{Netrusov05}.  Levitin and Vassiliev also considered
the asymptotic formulas for iterated sets with fractal boundary
\cite{Levitin96}.  Extensions to various manifolds and higher order
Laplacians were discussed
\cite{Desjardins94,Desjardins98}.

The high-frequency Weyl's law and the related short-time asymptotics
of the heat kernel have been thoroughly investigated \cite{Arendt09}.
The dependence of these asymptotic laws on the volume and surface of
the domain has found applications in physics.  For instance,
diffusion-weighted nuclear magnetic resonance experiments were
proposed and conducted to estimate the surface-to-volume ratio of
mineral samples and biological tissues
\cite{Mitra92,Mitra93,Latour93,Hurlimann94,Latour94,Helmer95,Sen04,Grebenkov07}.

The multiplicity of eigenvalues is yet a more difficult problem
\cite{Nadirashvili88b}.  From basic properties (see
Sec. \ref{sec:general}), the first eigenvalue $\lambda_1$ is simple.
Cheng proved that the multiplicity $m(\lambda_2)$ of the second
Dirichlet eigenvalue $\lambda_2$ is not greater than $3$
\cite{Cheng76}.  This inequality is sharp since an example of domain
with $m(\lambda_2) = 3$ was constructed.  For $k \geq 3$,
Hoffmann-Ostenhof {\it et al.} proved the inequality $m(\lambda_k)
\leq 2k-3$ \cite{Hoffmann99a,Hoffmann99b}.

\subsection{Isoperimetric inequalities for eigenvalues}

In the low-frequency limit, the relation between the shape of a domain
and the associated eigenvalues manifests in the form of isoperimetric
inequalities.  Since there are many excellent reviews on this topic,
we only provide a list of the best-known inequalities, while further
discussion and references can be found in
\cite{Polya,Payne67,Bandle,Hile80,Kuttler84,Hansen94,Ashbaugh00,Reuter06,Henrot,Ashbaugh07,Benguria11}.

(i) {\it The Rayleigh-Faber-Krahn inequality} states that the disk
minimizes the first Dirichlet eigenvalue $\lambda_1$ among all planar
domains of the same area $\mu_2(\Omega)$, i.e.
\begin{equation}
\label{eq:RFK1}
\lambda_1^D \geq \frac{\pi}{\mu_2(\Omega)} (j_{0,1})^2, 
\end{equation}
where $j_{\nu,1}$ is the first positive zero of $J_\nu(z)$ (e.g.,
$j_{0,1}\approx 2.4048...$).  This inequality was conjectured by Lord
Rayleigh and proven independently by Faber and Krahn
\cite{Faber23,Krahn25}.  The corresponding isoperimetric inequality in
$d$ dimensions,
\begin{equation}
\label{eq:RFK1d}
\lambda_1^D \geq \left(\frac{\omega_d}{\mu_d(\Omega)}\right)^{2/d} (j_{\frac{d}{2}-1,1})^2,
\end{equation}
was proven by Krahn \cite{Krahn26}.

Another lower bound for the first Dirichlet eigenvalue for a simply
connected planar domain was obtained by Makai \cite{Makai65} and later
rediscovered (in a weaker form) by Hayman \cite{Hayman78}
\begin{equation}
\label{eq:lambda1D_L}
\lambda_1^D \geq \frac{\alpha}{\rho^2}
\end{equation}
where $\alpha$ is a constant, and
\begin{equation}
\label{eq:inradius}
\rho = \max\limits_{x \in \Omega} \min\limits_{y\in \pa} \{ |x-y|\}
\end{equation}
is the inradius of $\Omega$ (i.e., the radius of the largest ball
inscribed in $\Omega$).  The above inequality means that the lowest
frequency (bass note) can be made arbitrarily small only if the domain
includes an arbitrarily large circular drum (i.e., $\rho$ goes to
infinity).  The constant $\alpha$, which was equal to $1/4$ in the
original Makai's proof (see also \cite{Osserman77}) and to $1/900$ in
the Hayman's proof, was gradually increased, to the best value (up to
date) $\alpha = 0.6197...$ by Banuelos and Carroll \cite{Banuelos94}.
For convex domains, the lower bound (\ref{eq:lambda1D_L}) with $\alpha
= \pi^2/4 \approx 2.4674$ was derived much earlier by Hersch
\cite{Hersch60}, with the equality if and only if $\Omega$ is an
infinite strip (see also a historical overview in
\cite{Ashbaugh07}).

An obvious upper bound for the first Dirichlet eigenvalue can be
obtained from the domain monotonicity (property (v) in
Sec. \ref{sec:general}):
\begin{equation}
\lambda_1^D \leq \lambda_1^D(B_\rho) =  \rho^{-2} ~ j_{\frac{d}{2}-1,1}^2 ,
\end{equation}
with the first Dirichlet eigenvalue $\lambda_1^D(B_\rho)$ for the
largest ball $B_{\rho}$ inscribed in $\Omega$ ($\rho$ is the
inradius).  However, this upper bound is not accurate in general.
P\'olya and Szeg\H{o} gave another upper bound for planar star-shaped
domains \cite{Polya}.  Freitas and Krej\v{c}i\v{r}\'ik extended their
result to higher dimensions \cite{Freitas08}: for a bounded strictly
star-shaped domain $\Omega\subset \R^d$ with locally Lipschitz
boundary, they proved
\begin{equation}
\lambda_1^D \leq \lambda_1^D(B_1) \frac{F(\Omega)}{d~ \mu_d(\Omega)} ,
\end{equation}
where the function $F(\Omega)$ is defined in \cite{Freitas08}. From
this inequality, they also deduced a weaker but more explicit upper
bound which is applicable to any bounded convex domain in $\R^d$:
\begin{equation}
\lambda_1^D \leq \lambda_1^D(B_1) ~ \frac{\mu_{d-1}(\pa)}{d~\rho ~\mu_d(\Omega)} .
\end{equation}

The second Dirichlet eigenvalue $\lambda_2^D$ is minimized by the
union of two identical balls 
\begin{equation}
\lambda_2^D \geq 2^{2/d} \left(\frac{\omega_d}{\mu_d(\Omega)}\right)^{2/d} (j_{\frac{d}{2}-1,1})^2 .
\end{equation}
This inequality, which can be deduced by looking at nodal domains for
$u_2$ and using Rayleigh-Faber-Krahn inequality (\ref{eq:RFK1d}) on
each nodal domain, was first established by Krahn \cite{Krahn26}.  It
is also sometimes attributed to Peter Szeg\H{o} (see \cite{Polya55}).
Note that finding the minimizer of $\lambda_2^D$ among convex planar
sets is still an open problem \cite{Henrot03}.  Bucur and Henrot
proved the existence of a minimizer for the third eigenvalue in the
family of domains in $\R^d$ of given volume, although its shape
remains unknown \cite{Bucur00}.  The range of the first two
eigenvalues was also investigated
\cite{Wolf94,Bucur99}.

The first nontrivial Neumann eigenvalue $\lambda_2^N$ (as $\lambda_1^N
= 0$) also satisfies the isoperimetric inequality
\begin{equation}
\label{eq:RFK_N}
\lambda_2^N \leq \left(\frac{\omega_d}{\mu_d(\Omega)}\right)^{2/d} (\tilde{j}_{\frac{d}{2},1})^2 ,
\end{equation}
which states that $\lambda_2^N$ is maximized by a $d$-dimensional ball
(here $\tilde{j}_{\nu,1}$ is the first positive zero of the function
$\frac{d}{dz}[z^{1-d/2} J_{\frac{d}{2}-1+\nu}(z)]$ which reduces to
$J'_{\nu}(z)$ and $\sqrt{2/\pi} ~ j'_{\nu}(z)$ for $d = 2$ and $d=3$,
respectively).  This inequality was proven for simply-connected planar
domains by Szeg\H{o} \cite{Szego54} and in higher dimensions by
Weinberger \cite{Weinberger56}.  P\'olya conjectured the following
upper bound for all Neumann eigenvalues \cite{Polya54} in planar
bounded regular domains (see also \cite{Schoen})
\begin{equation}
\label{eq:Polya_N}
\lambda_n^N \leq \frac{4(n-1)\pi}{\mu_2(\Omega)}    \qquad (n = 2,3,4,...) 
\end{equation}
(the domain is called regular if its Neumann eigenspectrum is
discrete, see \cite{Girouard09} for details).  This inequality is true
for all domains that tile the plane, e.g., for any triangle and any
quadrilateral \cite{Polya61}.  For $n = 2$, the inequality
(\ref{eq:Polya_N}) follows from (\ref{eq:RFK_N}).  For $n
\geq 3$, P\'olya's conjecture is still open, although Kr\"oger
proved a weaker estimate $\lambda_n^N \leq 8\pi (n-1)$
\cite{Kroger92}.  Recently, Girouard {\it et al.} obtained a sharp
upper bound for the second nontrivial Neumann eigenvalue $\lambda_3^N$
for a regular simply-connected planar domain
\cite{Girouard09}:
\begin{equation}
\lambda_3^N \leq \frac{2\pi (\tilde{j}_{0,1})^2}{\mu_2(\Omega)} ,
\end{equation}
with the equality attained in the limit by a family of domains
degenerating to a disjoint union of two identical disks.

Payne and Weinberger obtained the lower bound for the second Neumann
eigenvalue in $d$ dimensions \cite{Payne60}
\begin{equation}
\lambda_2^N \geq \frac{\pi^2}{\delta^2} ,
\end{equation}
where $\delta$ is the diameter of $\Omega$:
\begin{equation}
\label{eq:diameter}
\delta = \max\limits_{x,y \in \pa} \{|x-y|\} .
\end{equation}
This is the best bound that can be given in terms of the diameter
alone in the sense that $\lambda_2^N \delta^2$ tends to $\pi^2$ for a
parallelepiped all but one of whose dimensions shrink to zero.

Szeg\H{o} and Weinberger noticed that Szeg\H{o}'s proof of the
inequality (\ref{eq:RFK_N}) for planar simply connected domains
extends to prove the bound
\begin{equation}
\frac{1}{\lambda_2^N}+\frac{1}{\lambda_3^N} \ge \frac{2\mu_2(\Omega)}{\pi (\tilde{j}_{1,1})^2},
\end{equation}
with equality if and only if $\Omega$ is a disk
\cite{Szego54,Weinberger56}.  Ashbaugh and Benguria derived another
bound for arbitrary bounded domain in $\R^d$ \cite{Ashbaugh93a}
\begin{equation}
\frac{1}{\lambda_2^N} + ... + \frac{1}{\lambda_{d+1}^N} \ge \frac{d}{d+2} \left(\frac{\mu_d(\Omega)}{\omega_d}\right)^{2/d} 
\end{equation}
In particular, one gets $1/\lambda_2^N + 1/\lambda_3^N \ge
\frac{\mu_2(\Omega)}{2\pi}$ for $d = 2$ (see also extensions in
\cite{Hile93,Xia99}).

(ii) {\it The Payne-P\'olya-Weinberger inequality}, which can also be
called {\it Ashbaugh-Benguria inequality}, concerns the ratio between
first two Dirichlet eigenvalues and states that
\begin{equation}
\frac{\lambda_2^D}{\lambda_1^D} \leq \left(\frac{j_{\frac{d}{2},1}}{j_{\frac{d}{2}-1,1}}\right)^2 ,
\end{equation}
with equality if and only if $\Omega$ is the $d$-dimensional ball.
This inequality (in 2D form) was conjectured by Payne, P\'olya and
Weinberger \cite{Payne56} and proved by Ashbaugh and Benguria in 1990
\cite{Ashbaugh91,Ashbaugh92,Ashbaugh93a,Ashbaugh93b}.  A weaker estimate
$\lambda_2^D/\lambda_1^D \leq 1 + 4/d$ was proved for $d = 2$ in the
original paper by Payne, P\'olya and Weinberger \cite{Payne56}.

(iii) Singer {\it et al.} derived the upper and lower estimates for
the spectral (or fundamental) gap between the first two Dirichlet
eigenvalues for a smooth convex bounded domain $\Omega$ in $\R^d$ (in
\cite{Singer85}, a more general problem in the presence of a potential
was considered):
\begin{equation}
\frac{d\pi^2}{\rho^2} \ge \lambda_2^D - \lambda_1^D \ge \frac{\pi^2}{4\delta^2} ,
\end{equation}
where $\delta$ is the diameter of $\Omega$ and $\rho$ is the inradius
\cite{Singer85}.  For a convex planar domain, Donnelly proposed a
sharper lower estimate \cite{Donnelly11}
\begin{equation}
\lambda_2^D - \lambda_1^D \geq \frac{3\pi^2}{\delta^2} .
\end{equation}
However, Ashbaugh {\it et al.} pointed out on a flaw in the proof
\cite{Ashbaugh11}.  The estimate was later rigorously proved by
Andrews and Clutterbuck for any bounded convex domain $\Omega$ in
$\R^d$, even in the presence of a semi-convex potential
\cite{Andrews11} (for more background on the spectral gap, see
notes by Ashbaugh \cite{Ashbaugh06}).

(iv) The isoperimetric inequalities for Robin eigenvalues are less
known.  Daners proved that among all bounded domains $\Omega\subset
\R^d$ of the same volume, the ball $B$ minimizes the first Robin
eigenvalue \cite{Daners06,Bucur10}
\begin{equation}
\lambda_1^R(\Omega) \geq \lambda_1^R(B). 
\end{equation}
Kennedy showed that among all bounded domains in $\R^d$, a domain
$B_2$ composed of two disjoint balls minimizes the second Robin
eigenvalue \cite{Kennedy09}
\begin{equation}
\lambda_2^R(\Omega) \geq \lambda_2^R(B_2).
\end{equation}

(v) The minimax principle ensures that the Neumann eigenvalues are
always smaller than the corresponding Dirichlet eigenvalues:
$\lambda^N_n \leq \lambda^D_n$.  P\'olya proved $\lambda^N_2 <
\lambda^D_1$ \cite{Polya52} while Szeg\H{o} got a sharper inequality
$\lambda^N_2 \leq c\lambda^D_1$ for a planar domain bounded by an
analytic curve, where $c = (\tilde{j}_{1,1}/j_{0,1})^2 \approx
0.5862...$ \cite{Szego54} (note that this result also follows from
inequalities (\ref{eq:RFK1}, \ref{eq:RFK_N})).  Payne derived a
stronger inequality for a planar domain with a $C^2$ boundary:
$\lambda^N_{n+2} < \lambda^D_n$ for all $n$ \cite{Payne55}.  Levine
and Weinberger generalized this result for higher dimensions $d$ and
proved that $\lambda^N_{n+d} < \lambda^D_n$ for all $n$ when $\Omega$
is smooth and convex, and that $\lambda^N_{n+d} \leq \lambda^D_n$ if
$\Omega$ is merely convex \cite{Levine86}.  Friedlander proved the
inequality $\lambda^N_{n+1} \leq \lambda^D_n$ for a general bounded
domain with a $C^1$ boundary \cite{Friedlander91}.  Filonov found a
simpler proof of this inequality in a more general situation (see
\cite{Filonov05} for details).  

Many other inequalities can be found in several reviews
\cite{Ashbaugh00,Ashbaugh07,Benguria11}.  It is worth noting that
isoperimetric inequalities are related to shape optimization problems
\cite{Buttazzo91,Pironneau84,Sokolowski92,Allaire02,Bucur05,Bourdin09,Buttazzo11}.

\subsection{Kac's inverse spectral problem}

\begin{figure}
\begin{center}
\includegraphics[width=120mm]{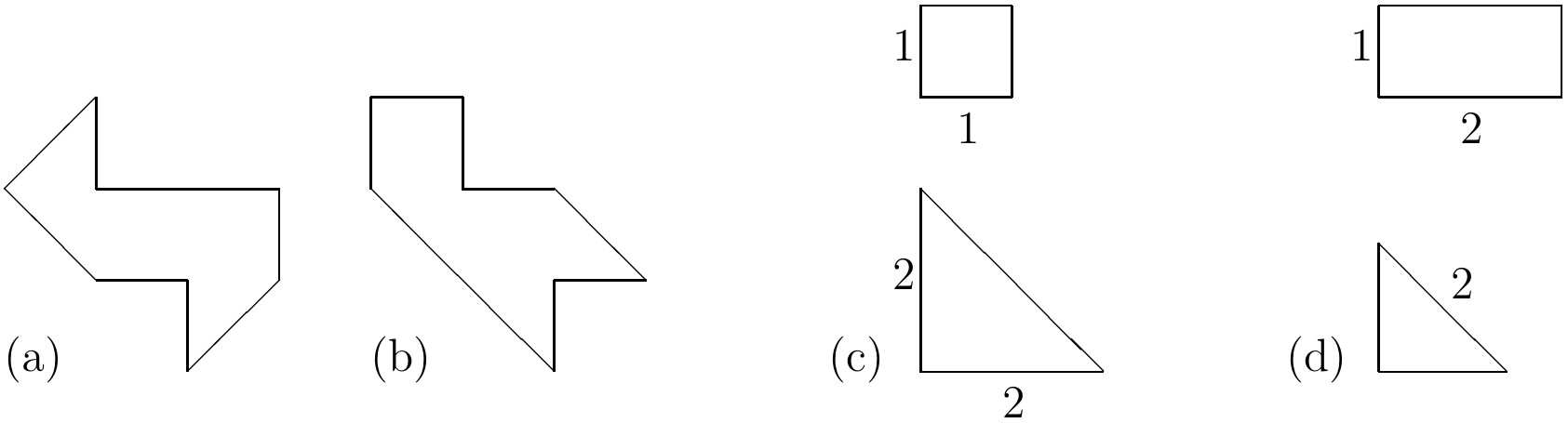}
\end{center}
\caption{ 
Two examples of nonisometric domains with the identical Laplace
operator eigenspectra (with Dirichlet or Neumann boundary conditions):
the original example (shapes 'a' and 'b') constructed by Gordon {\it
et al.} \cite{Gordon92}, and a simpler example with disconnected
domains (shapes 'c' and 'd') by Chapman \cite{Chapman95}.  In the
latter case, the eigenspectrum is simply obtained as the union of the
eigenspectra of two subdomains known explicitly.  For instance, the
Dirichlet eigenspectrum is $\{\pi^2(m^2+n^2)~:~ m,n\in \N\} \cup \{
\pi^2((i/2)^2+(j/2)^2)~:~ i,j\in \N,~ i>j\}$ .}
\label{fig:gordon}
\end{figure}

The problem of finding relations between the Laplacian eigenspectrum
and the shape of a domain was formulated in the famous Kac's question
``Can one hear the shape of a drum?'' \cite{Kac66}.  In fact, the
drum's frequencies are uniquely determined by the eigenvalues of the
Laplace operator in the domain of drum's shape.  By definition, the
shape of the domain fully determines the Laplacian eigenspectrum.  Is
the opposite true, i.e., does the set of eigenvalues which appear as
``fingerprints'' of the shape, uniquely identify the domain?  The
negative answer to this question for general planar domains was given
by Gordon and co-authors \cite{Gordon92} who constructed two different
(nonisometric) planar polygons (Fig.~\ref{fig:gordon}a,b) with the
identical Laplacian eigenspectra, both for Dirichlet and Neumann
boundary conditions (see also \cite{Berard92}).  Their construction
was based on Sunada's paper on isospectral manifolds \cite{Sunada85}.
An elementary proof, as well as many other examples of isospectral
domains, were provided by Buser and co-workers \cite{Buser94} and by
Chapman \cite{Chapman95} (see Fig. \ref{fig:gordon}c,d).  An
experimental evidence for this not ``hearing the shape'' of drums was
brought by Sridhar and Kudrolli \cite{Sridhar94} (see also
\cite{Cipra92}).  In all these examples, isospectral domains are
either non-convex, or disjoint.  Gordon and Webb addressed the
question of existence of isospectral convex connected domains and
answered this question positively (i.e., negatively to the original
Kac's question) for domains in Euclidean spaces of dimension $d \geq
4$ \cite{Gordon94}.  To our knowledge, this question remains open for
convex domains in two and three dimensions, as well as for domains
with smooth boundaries.  It is worth noting that the positive answer
to Kac's question can be given for some classes of domains.  For
instance, Zelditch proved that for domains that possess the symmetry
of an ellipse and satisfy some generic conditions on the boundary, the
spectrum of the Dirichlet Laplacian uniquely determines the shape
\cite{Zelditch00}.  Later, he extended this result to real analytic
planar domains with only one symmetry \cite{Zelditch04,Zelditch09}.

A somewhat similar problem was recently formulated for domains in
which one part of the boundary admits Dirichlet boundary condition and
the other Neumann boundary condition.  Does the spectrum of the
Laplace operator determine uniquely which condition is imposed on
which part?  Jakobsen and co-workers gave the negative answer to this
question by assigning Dirichlet and Neumann conditions onto different
parts of the boundary of the half-disk (and some other domains), in a
way to produce the same eigenspectra \cite{Jakobson06}.

The Kac's inverse spectral problem can also be seen from a different
point of view.  For a given sequence $0 \leq \lambda_1 < \lambda_2
\leq \lambda_3 \leq ...$, whether does exist a domain $\Omega$ in $\R^d$
for which the Laplace operator with Dirichlet or Neumann boundary
condition has the spectrum given by this sequence.  A similar problem
can be formulated for a compact Riemannian manifold with arbitrary
Riemannian metrics.  Colin de Verdi\`ere studied these problems for
finite sequences $\{\lambda_n\}_{n=1}^N$ and proved the existence of
such domains or manifolds under certain restrictions
\cite{ColindeVerdiere87}.  We also mention the work of Sleeman who
discusses the relationship between the inverse scattering theory
(i.e., the Helmholtz equation for an exterior domain) and the Kac's
inverse spectral problem (i.e., for an interior domain)
\cite{Sleeman09} (see \cite{Chu} for further discussion on inverse
eigenvalue problems).

\section{Nodal lines}
\label{sec:nodal}

The first insight onto the geometrical structure of eigenfunctions can
be gained from their nodal lines.  Kuttler and Sigillito gave a brief
overview of the basic properties of nodal lines for Dirichlet
eigenfunctions in two dimensions \cite{Kuttler84} that we partly
reproduce here:

``The set of points in $\Omega$ where $u_m = 0$ is the nodal set of
$u_m$.  By the unique continuation property, it consists of curves
that are $C^\infty$ in the interior of $\Omega$.  Where nodal lines
cross, they form equal angles \cite{Courant}.  Also, when nodal lines
intersect a $C^\infty$ portion of the boundary, they form equal
angles.  Thus, a single nodal line intersects the $C^\infty$ boundary
at right angles, two intersect it at $60^\circ$ angles, and so forth.
Courant's nodal line theorem \cite{Courant} states that the nodal
lines of the $m$-th eigenfunction divide $\Omega$ into no more than
$m$ subregions (called nodal domains): $\nu_m \leq m$, $\nu_m$ being
the number of nodal domains.  In particular, $u_1$ has no interior
nodes and so $\lambda_1$ is a simple eigenvalue (has multiplicity
one).''

It is worth noting that any eigenvalue $\lambda_m$ of the
Dirichlet-Laplace operator in $\Omega$ is the first eigenvalue for
each of its nodal domains.  This simple observation allows one to
construct specific domains with a prescribed eigenvalue (see
\cite{Kuttler84} for examples).  Eigenfunctions with few nodal domains
were constructed in \cite{Courant,Lewy77}.

\begin{figure}
\begin{center}
\includegraphics[width=62mm]{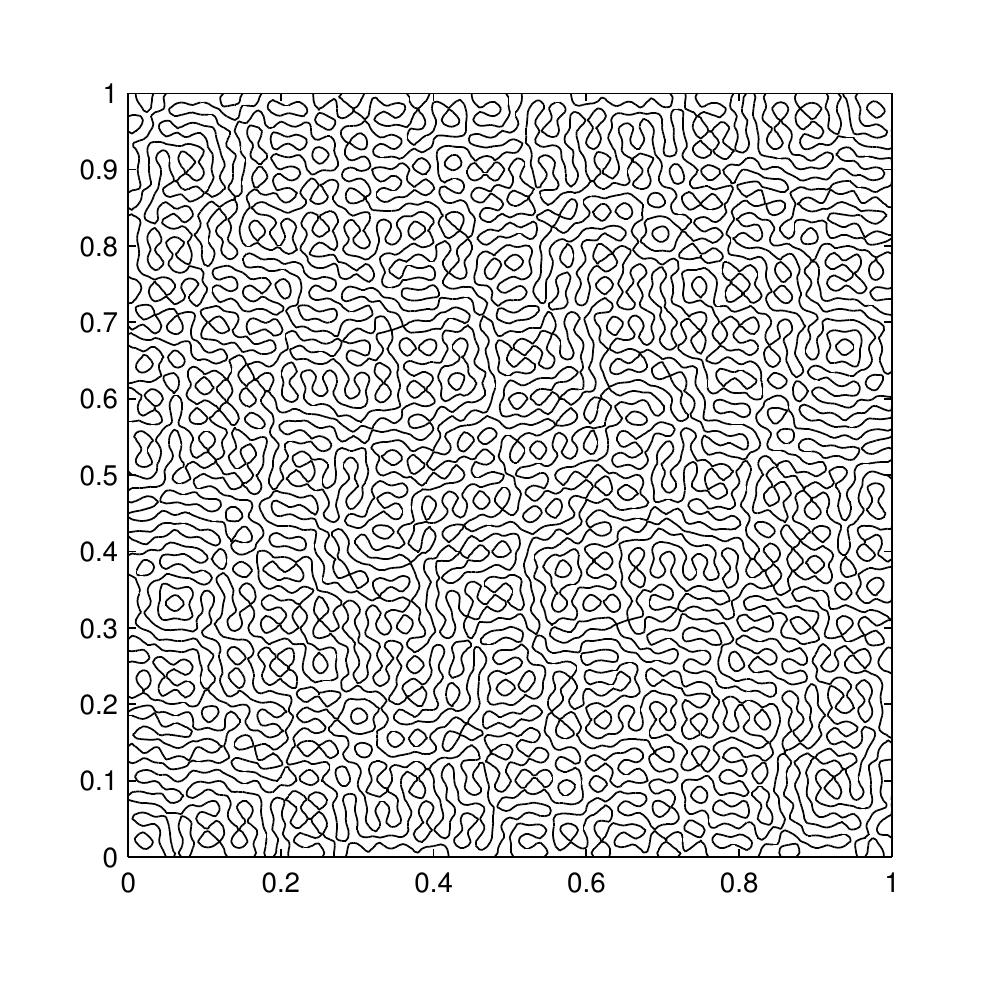}
\includegraphics[width=62mm]{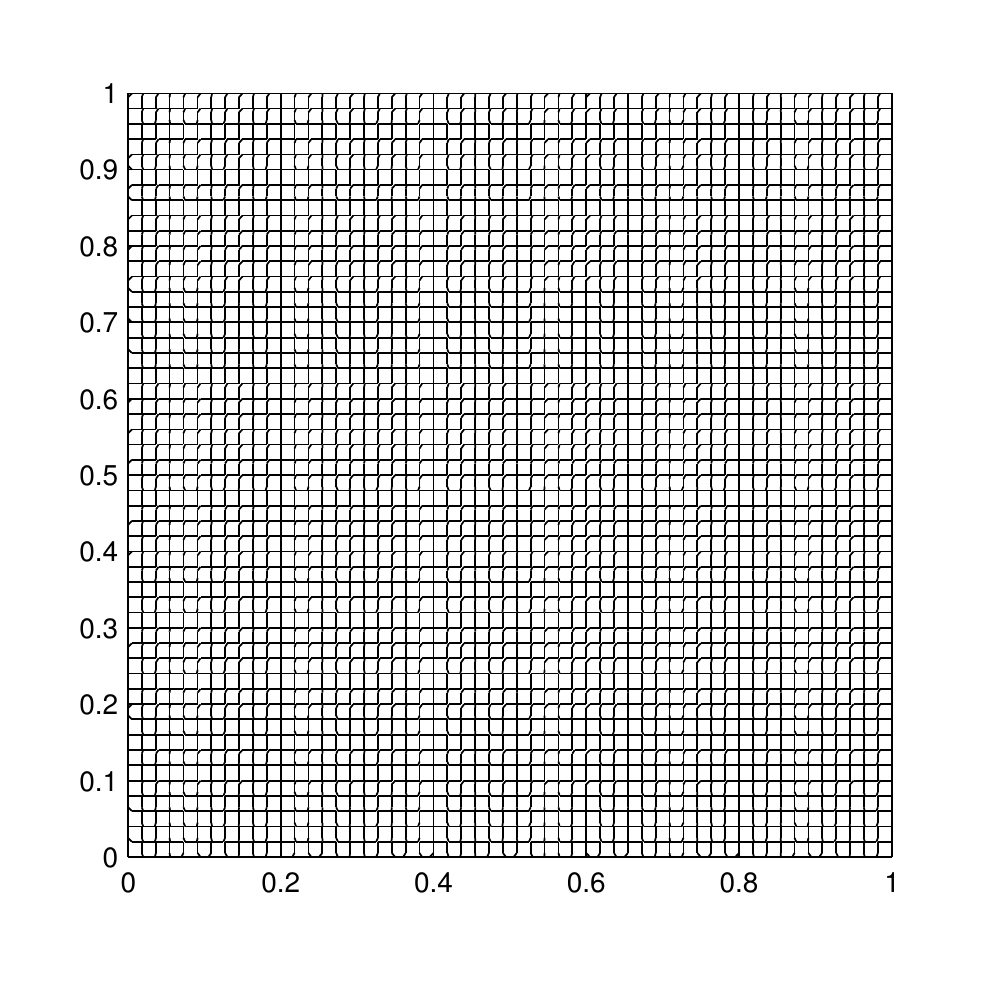}
\end{center}
\caption{
The nodal lines of a Dirichlet eigenfunction $u(x_1,x_2)$ on the unit
square, with the associated eigenvalue $\lambda = 5525\pi^2$ of
multiplicity 12.  The eigenfunction was obtained as a linear
combination of terms $\sin(\pi n_1 x_1)\sin(\pi n_2 x_2)$, with $n_1^2
+ n_2^2 = 5525$ and randomly chosen coefficients.  For comparison,
another eigenfunction with the same eigenvalue, $\sin(50\pi
x_1)\sin(55 \pi x_2)$, is shown.}
\label{fig:nodal}
\end{figure}

Even for such a simple domain as a square, the nodal lines and domains
may have complicated structure, especially for high-frequency
eigenfunctions (Fig. \ref{fig:nodal}).  This is particularly true for
degenerate eigenfunctions for which one can ``tune'' the coefficients
of the corresponding linear combination to modify continuously the
nodal lines.

Pleijel sharpened the Courant's theorem by showing that the upper
bound $m$ for the number $\nu_m$ of nodal domains is attained only for
a finite number of eigenfunctions \cite{Pleijel56}.  Moreover, he
obtained the upper limit: $\varlimsup\limits_{m\to\infty} \nu_m/m =
4/j_{0,1}^2 \approx 0.691...$.  Note that Lewy constructed spherical
harmonics of any degree $n$ whose nodal sets have one component for
odd $n$ and two components for even $n$ implying that no non-trivial
lower bound for $\nu_m$ is possible \cite{Lewy77}.  We also
mention that the counting of nodal domains can be viewed as
partitioning of the domain into a fixed number of subdomains and
minimizing an appropriate ``energy'' of the partition (e.g., the
maximum of the ground state energies of the subdomains).  When a
partition corresponds to an eigenfunction, the ground state energies
of all the nodal domains are the same, i.e., it is an equipartition
\cite{Berkolaiko12}.

Blum {\it et al.} considered the distribution of the (properly
normalized) number of nodal domains of the Dirichlet-Laplacian
eigenfunctions in 2D quantum billiards and showed the existence of the
limiting distribution in the high-frequency limit (i.e., when
$\lambda_m\to\infty$) \cite{Blum02}.  These distributions were argued
to be universal for systems with integrable or chaotic classical
dynamics that allows one to distinguish them and thus provides a new
criterion for quantum chaos (see Sec. \ref{sec:quantum_billiard}).  It
was also conjectured that the distribution of nodal domains for
chaotic systems coincides with that for Gaussian random functions.

Bogomolny and Schmit proposed a percolation-like model to describe the
nodal domains which permitted to perform analytical calculations and
agreed well with numerical simulations \cite{Bogomolny02}.  This model
allows one to apply ideas and methods developed within the percolation
theory \cite{Stauffer} to the field of quantum chaos.  Using the
analogy with Gaussian random functions, Bogomolny and Schmit obtained
that the mean and variance of the number $\nu_m$ of nodal domains grow
as $m$, with explicit formulas for the prefactors.  From the
percolation theory, the distribution of the area $s$ of the connected
nodal domains was conjectured to follow a power law, $n(s)\propto
s^{-187/91}$, as confirmed by simulations \cite{Bogomolny02}.  In the
particular case of random Gaussian spherical harmonics, Nazarov and
Sodin rigorously derived the asymptotic behavior for the number
$\nu_n$ of nodal domains of the harmonic of degree $n$
\cite{Nazarov09}.  They proved that as $n$ grows to infinity, the mean
of $\nu_n/n^2$ tends to a positive constant, and that $\nu_n/n^2$
exponentially concentrates around this constant (we recall that the
associate eigenvalue is $n(n+1)$).

The geometrical structure of nodal lines and domains has been
intensively studied (see \cite{Nazarov05,Polterovich07} for further
discussion of the asymptotic nodal geometry).  For instance, the
length of the nodal line of an eigenfunction of the Laplace operator
in two-dimensional Riemannian manifolds was separately investigated by
Br\"uning, Yao and Nadirashvili who obtained its lower and upper
bounds \cite{Bruning78,Yau,Nadirashvili88a}.  In addition, a number of
conjectures about the properties of particular eigenfunctions were
discussed in the literature.  We mention three of them:

(i) In 1967, Payne conjectured that the second Dirichlet eigenfunction
$u_2$ cannot have a closed nodal line in a bounded planar domain
\cite{Payne67,Payne73c}.  This conjecture was proved for convex domains
\cite{Alessandrini94,Melas92} and disproved by non-convex domains
\cite{Hoffmann97}, see also \cite{Jerison95,Grieser96}.

(ii) The hot spots conjecture formulated by J. Rauch in 1974 says that
the maximum of the second Neumann eigenfunction is attained at a
boundary point.  This conjecture was proved by Banuelos and Burdzy for
a class of planar domains \cite{Banuelos99} but in general the
statement is wrong, as shown by several counter-examples
\cite{Burdzy99,Jerison00,Bass00,Burdzy05}.

(iii) Liboff formulated several conjectures; one of them states that
the nodal surface of the first-excited state of a 3D convex domain
intersects its boundary in a single simple closed curve
\cite{Liboff94a}.

The analysis of nodal lines that describe zeros of eigenfunctions, can
be extended to other level sets.  For instance, a level set of the
first Dirichlet eigenfunction $u_1$ on a bounded convex domain $\Omega
\in \R^d$ is itself convex \cite{Kawohl}.  Grieser and Jerison
estimated the size of the first eigenfunction uniformly for all convex
domains \cite{Grieser98}.  In particular, they located the place where
$u_1$ achieves its maximum to within a distance comparable to the
inradius, uniformly for arbitrarily large diameter.  Other geometrical
characteristics (e.g., the volume of a set on which an eigenfunction
is positive) can also be analyzed \cite{Nadirashvili91}.

\section{Estimates for Laplacian eigenfunctions}
\label{sec:estimates}

The ``amplitudes'' of eigenfunctions can be characterized either
globally by their $L_p$ norms
\begin{equation}
\|u\|_p \equiv \left(\int\limits_\Omega dx ~|u(x)|^p \right)^{1/p}  \qquad (p \geq 1) ,
\end{equation}
or locally by pointwise estimates.  Since eigenfunctions are defined
up to a multiplicative constant, one often uses $L_2(\Omega)$
normalization: $\|u\|_2 = 1$.  Note also the limiting case of
$L_\infty$-norm
\begin{equation}
\|u\|_\infty \equiv \esssup\limits_{x\in \Omega} |u(x)| = \max\limits_{x\in \Omega} |u(x)| 
\end{equation}
(the first equality is the general definition, while the second
equality is applicable for eigenfunctions).  It is worth recalling
H\"older's inequality for any two measurable functions $u$ and $v$ and
for any positive $p$, $q$ such that $1/p + 1/q = 1$:
\begin{equation}
\label{eq:Holder}
\|u v\|_1 \leq \|u\|_p \|v\|_q  .
\end{equation}
For a bounded domain $\Omega\subset\R^d$ (with a finite Lebesgue
measure $\mu_d(\Omega)$), H\"older's inequality implies
%
\begin{equation}
\|u\|_p \leq [\mu_d(\Omega)]^{\frac{1}{p} - \frac{1}{p'}} \|u\|_{p'}  \qquad (1 \leq p \leq p') .
\end{equation}
We also mention Minkowski's inequality for two measurable
functions and any $p \geq 1$:
\begin{equation}
\label{eq:Minkowski}
\|u + v\|_p \leq \|u\|_p  + \|v\|_p  .
\end{equation}

\subsection{First (ground) Dirichlet eigenfunction}
\label{sec:first_eigenfunction}

The Dirichlet eigenfunction $u_1$ associated with the first eigenvalue
$\lambda_1 > 0$ does not change the sign in $\Omega$ and may be taken
to be positive.  It satisfies the following inequalities.

(i) Payne and Rayner showed in two dimensions that
\begin{equation}
\label{eq:Payne-Rayner}
\|u_1\|_2 \leq \frac{\sqrt{\lambda_1}}{\sqrt{4\pi}} ~\|u_1\|_1 ,
\end{equation}
with equality if and only if $\Omega$ is a disk
\cite{Payne72,Payne73}.  Kohler-Jobin gave an extension of this
inequality to higher dimensions
\cite{Kohler77} (see \cite{Payne73,Kohler81,Chiti82} for other extensions):
\begin{equation}
\|u_1\|_2 \leq \frac{\lambda_1^{d/4}}{\sqrt{2d\omega_d [j_{\frac{d}{2}-1,1}]^{d-2}}} ~ \| u_1\|_1 .
\end{equation}

(ii) Payne and Stakgold derived two inequalities for a convex domain
in 2D
\begin{equation}
\frac{\pi}{2\mu_2(\Omega)} \|u_1\|_1 \leq \|u_1\|_\infty
\end{equation}
and
\begin{equation}
u_1(x) \leq |x-\pa| \frac{\sqrt{\lambda_1}}{\mu_2(\Omega)} \|u_1\|_1  \qquad (x\in \Omega),
\end{equation}
where $|x-\pa|$ is the distance from a point $x$ in $\Omega$ to the
boundary $\pa$ \cite{Payne73b}.

(iii) van den Berg proved the following inequality for
$L_2$-normalized eigenfunction $u_1$ when $\Omega$ is an open, bounded
and connected set in $\R^d$ ($d=2,3, \dots$):
\begin{equation}
\|u_1\|_{\infty} \leq \frac{2^{\frac{2-d}{2}}}{\pi^{d/4} \sqrt{\Gamma(d/2)}} ~ 
 \frac{\left({j_{\frac{d-2}{2},1}}\right)^{\frac{d-2}{2}}}{|{J_{\frac{d}{2}}\bigl(j_{\frac{d-2}{2},1}\bigr)}|}~ \rho^{-d/2} ,
\end{equation}
with equality if and only if $\Omega$ is a ball, where $\rho$ is the
inradius (Eq. (\ref{eq:inradius})) \cite{VanDenBerg00}.  van den Berg
also conjectured the stronger inequality for an open bounded convex
domain $\Omega\subset \R^d$:
\begin{equation}
\|u_1\|_{\infty} \leq C_d \rho^{-d/2} (\rho/\delta)^{1/6} ,
\end{equation}
where $\delta$ is the diameter of $\Omega$, and $C_d$ is a universal
constant independent of $\Omega$.

(iv) Kr\"oger obtained another upper bound for $\|u_1\|_\infty$ for a
convex domain $\Omega \subset \R^d$.  Suppose that $\lambda_1(D) \geq
\Lambda(\delta)$ for every convex subdomain $D\subset \Omega$ with
$\mu_d(D) \leq \delta \mu_d(\Omega)$ and positive numbers $\delta$ and
$\Lambda(\delta)$.  The first eigenfunction $u_1$ which is normalized
such that $\|u_1\|_2^2 = \mu_d(\Omega)$, satisfies
\begin{equation}
\|u_1\|_\infty \leq C_d \delta^{-1/2} \bigl[1 + \ln \|u_1\|_\infty - \ln(1 - \lambda_1/\Lambda(\delta))\bigr]^{d/2} ,
\end{equation}
with a universal positive constant $C_d$ which depends only on the
dimension $d$ \cite{Kroger96}.

(v) Pang investigated how the first Dirichlet eigenvalue and
eigenfunction would change when the domain slightly shrinks
\cite{Pang96,Pang97}.  For a bounded simply connected open set
$\Omega\subset \R^2$, let
\begin{equation*}
\Omega_\epsilon \supseteq \left\{ x \in \Omega : |x - \pa| \ge \epsilon \right\}
\end{equation*}
be its interior, i.e., $\Omega$ without an $\epsilon$ boundary layer.
Let $\lambda_m^\epsilon$ and $u_m^\epsilon$ be the Dirichlet
eigenvalues and $L_2$-normalized eigenfunctions in $\Omega_\epsilon$
(with $\lambda_m^0 = \lambda_m$ and $u_m^0 = u_m$ referring to the
original domain $\Omega$).  Then, for all $\epsilon \in (0, \rho/2)$,
\begin{equation}
\begin{split}
|\lambda_1^\epsilon - \lambda_1| & \leq C_1 \epsilon^{1/2}, \\
\|u_1 - T_{\epsilon} u_1^\epsilon\|_{L_\infty(\Omega)} & \leq 
\left[{C_2 + C_3(\lambda_2 - \lambda_1)^{-1/2} + C_4(\lambda_2 - \lambda_1)^{-1}}\right] \epsilon^{1/2} , \\
\end{split}
\end{equation}
where $\rho$ is the inradius of $\Omega$ (Eq. (\ref{eq:inradius})),
$T_\epsilon$ is the extension operator from $\Omega_\epsilon$ to
$\Omega$, and
\begin{eqnarray*}
C_1 &=& \rho^{-3/2} \beta^{9/4} ~ \frac{2^9\gamma_1^4}{3\pi^{9/4}} ,   \qquad
C_2 = \rho^{-3/2} \beta^{13/4} ~ \frac{2^{12}\gamma_1^{5}}{\pi^{15/4}} , \\
C_3 &=& \rho^{-5/2} \beta^4 \biggl(\frac{2^{15}\gamma_1^6\gamma_2}{3\sqrt{2\alpha} \pi^{9/2}}\biggr)
\biggl[1 + \frac{9\gamma_1}{\pi^{3/4}} \beta^{3/4}\biggr], \\
C_4 &=& \rho^{-7/2} \beta^7 \biggl(\frac{2^{26} \gamma_1^{10} \gamma_2^2}{81\sqrt{2} ~ \alpha ~ \pi^{15/2}} \biggr)
\biggl[1 + 18 \gamma_1 \beta^{3/4} + \frac{81\gamma_1^2}{\pi^{3/2}} \beta^{3/2}\biggr] ,
\end{eqnarray*}
where $\beta = \mu_2(\Omega)/\rho^2$, $\alpha$ is the constant from
Eq. (\ref{eq:lambda1D_L}) (for which one can use the best known
estimate $\alpha = 0.6197...$ from \cite{Banuelos94}), and $\gamma_1$
and $\gamma_2$ are the first and second Dirichlet eigenvalues for the
unit disk: $\gamma_1 = j_{0,1}^2 \approx 5.7832$ and $\gamma_2 =
j_{1,1}^2
\approx 14.6820$.  Moreover, when $\Omega$ is the cardioid in $\R^2$,
the term $\epsilon^{1/2}$ cannot be improved.%
\footnote{
In the original paper \cite{Pang97}, the coefficient $C_4$ in
Eq. (1.5) should be multiplied by the omitted prefactor $\sqrt{2}
|\Omega|$ that follows from the derivation. }

In addition, Davies proved that for a bounded simply connected open
set $\Omega\in \R^2$ and for any $\beta \in (0,1/2)$, there exists
$c_\beta\ge 1$ such that \cite{Davies93}
\begin{equation}
|\lambda_1^{\epsilon} - \lambda_1| \leq c_\beta \epsilon^{\beta}
\end{equation}
for all sufficiently small $\epsilon > 0$.  The estimate also holds
for higher Dirichlet eigenvalues.

\subsection{Estimates applicable for all eigenfunctions}

\subsubsection{Estimates through the Green function}
\label{sec:estimate_Green}

Using the spectral decomposition (\ref{eq:G_decomp}) of the Green
function $G(x,y)$, one can rewrite Eq. (\ref{eq:eigen}) as
\begin{equation*}
u_m(x) = \lambda_m \int\limits_\Omega G(x,y) u_m(y) dy ,
\end{equation*}
from which H\"older inequality (\ref{eq:Holder}) yields a family
of simple pointwise estimates
\begin{equation}
|u_m(x)| \leq \lambda_m \|u_m\|_{\frac{p}{p-1}} \left(\int\limits_\Omega |G(x,y)|^p dy \right)^{1/p},
\end{equation}
with any $p \geq 1$.  Here, a single function of $x$ in the right-hand
side bounds all the eigenfunctions.  In particular, for $p = 1$, one
gets
\begin{equation}
\label{eq:estim_U1}
|u_m(x)| \leq \lambda_m \|u_m\|_\infty \int\limits_\Omega |G(x,y)| dy .
\end{equation}
For Dirichlet boundary condition, $G(x,y)$ is positive everywhere in
$\Omega$ so that
\begin{equation}
\label{eq:estim_U}
|u_m(x)| \leq \lambda_m \|u_m\|_\infty ~U(x),  \qquad  U(x) = \int\limits_\Omega G(x,y) dy,
\end{equation}
where $U(x)$ solves the boundary value problem
\begin{equation}
-\Delta U(x) = 1  \quad (x\in\Omega), \qquad U(x) = 0 \quad (x\in\pa) .
\end{equation}
The solution of this equation is known to be the mean first passage
time to the boundary $\pa$ from an interior point $x$ \cite{Redner}.
The inequalities (\ref{eq:estim_U1}, \ref{eq:estim_U}) (or their
extensions) were reported by Moler and Payne \cite{Moler68}
(Sect. \ref{sec:Moler}) and were used by Filoche and Mayboroda for
determining the geometrical structure of eigenfunctions
\cite{Filoche12} (Sect. \ref{sec:Filoche}).  Note that the function
$U(x)$ was also considered by Gorelick {\it et al.} for a reliable
extraction of various shape properties of a silhouette, including part
structure and rough skeleton, local orientation and aspect ratio of
different parts, and convex and concave sections of the boundaries
\cite{Gorelick06}.

\subsubsection{Bounds for eigenvalues and eigenfunctions of symmetric operators}
\label{sec:Moler}

Moler and Payne derived simple bounds for eigenvalues and
eigenfunctions of symmetric operators by considering their extensions
\cite{Moler68}.  As a typical example, one can think of the
Dirichlet-Laplace operator in a bounded domain $\Omega$ (symmetric
operator $A$) and of the Laplace operator without boundary conditions
(extension $A_*$).  An approximation to an eigenvalue and
eigenfunction of $A$ can be obtained by solving a simpler eigenvalue
problem $A_* u_* = \lambda_* u_*$ without boundary condition.  If
there exists a function $w$ such that $A_* w = 0$ and $w = u_*$ at the
boundary of $\Omega$ and if
$\displaystyle \ve = \frac{\| w\|_{L_2(\Omega)}}{\|u_*\|_{L_2(\Omega)}} < 1$, 
then there exists an eigenvalue $\lambda_k$ of $A$ satisfying
\begin{equation}
\label{eq:Moler}
\frac{|\lambda_*|}{1+\ve} \leq |\lambda_k| \leq \frac{|\lambda_*|}{1-\ve} .
\end{equation}
Moreover, if $\|u_*\|_{L_2(\Omega)} = 1$ and $u_k$ is the
$L_2$-normalized projection of $u_*$ onto the eigenspace of
$\lambda_k$, then
\begin{equation}
\label{eq:Moler2}
\| u_* - u_k\|_{L_2(\Omega)} \leq \frac{\ve}{\alpha}\left(1 + \frac{\ve^2}{\alpha^2}\right)^{1/2} ,
\qquad {\rm with} \quad \alpha = \min\limits_{\lambda_n \ne \lambda_k} \frac{|\lambda_n - \lambda_*|}{|\lambda_n|} .
\end{equation}

If $u_*$ is a good approximation to an eigenfunction of the
Dirichlet-Laplace operator, then it must be close to zero on the
boundary of $\Omega$, yielding small $\ve$ and thus accurate lower and
upper bounds in (\ref{eq:Moler}).  The accuracy of the bound
(\ref{eq:Moler2}) also depends on the separation $\alpha$ between
eigenvalues.

In the same work, Moler and Payne also provided pointwise bounds for
eigenfunctions that rely on Green's functions (an extension of
Sec. \ref{sec:estimate_Green}).

\subsubsection{Estimates for $L_p$-norms}

Chiti extended Payne-Rayner's inequality (\ref{eq:Payne-Rayner}) to
the eigenfunctions of linear elliptic second order operators in
divergent form, with Dirichlet boundary condition \cite{Chiti82}.  For
the Laplace operator in a bounded domain $\Omega\subset \R^d$, Chiti's
inequality for any real numbers $q \geq p > 0$ states:
\begin{equation}
\|u \|_q \leq  \|u \|_p ~ (d\omega_d)^{\frac{1}{q}-\frac{1}{p}} \lambda^{\frac{q-p}{2pq}d} 
\frac{\biggl(\int\limits_0^{j_{\frac{d}{2}-1,1}} dr ~ r^{d-1 + q(1-d/2)} [J_{\frac{d}{2}-1}(r)]^q\biggr)^{1/q}}
{\biggl(\int\limits_0^{j_{\frac{d}{2}-1,1}} dr ~ r^{d-1 + p(1-d/2)} [J_{\frac{d}{2}-1}(r)]^p\biggr)^{1/p}} ,
\end{equation}
where $\omega_d$ is given by Eq. (\ref{eq:omega_d}).

\subsubsection{Pointwise bounds for Dirichlet eigenfunctions}

Banuelos derived a pointwise upper bound for $L_2$-normalized
Dirichlet eigenfunctions \cite{Banuelos96}
\begin{equation}
|u_m(x)| \leq \lambda_m^{d/4}  \qquad (x\in \Omega).
\end{equation}

van den Berg and Bolthausen proved several estimates for
$L_2$-normalized Dirichlet eigenfunctions \cite{VanDenBerg99}.  Let
$\Omega \subset \R^d$ $(d = 2,3,\dots)$ be an open bounded domain with
boundary $\pa$ which satisfies an $\alpha$-uniform capacitary density
condition with some $\alpha \in (0,1]$, i.e.
\begin{equation}
\label{eq:cap_cond}
\textrm{Cap}\{ \pa \cap B(x;r)\} \geq \alpha \textrm{Cap}\{B(x;r)\}, \quad x \in \pa,~ 0<r<\delta ,
\end{equation}
where $B(x,r)$ is the ball of radius $r$ centered at $x$, $\delta$ is
the diameter of $\Omega$ (Eq. (\ref{eq:diameter})), and $\textrm{Cap}$
is the logarithmic capacity for $d = 2$ and the Newtonian (or
harmonic) capacity for $d > 2$ (the harmonic capacity of an Euclidean
domain presents a measure of its ``size'' through the total charge the
domain can hold at a given potential energy \cite{Doob}).  The
condition (\ref{eq:cap_cond}) guarantees that all points of $\pa$ are
regular.  The following estimates hold

(i) in two dimensions ($d = 2$), for all $m=1,2,\dots$ and all $x \in
\Omega$ such that $|x-\pa| \sqrt{\lambda_m} < 1$, one has 
\begin{equation}
|u_m(x)| \leq \left\{\frac{6\lambda_m \ln(\alpha^{2\pi}/2)}{\ln{\left({|x - \pa|\sqrt{\lambda_m}}\right)}}\right\}^{1/2} .
\end{equation}

(ii) in higher dimensions ($d > 2$), for all $m=1,2,\dots$ and all $x
\in \Omega$ such that
\begin{equation}
|x-\pa| \sqrt{\lambda_m} \le \left({\frac{\alpha^6}{2^{13}}}\right)^{1+\gamma(d-1)/(d-2)} ,
\end{equation}
with
$\displaystyle \gamma = \frac{3^{-d-1}\alpha}{\ln(2(2/\alpha)^{1/(d-2)})}$,
one has
\begin{equation}
|u_m(x)| \leq 2 \lambda_m^{d/4} \biggl(|x-\pa|\sqrt{\lambda_m}\biggr)^{\frac12((1/\gamma)+(d-1)/(d-2))^{-1}} .
\end{equation}

(iii) for a planar simply connected domain and all $m= 1, 2, \dots$,
\begin{equation}
|u_m(x)| \leq m~ 2^{9/2} \pi^{1/4} \frac{(\mu_2(\Omega))^{1/4}}{\rho^2}~ |x-\pa|^{1/2}  \qquad (x\in\Omega),
\end{equation}
where $\rho$ is the inradius of $\Omega$ (see
Eq. (\ref{eq:inradius})), and the inequality is sharp.

\subsubsection{Upper and lower bounds for normal derivatives of Dirichlet eigenfunctions}

Suppose that $M$ is a compact Riemannian manifold with boundary and
$u$ is an $L_2$-normalized Dirichlet eigenfunction with eigenvalue
$\lambda$.  Let $\psi$ be its normal derivative at the boundary.  A
scaling argument suggests that the $L_2$-norm of $\psi$ will grow as
$\sqrt{\lambda}$ as $\lambda \to \infty$.  Hassell and Tao proved that
\begin{equation}
 c_M \sqrt{\lambda} \leq \|\psi \|_{L_2(\partial M)} \leq C_M \sqrt{\lambda} ,
\end{equation}
where the upper bound holds for any Riemannian manifold, while the
lower bound is valid provided that $M$ has no trapped geodesics
\cite{Hassell02}.  The positive constants $c_M$ and $C_M$ depend
on $M$, but not on $\lambda$.

\subsubsection{Estimates for restriction onto a subdomain}
\label{sec:Filoche}

For a bounded domain $\Omega\subset\R^d$, Filoche and Mayboroda
obtained the upper bound for the $L_2$-norm of a Dirichlet-Laplacian
eigenfunction $u$ associated to $\lambda$, in any open subset
$D\subset \Omega$ \cite{Filoche12}:
\begin{equation}
\label{eq:Mayboroda}
\|u\|_{L_2(D)} \leq \left(1 + \frac{\lambda}{d_D(\lambda)} \right) \|v\|_{L_2(D)} ,
\end{equation}
where the function $v$ solves the boundary value problem in $D$:
\begin{equation*}
\Delta v = 0  \quad (x\in D), \qquad    v = u \quad (x\in \paD).
\end{equation*}
and $d_D(\lambda)$ is the distance from $\lambda$ to the spectrum of
the Dirichlet-Laplace operator in $D$.  Note also that the above bound
was proved for general self-adjoint elliptic operators
\cite{Filoche12}.  When combined with Eq.  (\ref{eq:estim_U}), this
inequality helps one investigate the spatial distribution of
eigenfunctions because harmonic functions are in general much easier
to analyze or estimate than eigenfunctions.

We complete the above estimate by a lower bound \cite{Nguyen}
\begin{equation}
\label{eq:estim_dg1}
\|u\|_{L_2(D)} \geq \frac{\lambda_1(D)}{\lambda + \lambda_1(D)} \|v\|_{L_2(D)} ,
\end{equation}
where $\lambda_1(D)$ is the first Dirichlet eigenvalue of the
subdomain $D$.


\section{Localization of eigenfunctions}
\label{sec:localization}

``Localization'' is defined in the Webster's dictionary as ``act of
localizing, or state of being localized''.  The notion of localization
appears in various fields of science and often has different meanings.
Throughout this review, a function $u$ defined on a domain
$\Omega\subset \R^d$, is called $L_p$-\emph{localized} (for $p\geq 1$)
if there exists a bounded subset $\Omega_0\subset\Omega$ which
supports almost all $L_p$-norm of $u$, i.e.
\begin{equation}
\label{eq:def_loc}
\frac{\| u\|_{L_p(\Omega\setminus \Omega_0)}}{\| u\|_{L_p(\Omega)}} \ll 1  \qquad  \textrm{and} \qquad
\frac{\mu_d(\Omega_0)}{\mu_d(\Omega)} \ll 1 .
\end{equation}
Qualitatively, a localized function essentially ``lives'' on a small
subset of the domain and takes small values on the remaining part.
For instance, a Gaussian function $\exp(-x^2)$ on $\Omega = \R$ is
$L_p$-localized for any $p\geq 1$ since one can choose $\Omega_0 =
[-a,a]$ with large enough $a$ so that the ratio of $L_p$-norms can be
made arbitrarily small, while the ratio of lengths
$\mu_1(\Omega_0)/\mu_1(\Omega)$ is strictly $0$.  In turn, when
$\Omega = [-A,A]$, the localization character of $\exp(-x^2)$ on
$\Omega$ becomes dependent on $A$ and thus conventional.  A simple
calculation shows that both ratios in (\ref{eq:def_loc}) cannot be
simultaneously made smaller than $1/(A+1)$ for any $p\geq 1$.  For
instance, if $A = 3$ and the ``threshold'' 1/4 is viewed small enough,
then we are justified to call $\exp(-x^2)$ localized on $[-A, A]$.
This example illustrates that the above inequalities do not provide a
universal quantitative criterion to distinguish localized from
non-localized (or extended) functions.  In this section, we will
describe various kinds of localization for which some quantitative
criteria can be formulated.  We will also show that the choice of the
norm (i.e., $p$) may be important.

Another ``definition'' of localization was given by Felix {\it et al.}
who combined $L_2$ and $L_4$ norms to define the ``existence area'' as
\cite{Felix07}
\begin{equation}
S(u) = \frac{\|u\|_{L_2(\Omega)}^4}{\|u\|_{L_4(\Omega)}^4} .
\end{equation}
A function $u$ was called localized when its existence area $S(u)$ was
much smaller than the area $\mu_2(\Omega)$ \cite{Felix07} (this
definition trivially extends to other dimensions).  In fact, if a
function is small in a subdomain, the fourth power diminishes it
stronger than the second power.  For instance, if $\Omega = (0,1)$ and
$u$ is $1$ on the subinterval $\Omega_0 = (1/4,1/2)$ and $0$
otherwise, one has $\|u\|_{L_2(\Omega)} = 1/2$ and
$\|u\|_{L_4(\Omega)} = 1/\sqrt{2}$ so that $S(u) = 1/4$, i.e., the
length of the subinterval $\Omega_0$.  This definition is still
qualitative: e.g., in the above example, is the ratio
$S(u)/\mu_1(\Omega) = 1/4$ small enough to call $u$ localized?  Note
that a whole family of ``existence areas'' can be constructed by
comparing $L_p$ and $L_q$ norms (with $p < q$),
\begin{equation}
S_{p,q}(u) = \left(\frac{\|u\|_{L_p(\Omega)}}{\|u\|_{L_q(\Omega)}}\right)^{\frac{1}{\frac{1}{p} - \frac{1}{q}}} .
\end{equation}

\subsection{Bound quantum states in a potential}
\label{sec:potential}

The notion of bound, trapped or localized quantum states is known for
a long time \cite{Reed,Blank}.  The simplest ``canonical'' example is
the quantum harmonic oscillator, i.e., a particle of mass $m$ in a
harmonic potential of frequency $\omega$ which is described by the
Hamiltonian
\begin{equation}
H = \frac{\hat{p}^2}{2m} + \frac{m \omega^2}{2} \hat{x}^2 = - \frac{\hbar^2}{2m}\partial_x^2 + \frac{m \omega^2}{2} x^2 ,
\end{equation}
where $\hat{p} = -i\hbar \partial_x$ is the momentum operator, and
$\hat{x} = x$ is the position operator ($\hbar$ being the Planck's
constant).  The eigenfunctions of this operator are well known:
\begin{equation}
\psi_n(x) = \sqrt{\frac{1}{2^n~ n!}} \left(\frac{m\omega}{\pi \hbar}\right)^{1/4} \exp\left(-\frac{m\omega x^2}{2\hbar}\right) 
H_n\bigl(\sqrt{m\omega/\hbar} x\bigr),
\end{equation}
where $H_n(x)$ are the Hermite polynomials.  All these functions are
localized in a region around the minimum of the harmonic potential
(here, $x = 0$), and rapidly decay outside this region.  For this
example, the definition (\ref{eq:def_loc}) of localization is
rigorous.  In physical terms, the presence of a strong potential
forbids the particle to travel far from the origin, the size of the
localization region being $\sqrt{\hbar/(m\omega)}$.  This so-called
strong localization has been thoroughly investigated in physics and
mathematics \cite{Maslov,Maslov64,Agmon,Reed,Simon,Schnol57,Schwinger61,Orocko74}.

\subsection{Anderson localization}
\label{sec:Anderson}

\begin{figure}
\begin{center}
\includegraphics[width=130mm]{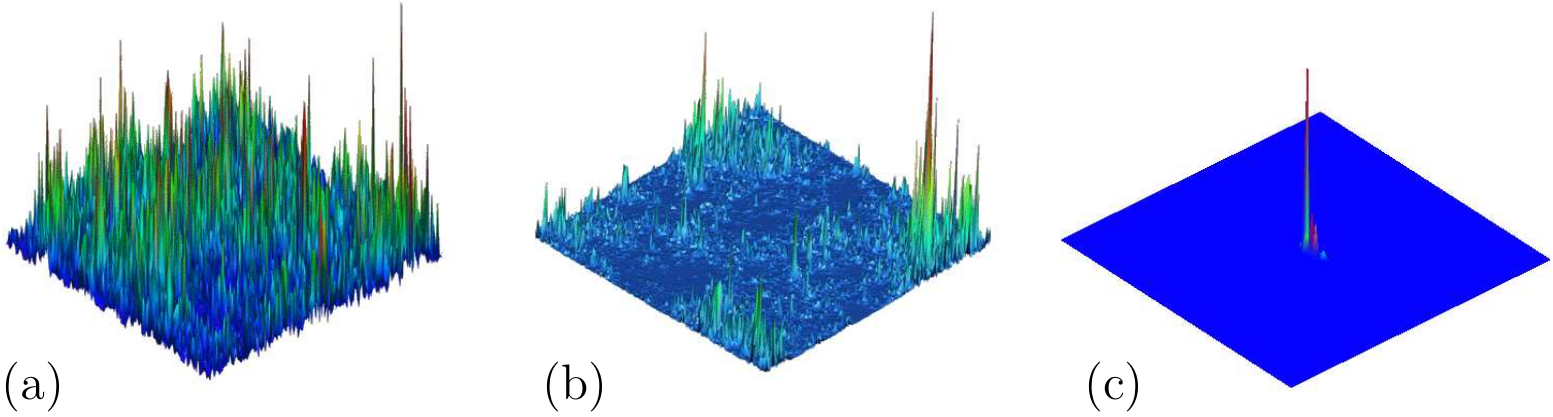}
\end{center}
\caption{
Illustration of the Anderson transition in a tight-binding model (or
so-called SU(2) model) in the two-dimensional symplectic class
\cite{Asada02,Asada04,Obuse04,Obuse07}.  Three shown eigenfunctions
(with the energy close to $1$) were computed for three disorder
strengths $W$ that correspond to (a) metallic state ($W < W_0$), (b)
critical state ($W = W_0$), and (c) insulating state ($W > W_0$), $W_0
= 5.952$ being the critical disorder strength.  The latter
eigenfunction is strongly localized that prohibits diffusion of charge
carriers (i.e., no electric current).  The eigenfunctions were
computed and provided by H. Obuse (unpublished earlier). }
\label{fig:Anderson}
\end{figure}

The previous example of a single quantum harmonic well is too
idealized.  A piece of matter contains an extremely large number of
interacting atoms.  Even if one focuses onto a single atom in an
effective potential, the form of this potential may be so complicated
that the study of the underlying eigenfunctions would in general be
intractable.  In 1958, Anderson considered a lattice model for a
charge carrier in a random potential and proved the localization of
eigenfunctions under certain conditions \cite{Anderson58}.  The
localization of charge carriers means no electric current through the
medium (insulating state), in contrast to metallic or conducting state
when the charge carriers are not localized.  The Anderson transition
between insulating and conducting states is illustrated for the
tight-binding model on Fig. \ref{fig:Anderson}.  The shown
eigenfunctions were computed by Obuse for three disorder strengths $W$
that correspond to metallic ($W < W_0$), critical ($W = W_0$), and
insulating ($W > W_0$) states, $W_0 = 5.952$ being the critical
disorder strength.  The latter eigenfunction is strongly localized
that prohibits diffusion of charge carriers (i.e., no electric
current).  The Anderson localization which explains the
metal-insulator transitions in semiconductors, was thoroughly
investigated during the last fifty years (see
\cite{Thouless74,Lee85,Kramer93,Belitz94,delRio95,Mirlin00,Evers08,Stolz11,Stollman}
for details and references).  Similar localization phenomena were
observed for microwaves with two-dimensional random scattering
\cite{Dalichaouch91}, for light in a disordered medium
\cite{Wiersma97} and in disordered photonic crystals
\cite{Schwartz07,Sapienza10}, for matter waves in a controlled
disorder \cite{Billy08} and in non-interacting Bose-Einstein
condensate \cite{Roati08}, and for ultrasound \cite{Hu08}.  The
multifractal structure of the eigenfunctions at the critical point
(illustrated by Fig. \ref{fig:Anderson}b) has also been intensively
investigated (see \cite{Evers08,Gruzberg11} and references therein).
Localization of eigenstates and transport phenomena in
one-dimensional disordered systems are reviewed in \cite{Ishii73}.  An
introduction to wave scattering and the related localization is given
in \cite{Sheng}.

\subsection{Trapping in infinite waveguides}
\label{sec:trapping}

In both previous cases, localization of eigenfunctions was related to
an external potential.  In particular, if the potential was not strong
enough, Anderson localization could disappear
(Fig. \ref{fig:Anderson}a).  Is the presence of a potential necessary
for localization?  The formal answer is positive because the
eigenstates of the Laplace operator in the whole space $\R^d$ are
simply $e^{i(k\cdot x)}$ (parameterized by the vector $k$) which are
all extended in $\R^d$.  These waves are called ``resonances'' (not
eigenfunctions) of the Laplace operator, as their $L_2$-norm is
infinite.

The situation is different for the Laplace operator in a bounded
domain with Dirichlet boundary condition.  In quantum mechanics, such
a boundary presents a ``hard wall'' that separates the interior of the
domain with zero potential from the exterior of the domain with
infinite potential.  For instance, this ``model'' was employed by
Crommie {\it et al.} to describe the confinement of electrons to
quantum corrals on a metallic surface \cite{Crommie93} (see also their
figure 2 that shows the experimental spatial structure of the
electron's wavefunction).  Although the physical interpretation of a
boundary through an infinite potential is instructive, we will use the
mathematical terminology and speak about the eigenvalue problem for
the Laplace operator in a bounded domain without potential.

For unbounded domains, the spectrum of the Laplace operator consists
of two parts: (i) the discrete (or point-like) spectrum, with
eigenfunctions of finite $L_2$ norm that are necessarily ``trapped''
or ``localized'' in a bounded region of the waveguide, and (ii) the
continuous spectrum, with associated functions of infinite $L_2$ norm
that are extended over the whole domain.  The continuous spectrum may
also contain embedded eigenvalues whose eigenfunctions have finite
$L_2$ norm.  A wave excited at the frequency of the trapped eigenmode
remains in the localization region and does not propagate.  In this
case, the definition (\ref{eq:def_loc}) of localization is again
rigorous, as for any bounded subset $\Omega_0$ of an unbounded domain
$\Omega$, one has $\mu_d(\Omega_0)/\mu_d(\Omega) = 0$, while the ratio
of $L_2$ norms can be made arbitrarily small by expanding $\Omega_0$.

This kind of localization in classical and quantum waveguides has been
thoroughly investigated (see reviews \cite{Duclos95,Linton07} and also
references in \cite{Olendski10}).  In the seminal paper, Rellich
proved the existence of a localized eigenfunction in a deformed
infinite cylinder \cite{Rellich}.  His results were significantly
extended by Jones \cite{Jones53}.  Ursell reported on the existence of
trapped modes in surface water waves in channels
\cite{Ursell51,Ursell87,Ursell91}, while Parker observed
experimentally the trapped modes in locally perturbed acoustic
waveguides \cite{Parker66,Parker67}.  Exner and Seba considered an
infinite bent strip of smooth curvature and showed the existence of
trapped modes by reducing the problem to Schr\"odinger operator in the
straight strip, with the potential depending on the curvature
\cite{Exner89}.  Goldstone and Jaffe gave the variational proof that the
wave equation subject to Dirichlet boundary condition always has a
localized eigenmode in an infinite tube of constant cross-section in
any dimension, provided that the tube is not exactly straight
\cite{Goldstone92}.  This result was further extended by Chenaud {\it
et al.} to arbitrary dimension \cite{Chenaud05}.  The problem of
localization in acoustic waveguides with Neumann boundary condition
has also been investigated \cite{Evans92,Evans94}.  For instance,
Evans {\it et al.}  considered a straight strip with an inclusion of
arbitrary (but symmetric) shape \cite{Evans94} (see \cite{Davies98}
for further extensions).  Such an inclusion obstructed the propagation
of waves and was shown to result in trapped modes.  The effect of
mixed Dirichlet, Neumann and Robin boundary conditions on the
localization was also investigated (see
\cite{Olendski10,Bulla97,Dittrich02,Freitas06} and references
therein).  A mathematical analysis of guided water waves was developed
by Bonnet-Ben Dhia and Joly \cite{Bonnet-Ben93} (see also
\cite{Bonnet-Ben99}).  Lower bounds for the eigenvalues below the
cut-off frequency (for which the associated eigenfunctions are
localized) were obtained by Ashbaugh and Exner for infinite thin tubes
in two and three dimensions \cite{Ashbaugh90}.  In addition, these
authors derived an upper bound for the number of the trapped modes.
Exner {\it et al.} considered the Laplacian in finite-length curved
tubes of arbitrary cross-section, subject to Dirichlet boundary
conditions on the cylindrical surface and Neumann conditions at the
ends of the tube.  They expressed a lower bound for the spectral
threshold of the Laplacian through the lowest eigenvalue of the
Dirichlet Laplacian in a torus determined by the geometry of the tube
\cite{Exner04}.  In a different work, Exner and co-worker investigated
bound states and scattering in quantum waveguides coupled laterally
through a boundary window \cite{Exner96}.

Examples of waveguides with {\it numerous} localized states were
reported in the literature.  For instance, Avishai {\it et al.}
demonstrated the existence of many localized states for a sharp
``broken strip'', i.e., a waveguide made of two channels of equal width
intersecting at a small angle $\theta$ \cite{Avishai91}.  Carini and
co-workers reported an experimental confirmation of this prediction
and its further extensions \cite{Carini93,Carini97,Londergan}.
Bulgakov {\it et al.}  considered two straight strips of the same
width which cross at an angle $\theta \in (0,\pi/2)$ and showed that,
for small $\theta$, the number of localized states is greater than $(1
- 2^{-2/3})^{3/2}/\theta$ \cite{Bulgakov02}.  Even for the simple case
of two strips crossed at the right angle $\theta = \pi/2$, Schult {\it
et al.} showed the existence of two localized states, one lying below
the cut-off frequency and the other being embedded into the continuous
spectrum \cite{Schult89}.

\subsection{Exponential estimate for eigenfunctions}
\label{sec:expon}

Qualitatively, an eigenmode is trapped when it cannot ``squeeze''
outside the localization region through narrow channels or branches of
the waveguide.  This happens when typical spatial variations of the
eigenmode, which are in the order of a wavelength $\pi
\lambda^{-1/2}$, are larger than the size $a$ of the narrow part,
i.e., $\pi \lambda^{-1/2} \geq a$ or $\lambda \leq \pi^2/a^2$
\cite{Jackson}.  This simplistic argument suggests that there exists a
threshold value $\mu$ (which may eventually be $0$), or so-called
cut-off frequency, such that the eigenmodes with $\lambda \leq \mu$
are localized.  Moreover, this qualitative geometrical interpretation
is well adapted for both unbounded and bounded domains.  While the
former case of infinite waveguides was thoroughly investigated, the
existence of trapped or localized eigenmodes in bounded domains has
attracted less attention.  Even the definition of localization in
bounded domains remains conventional because all eigenfunctions have
finite $L_2$ norm.

\begin{figure}
\begin{center}
\includegraphics[width=120mm]{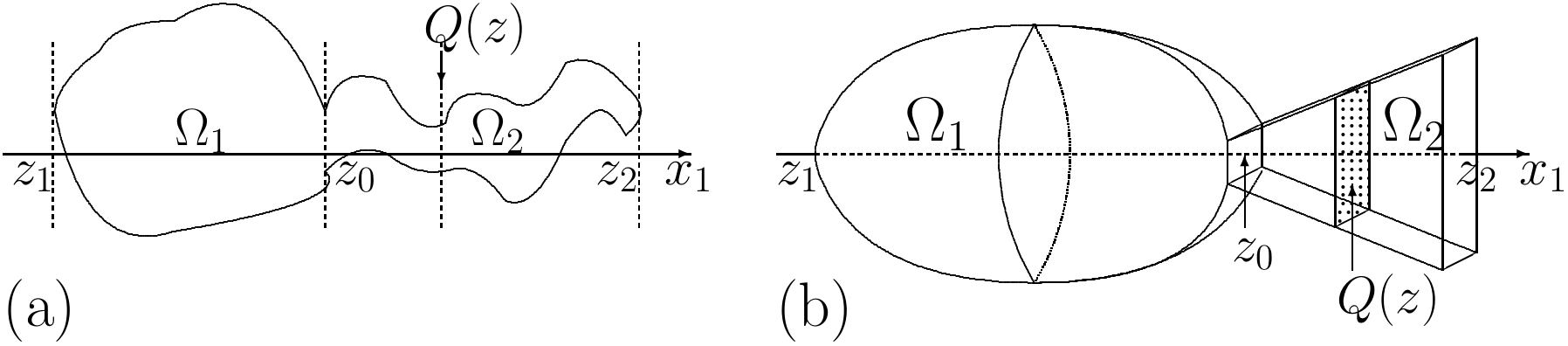}
\end{center}
\caption{
Two examples of a bounded domain $\Omega = \Omega_1 \cup \Omega_2$
with a branch of variable cross-sectional profile.  When the
eigenvalue $\lambda$ is smaller than the cut-off ``frequency'' $\mu$,
the associated eigenfunction exponentially decays in the branch
$\Omega_2$ and is thus mainly localized in $\Omega_1$.  Note that the
branch itself may even be increasing.}
\label{fig:branches}
\end{figure}

\begin{figure}
\begin{center}
\includegraphics[width=40mm]{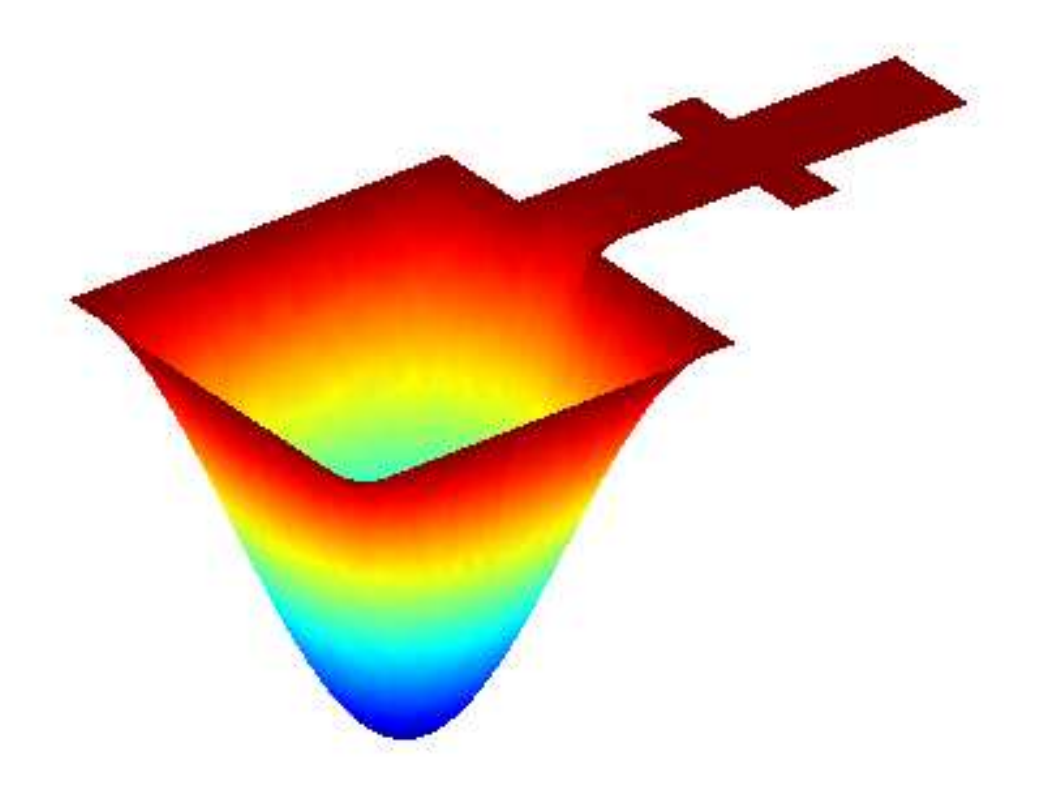}
\includegraphics[width=40mm]{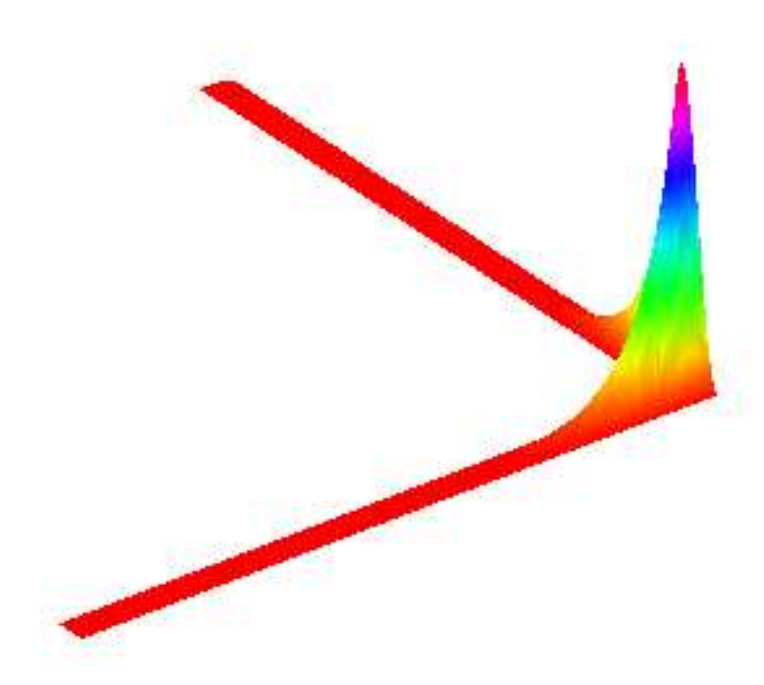}
\includegraphics[width=40mm]{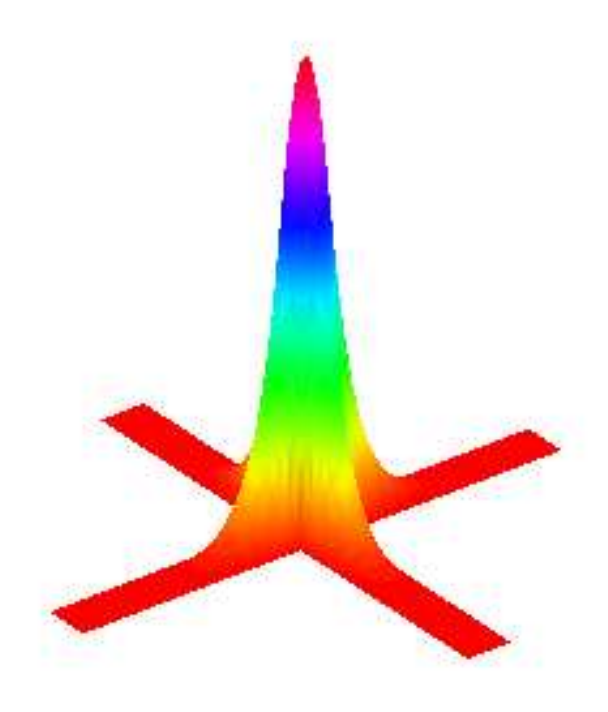}
\end{center}
\caption{
Examples of localized Dirichlet Laplacian eigenfunctions with an
exponential decay: square with a branch (from \cite{Delitsyn11a}),
L-shape and crossing of two stripes (from \cite{Delitsyn11b}). }
\label{fig:eigen_expon}
\end{figure}

This problem was studied by Delitsyn and co-workers for domains with
branches of variable cross-sectional profiles \cite{Delitsyn11a}.
More precisely, one considers a bounded domain $\Omega\subset \R^d$
($d = 2,3,...$) with a piecewise smooth boundary $\pa$ and denote
$Q(z) = \Omega\cap \{ x\in\R^d~:~ x_1 = z\}$ the cross-section of
$\Omega$ at $x_1 = z \in \R$ by a hyperplane perpendicular to the
coordinate axis $x_1$ (Fig. \ref{fig:branches}).  Let
\begin{equation*}
z_1 = \inf\{ z\in \R~:~ Q(z) \ne \emptyset\} , \qquad z_2 = \sup\{ z\in \R~:~ Q(z) \ne \emptyset\} ,
\end{equation*}
and we fix some $z_0$ such that $z_1 < z_0 < z_2$.  Let $\mu(z)$ be
the first eigenvalue of the Laplace operator in $Q(z)$, with Dirichlet
boundary condition on $\partial Q(z)$, and $\mu = \inf\limits_{z\in
(z_0,z_2)} \mu(z)$.  Let $u$ be a Dirichlet-Laplacian eigenfunction in
$\Omega$, and $\lambda$ the associate eigenvalue.  If $\lambda < \mu$,
then
\begin{equation}
\label{eq:expon_estimate}
\|u\|_{L_2(Q(z))} \leq \|u\|_{L_2(Q(z_0))} \exp(-\beta \sqrt{\mu-\lambda}~(z - z_0))   \quad  (z \geq z_0),
\end{equation} 
with $\beta = 1/\sqrt{2}$.  Moreover, if $(e_1 \cdot n(x)) \geq 0$ for
all $x\in\pa$ with $x_1 > z_0$, where $e_1$ is the unit vector
$(1,0,...,0)$ in the direction $x_1$, and $n(x)$ is the normal vector
at $x\in\pa$ directed outwards the domain, then the above inequality
holds with $\beta = 1$.

\begin{figure}
\begin{center}
\includegraphics[width=120mm]{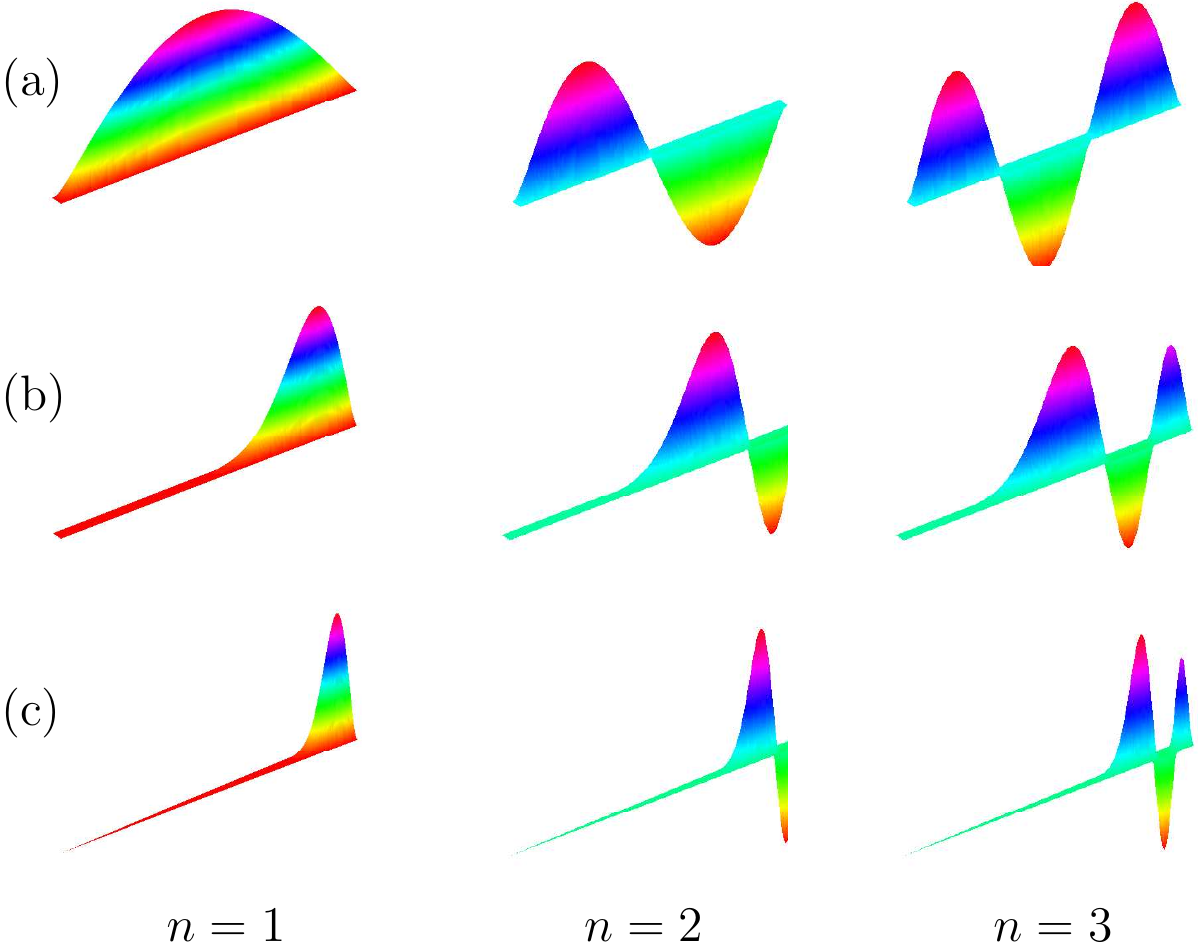}
\end{center}
\caption{
The first three Dirichlet Laplacian eigenfunctions for three elongated
domains: (a) rectangle of size $25\times 1$; (b) right trapezoid with
bases $1$ and $0.9$ and height $25$ which is very close to the above
rectangle; and (c) right triangle with edges $25$ and $1$ (half of the
rectangle).  There is no localization for the first shape, while the
first eigenfunctions for the second and third domains tend to be
localized. }
\label{fig:eigen_triangle}
\end{figure}

In this statement, a domain $\Omega$ is arbitrarily split into two
subdomains, $\Omega_1$ (with $x_1 < z_0$) and $\Omega_2$ (with $x_1 >
z_0$), by the hyperplane at $x_1 = z_0$ (the coordinate axis $x_1$ can
be replaced by any straight line).  Under the condition $\lambda <
\mu$, the eigenfunction $u$ exponentially decays in the subdomain
$\Omega_2$ which is loosely called ``branch''.  Note that the choice
of the splitting hyperplane (i.e., $z_0$) determines the threshold
$\mu$.

The theorem formalizes the notion of the cut-off frequency $\mu$ for
branches of variable cross-sectional profiles and provides a
constructive way for its computation.  For instance, if $\Omega_2$ is
a rectangular channel of width $a$, the first eigenvalue in all
cross-sections $Q(z)$ is $\pi^2/a^2$ (independently of $z$) so that
$\mu = \pi^2/a^2$, as expected.  The exponential estimate quantifies
the ``difficulty'' of penetration, or ``squeezing'', into the branch
$\Omega_2$ and ensures the localization of the eigenfunction $u$ in
$\Omega_1$.  Since the cut-off frequency $\mu$ is independent of the
subdomain $\Omega_1$, one can impose any boundary condition on
$\partial \Omega_1$ (that still ensures the self-adjointness of the
Laplace operator).  In turn, the Dirichlet boundary condition on the
boundary of the branch $\Omega_2$ is relevant, although some
extensions were discussed in
\cite{Delitsyn11a}.  It is worth noting that the theorem also applies
to infinite branches $\Omega_2$, under supplementary condition
$\mu(z)\to\infty$ to ensure the existence of the discrete spectrum.

According to this theorem, the $L_2$-norm of an eigenfunction with
$\lambda < \mu$ in $\Omega(z) = \Omega \cap \{ x\in\R^d ~:~ x_1 > z\}$
can be made exponentially small provided that the branch $\Omega_2$ is
long enough.  Taking $\Omega_0 = \Omega \setminus \Omega(z)$, the
ratio of $L_2$-norms in Eq. (\ref{eq:def_loc}) can be made arbitrarily
small.  However, the second ratio may not be necessarily small.  In
fact, its smallness depends on the shape of the domain $\Omega$.  This
is once again a manifestation of the conventional character of
localization in bounded domains.

Figure \ref{fig:eigen_expon} presents several examples of localized
Dirichlet Laplacian eigenfunctions showing an exponential decay along
the branches.  Since an increase of branches diminishes the eigenvalue
$\lambda$ and thus further enhances the localization, the area of the
localized region $\Omega_1$ can be made arbitrarily small with respect
to the total area (one can even consider infinite branches).  Examples
of an L-shape and a cross illustrate that the linear sizes of the
localized region do not need to be large in comparison with the branch
width (a sufficient condition for getting this kind of localization
was proposed in \cite{Delitsyn11b}).  It is worth noting that the
separation into the localized region and branches is arbitrary.  For
instance, Fig. \ref{fig:eigen_triangle} shows several localized
eigenfunctions for elongated triangle and trapezoid, for which there
is no explicit separation.

\begin{figure}
\begin{center}
\includegraphics[width=120mm]{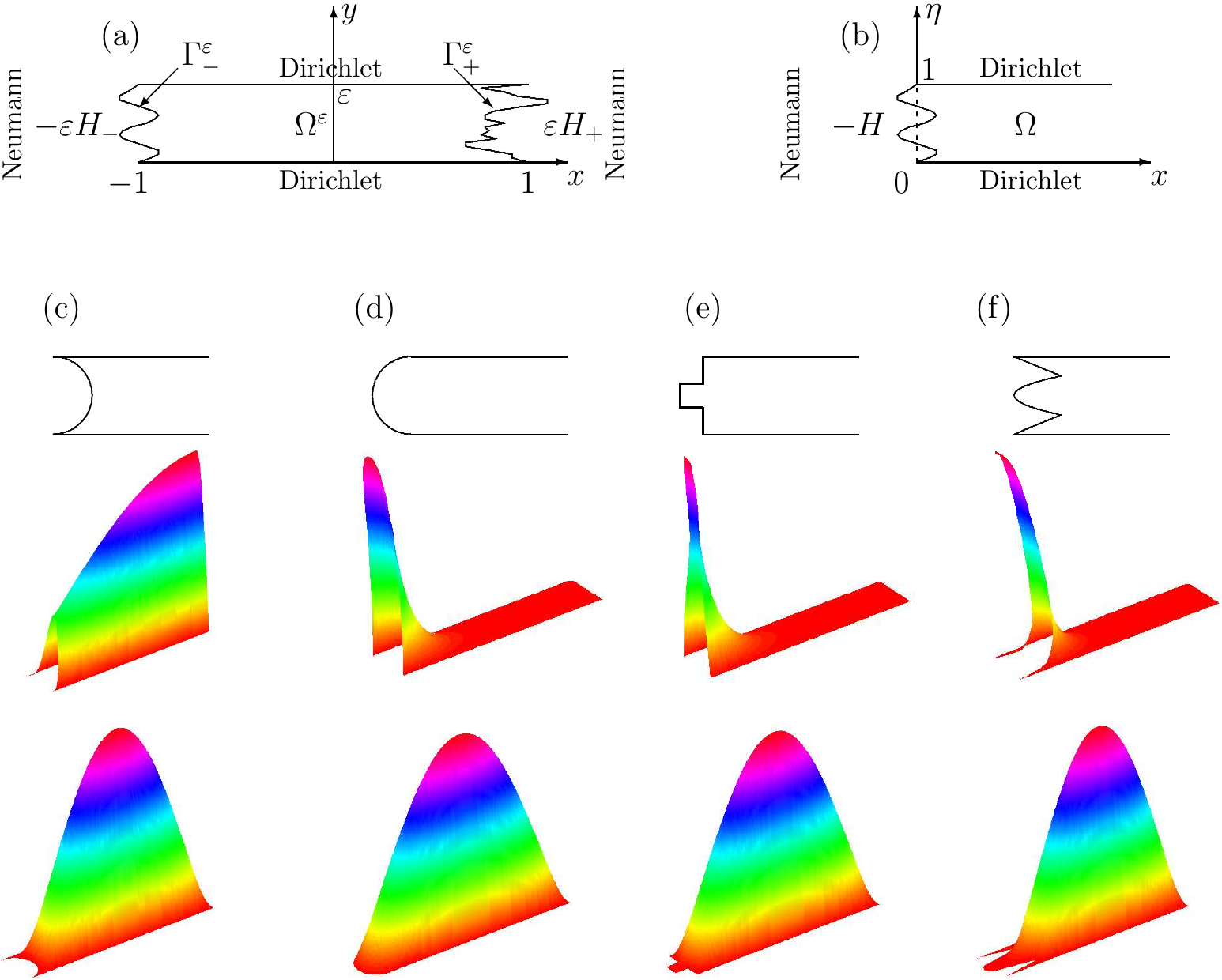}
\end{center}
\caption{
{\bf (a)} A thin cylinder $\Omega^\ve \in \R^2$ of width $\ve$ with
two distorted ends defined by functions $H_{\pm}(\eta)$,
$\eta\in\omega = (0,1)$.  {\bf (b)} In the limit $\ve\to 0$, the
analysis is reduced to a semi-infinite cylinder $\Omega$ with one
distorted end.  {\bf (c,d,e,f)} Four semi-infinite cylinders with
various distorted ends (top), the first eigenfunction for the Laplace
operator in these shapes with mixed Dirichlet-Neumann boundary
condition (middle), and the first eigenfunction for the Laplace
operator with purely Dirichlet boundary condition.  According to the
sufficient conditions (\ref{eq:Nazarov_cond1},\ref{eq:Nazarov_cond2}),
the first eigenfunction is localized near the distorted end for cases
'd', 'e' and 'f', and not localized for the case 'c'.  No localization
happens when the Dirichlet boundary condition is set over the whole
domain.  For numerical computation, semi-infinite cylinders were
``truncated'' and auxilary Dirichlet boundary condition was set at the
right straight end. }
\label{fig:nazarov}
\end{figure}

Localization and exponential decay of Laplacian eigenfunctions were
observed for various perturbations of cylindrical domains
\cite{Kamotskii00,Nazarov02,Cardone10}.  For instance, Kamotskii and
Nazarov studied localization of eigenfunctions in a neighborhood of
the lateral surface of a thin domain \cite{Kamotskii00}.  Nazarov and
co-workers analyzed the behavior of eigenfunctions for thin cylinders
with distorted ends \cite{Nazarov02,Cardone10}.  For a bounded domain
$\omega\subset \R^{n-1}$ ($n\geq 2$) with a simple closed Lipschitz
contour $\partial\omega$ and Lipschitz functions $H_\pm(\eta)$ in
$\bar{\omega} = \omega\cup \partial\omega$, the thin cylinder with
distorted ends is defined for a given small $\ve > 0$ as
\begin{equation*}
\Omega^\ve = \{ (x,y)\in \R\times \R^{n-1}~:~ y = \ve \eta, ~ -1- \ve H_-(\eta) < x < 1 + \ve H_+(\eta),~ \eta \in\omega\} .
\end{equation*}
One can view this domain as a thin cylinder $[-1,1]\times (\ve
\omega)$ to which two distorted ``cups'' characterized by functions
$H_\pm$, are attached (Fig. \ref{fig:nazarov}a).  The Neumann boundary
condition is imposed on the curved ends $\Gamma^\ve_\pm$:
\begin{equation*}
\Gamma^\ve_\pm = \{ (x,y)\in \R\times \R^{n-1}~:~ y = \ve \eta, ~ x = \pm 1 \pm \ve H_\pm(\eta),~ \eta \in\omega\} ,
\end{equation*}
while the Dirichlet boundary condition is set on the remaining lateral
side of the domain: $\Sigma^\ve = \partial \Omega^\ve \backslash
(\overline{\Gamma^\ve_+ \cup \Gamma^\ve_+})$.  When the ends of the
cylinder are straight ($H_\pm \equiv 0$), the eigenfunctions are
factored as $u_{mn}(x,y) = \cos(\pi m(x+1)/2) \varphi_n(y)$, where
$\varphi_n(y)$ are the eigenfunctions of the Laplace operator in the
cross-section $\omega$ with Dirichlet boundary condition.  These
eigenfunctions are extended over the whole cylinder, due to the cosine
factor.  Nazarov and co-workers showed that distortion of the ends
(i.e., $H_\pm \ne 0$) may lead to localization of the ground
eigenfunction in one (or both) ends, with an exponential decay toward
the central part.  In the limit $\ve\to 0$, the thinning of the
cylinder can be alternatively seen as its outstretching, allowing one
to reduce the analysis to a semi-infinite cylinder with one distorted
end (Fig. \ref{fig:nazarov}b), described by a single function
$H(\eta)$:
\begin{equation*}
\Omega = \{ (x,\eta)\in \R\times \R^{n-1}~:~ - H(\eta) < x ,~ \eta \in\omega\} .
\end{equation*}
Two sufficient conditions for getting the localized ground
eigenfunction at the distorted end were proposed in
\cite{Cardone10}:

(i) For $H\in C(\bar{\omega})$, the sufficient condition reads
\begin{equation}
\label{eq:Nazarov_cond1}
\int\limits_\omega d\eta ~ H(\eta) \biggl(|\nabla \varphi_1(\eta)|^2 - \mu_1 [\varphi_1(\eta)]^2\biggr) < 0 ,
\end{equation}
where $\mu_1$ is the smallest eigenvalue corresponding to $\varphi_1$
in the cross-section $\omega$ (in two dimensions, when $\omega =
[-1/2,1/2]$, one has $\varphi_1(\eta) = \sin \pi (\eta+1/2)$ and
$\mu_1 = \pi^2$ so that this condition reads $\int_{-1/2}^{1/2} d\eta
~ H(\eta) \cos(2\pi \eta) > 0$ \cite{Nazarov02});

(ii) Under a stronger assumption $H\in C^2(\bar{\omega})$, the
sufficient condition simplifies to
\begin{equation}
\label{eq:Nazarov_cond2}
\int\limits_\omega d\eta ~ [\varphi_1(\eta)]^2 \Delta H(\eta) < 0 .
\end{equation}
This inequality becomes true for a subharmonic profile $-H$ (i.e., for
$\Delta H(\eta) < 0$) with but false for superharmonic.%
\footnote{
In the discussion after Theorem 3 in \cite{Cardone10}, the sign minus
in front of $H$ was omitted.}

Figure \ref{fig:nazarov} shows several examples for which the
sufficient condition is either satisfied (\ref{fig:nazarov}d,e,f), or
not (\ref{fig:nazarov}c).  Nazarov and co-workers showed that these
results are applicable to bounded thin cylinders for small enough
$\ve$.  In addition, they found out a domain where the first
eigenfunction concentrates at the both ends simultaneously.  Finally,
they showed that no localization happens in the case in which the
mixed Dirichlet-Neumann boundary condition is replaced by the
Dirichlet boundary condition onto the whole boundary, as illustrated
on Fig. \ref{fig:nazarov} (see \cite{Cardone10} for further
discussions and results).

Friedlander and Solomyak studied the spectrum of the Dirichlet
Laplacian in a family of narrow strips of variable profile: $\Omega =
\{(x,y)\in\R^2~:~-a < x < b,~ 0 < y < \ve h(x)\}$
\cite{Friedlander09,Friedlander08}.  The main assumption was that $x
= 0$ is the only point of global maximum of the positive, continuous
function $h(x)$.  In the limit $\ve \to 0$, they found the two-term
asymptotics of the eigenvalues and the one-term asymptotics of the
corresponding eigenfunctions.  The asymptotic formulas obtained
involve the eigenvalues and eigenfunctions of an auxiliary ODE on $\R$
that depends only on the behavior of $h(x)$ as $x \to 0$, i.e., in the
vicinity of the widest cross-section of the strip.

\subsection{Dumbbell domains}
\label{sec:dumbbell}

\begin{figure}
\begin{center}
\includegraphics[width=120mm]{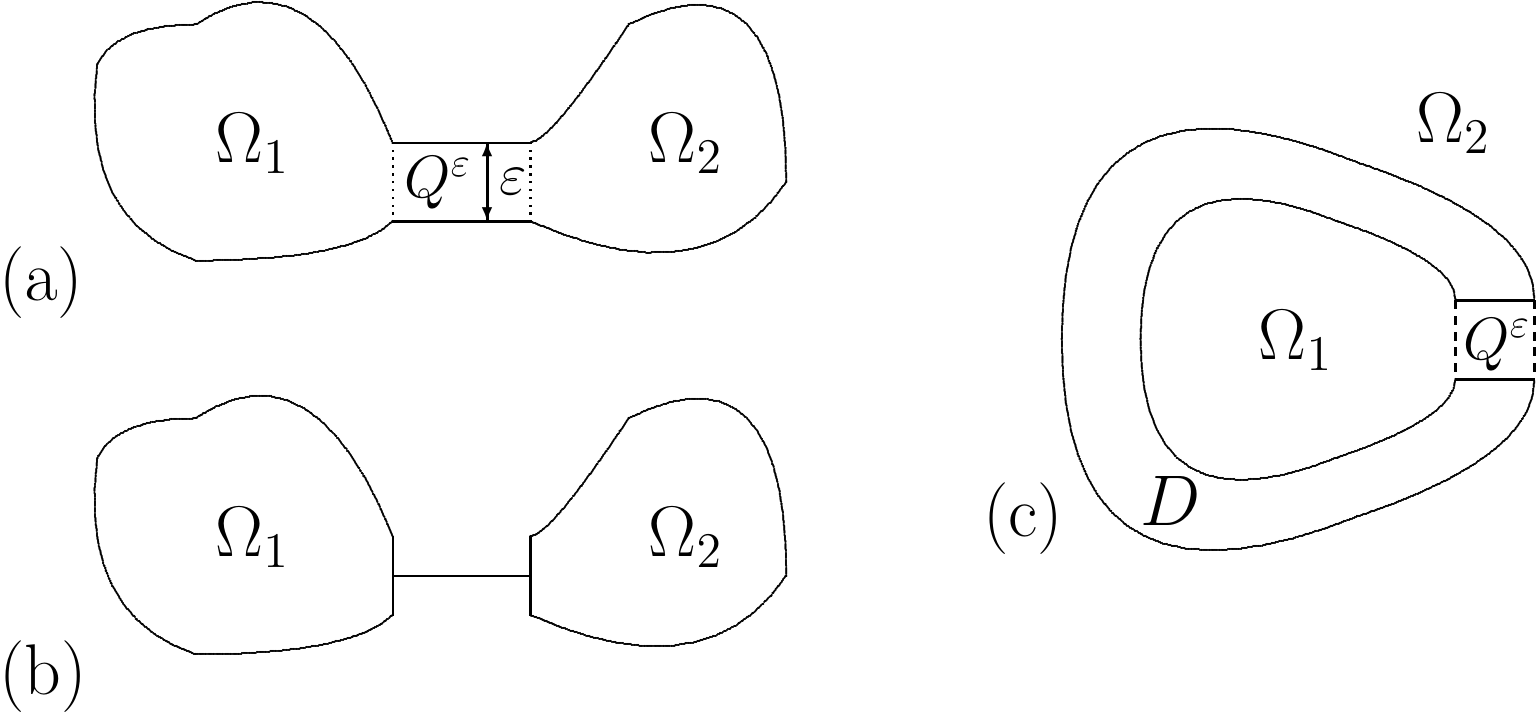}
\end{center}
\caption{
{\bf (a)} A dumbbell domain $\Omega^\ve$ is the union of two bounded
domains $\Omega_1$, $\Omega_2$ and a narrow ``connector'' $Q^\ve$ of
width $\ve$.  {\bf (b)} In the limit $\ve\to 0$, the connector
degenerates to a curve (here, an interval) so that the subdomains
$\Omega_1$ and $\Omega_2$ become disconnected.  {\bf (c)} In Beale's
work, $\Omega^\ve$ is composed of two components, a bounded domain
$\Omega_1$ and unbounded domain $\Omega_2$, which are connected by a
narrow channel $Q^\ve$. }
\label{fig:dumbbell}
\end{figure}

Yet another type of localization emerges for domains that can be split
into two or several subdomains with narrow connections (of ``width''
$\ve$) \cite{Raugel}, a standard example being a dumbbell: $\Omega^\ve
= \Omega_1 \cup Q^\ve \cup \Omega_2$ (Fig. \ref{fig:dumbbell}a).  The
asymptotic behavior of eigenvalues and eigenfunctions in the limit
$\ve\to 0$ was thoroughly investigated for both Dirichlet and Neumann
boundary conditions \cite{Jimbo09}.  We start by considering Dirichlet
boundary condition.

\begin{figure}
\begin{center}
\includegraphics[width=125mm]{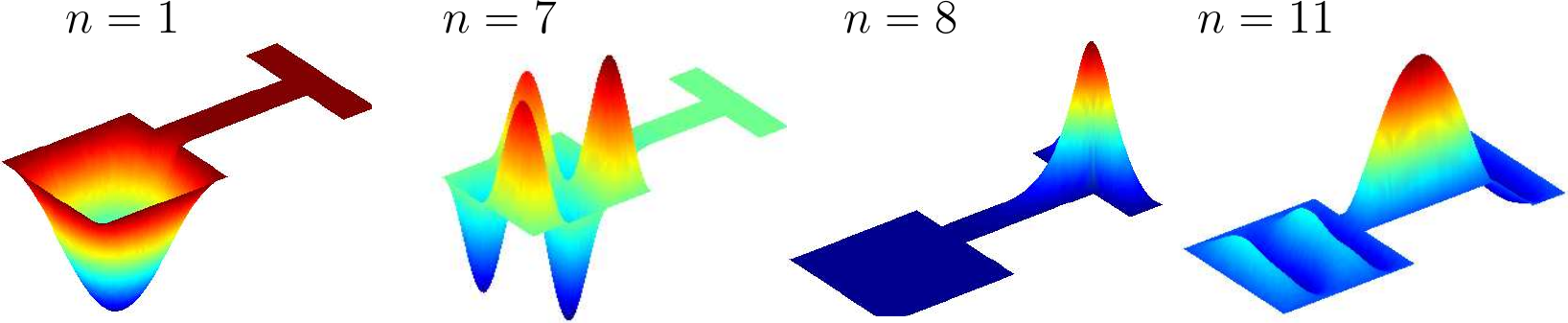}
\end{center}
\caption{
Several Dirichlet Laplacian eigenfunctions for a dumbbell domain which
is composed of two rectangles and connected by the third rectangle
(from \cite{Delitsyn11a}).  The 1st and 7th eigenfunctions are
localized in the larger subdomain, the 8th eigenfunction is localized
in the smaller subdomain, while the 11th eigenfunction is not
localized at all.  Note that the width of connection is not small
($1/4$ of the width of both subdomains). }
\label{fig:eigen_dumbbell}
\end{figure}

In the limiting case of zero width connections, the subdomains
$\Omega_i$ ($i=1,...,N$) become disconnected, and the eigenvalue
problem can be independently formulated for each subdomain.  Let
$\Lambda_i$ be the set of eigenvalues for the subdomain $\Omega_i$.
Each Dirichlet eigenvalue $\lambda^\ve$ of the Laplace operator in the
domain $\Omega^\ve$ approaches to an eigenvalue $\lambda^0$
corresponding to one limiting subdomain $\Omega_i \subset \Omega^0$:
$\lambda^0\in \Lambda_i$ for certain $i$.  Moreover, if
\begin{equation}
\label{eq:cond_bottle}
\Lambda_i \cap \Lambda_j = \emptyset  \quad \forall~ i\ne j, 
\end{equation}
the space of eigenfunctions in the limiting (disconnected) domain
$\Omega^0$ is the direct product of spaces of eigenfunctions for each
subdomain $\Omega_i$ (see \cite{Daners03} for discussion on
convergence and related issues).  This is a basis for what we will
call ``bottle-neck localization''.  In fact, each eigenfunction
$u_m^\ve$ on the domain $\Omega^\ve$ approaches an eigenfunction
$u_m^{0}$ of the limiting domain $\Omega^0$ which is fully localized
in one subdomain $\Omega_i$ and zero in the others.  For a small
$\ve$, the eigenfunction $u_m^\ve$ is therefore mainly localized in
the corresponding $i$-th subdomain $\Omega_i$, and is {\it almost}
zero in the other subdomains.  In other words, for any eigenfunction,
one can take the width $\ve$ small enough to ensure that the
$L_2$-norm of the eigenfunction in the subdomain $\Omega_i$ is
arbitrarily close to that in the whole domain $\Omega^\ve$:
\begin{equation}
\label{eq:ineq_dumbbel}
\forall~ m \geq 1  \quad \exists i\in \{1,...,N\}  \quad   \forall~ \delta \in (0,1)
\quad   \exists ~\ve > 0 ~:~  \|u_m^\ve\|_{L_2(\Omega_i)} > (1-\delta) \|u_m^\ve\|_{L_2(\Omega^\ve)} .
\end{equation}
This behavior is exemplified for a dumbbell domain which is composed
of two rectangles and connected by the third rectangle
(Fig. \ref{fig:eigen_dumbbell}).  The 1st and 7th eigenfunctions are
localized in the larger rectangle, the 8th eigenfunction is localized
in the smaller rectangle, while the 11th eigenfunction is not
localized at all.  Note that the width of connection is not too small
($1/4$ of the width of both subdomains).

It is worth noting that, for a small fixed width $\ve$ and a small
fixed threshold $\delta$, there may be infinitely many high-frequency
``non-localized'' eigenfunctions, for which the inequality
(\ref{eq:ineq_dumbbel}) is not satisfied.  In other words, for a given
connected domain with a narrow connection, one can only expect to
observe a finite number of low-frequency localized eigenfunctions.  We
note that the condition (\ref{eq:cond_bottle}) is important to ensure
that limiting eigenfunctions are fully localized in their respective
subdomains.  Without this condition, a limiting eigenfunction may be a
linear combination of eigenfunctions in different subdomains with the
same eigenvalue that would destroy localization.  Note that the
asymptotic behavior of eigenfunctions at the ``junction'' was studied
by Felli and Terracini \cite{Felli12}.

For Neumann boundary condition, the situation is more complicated, as
the eigenvalues and eigenfunctions may also approach the eigenvalues
and eigenfunctions of the limiting connector (in the simplest case,
the interval).  Arrieta considered a planar dumbbell domain
$\Omega_\ve$ consisted of two disjoint domains $\Omega_1$ and
$\Omega_2$ connected by a channel $Q^\ve$ of variable profile $g(x)$:
$Q^\ve = \{{ x\in\R^2~:~ 0 < x_1 < 1,~ 0 < x_2 < \ve g(x_1)}\}$, where
$g \in C^1(0,1)$ and $g(x_1) > 0$ for all $x_1 \in [0,1]$.  In the
limit $\ve\to 0$, each eigenvalue of the Laplace operator in
$\Omega^\ve$ with Neumann boundary condition was shown to converge
either to an eigenvalue $\mu_k$ of the Neumann-Laplace operator in
$\Omega_1\cup \Omega_2$, or to an eigenvalue $\nu_k$ of the
Sturm-Liouville operator $-\frac{1}{g}(gu_x)_x$ acting on a function
$u$ on $(0,1)$, with Dirichlet boundary condition
\cite{Arrieta95,Arrieta95b}.  The first-order term in the small
$\ve$-asymptotic expansion was obtained.  The special case of
cylindrical channels (of constant profile) in higher dimensions was
studied by Jimbo \cite{Jimbo89} (see also results by Hempel {\it et
al.} \cite{Hempel91}).  Jimbo and Morita studied an $N$-dumbbell
domain, i.e., a family of $N$ pairwise disjoint domains connected by
thin channels \cite{Jimbo92}.  They proved that $\lambda_m^{\ve} = C_m
\ve^{d-1} + o(\ve^{d-1})$ as $\ve \to 0$ for $m=1, 2, \dots, N$, while
$\lambda^{\ve}_{N+1}$ is uniformly bounded away from zero, where $d$
is the dimension of the embedding space, and $C_m$ are shape-dependent
constants.  Jimbo also analyzed the asymptotic behavior of the
eigenvalues $\lambda^{\ve}_m$ with $m > N$ under the condition that
the sets $\{\mu_k\}$ and $\{\nu_k\}$ do not intersect \cite{Jimbo93}.
In particular, for an eigenvalue $\lambda_m^{\ve}$ that converges to
an element of $\{\mu_k\}$, the asymptotic behavior is $\lambda_m^{\ve}
= \mu_k + C_m \ve^{d-1}+ o(\ve^{d-1})$.

Brown and co-workers studied upper bounds for $|\lambda_m^{\ve} -
\lambda_m^0|$ and showed \cite{Brown95}:

(i) If $\lambda_m^0 \in \{\mu_k\} \setminus \{\nu_k\}$,
\begin{eqnarray*}
|\lambda^{\ve}_m - \lambda_m^0| &\le& C |\ln \ve|^{-1/2} \hskip 5mm  (d=2),\\
|\lambda^{\ve}_m - \lambda_m^0| &\le& C \ve^{(d-2)/d} \qquad (d\ge 3). 
\end{eqnarray*}

(ii) If $\lambda_m^0 \in \{\nu_k\} \setminus \{\mu_k\}$, 
\begin{eqnarray*}
|\lambda^{\ve}_m - \lambda_m^0| &\le& C \ve^{1/2} |\ln \ve| \qquad (d=2),\\
|\lambda^{\ve}_m - \lambda_m^0| &\le& C \ve^{1/2}  \hskip 15mm (d\ge 3). 
\end{eqnarray*}

For a dumbbell domain in $\R^d$ with a thin cylindrical channel of a
smooth profile, Gadyl'shin obtained the complete small $\ve$
asymptotics of the Neumann-Laplace eigenvalues and eigenfunctions and
explicit formulas for the first term of these asymptotics, including
multiplicities \cite{Gadylshin93,Gadylshin94,Gadylshin05}.  

Arrieta and Krej\v{c}i\v{r}\'ik considered the problem of spectral
convergence from another point of view \cite{Arrieta10}.  They showed
that if $\Omega_0 \subset \Omega_\ve$ are bounded domains and if the
eigenvalues and eigenfunctions of the Laplace operator with Neumann
boundary condition in $\Omega_\ve$ converge to the ones in $\Omega_0$,
then necessarily $\mu_d(\Omega_\ve \backslash \Omega_0)\to 0$ as $\ve
\to 0$, while it is not necessarily true that ${\rm dist}(\Omega_\ve,
\Omega_0)\to 0$.  As a matter of fact, they constructed an example of
a perturbation where the spectra behave continuously but ${\rm
dist}(\Omega_\ve, \Omega_0)\to \infty$ as $\ve
\to 0$.

A somewhat related problem of scattering frequencies of the wave
equation associated to an exterior domain in $\R^3$ with an
appropriate boundary condition was investigated by Beale
\cite{Beale73} (for more general aspects of geometric scattering
theory, see \cite{Melrose}).  We recall that a scattering frequency
$\sqrt{\lambda}$ of an unbounded domain $\Omega$ is a (complex) number
for which there exists a nontrivial solution of $\Delta u+\lambda u =
0$ in $\Omega$, subject to Dirichlet, Neumann, or Robin boundary
condition and to an ``outgoing'' condition at infinity.  In Beale's
work, a bounded cavity $\Omega_1$ was connected by a thin channel to
the exterior (unbounded) space $\Omega_2$.  More specifically, he
considered a bounded domain $D$ such that its complement in $\R^3$ has
a bounded component $\Omega_1$ and an unbounded component $\Omega_2$.
After that, a thin ``hole'' $Q^\ve$ in $D$ was made to connect both
components (Fig. \ref{fig:dumbbell}c).  Beale showed that the joint
domain $\Omega^\ve = \Omega_1\cup Q^\ve\cup \Omega_2$ with Dirichlet
boundary condition has a scattering frequency which is arbitrarily
close either to an eigenfrequency (i.e., the square root of the
eigenvalue) of the Laplace operator in $\Omega_1$, or to a scattering
frequency in $\Omega_2$, provided the channel $Q^\ve$ is narrow
enough.  The same result was extended to Robin boundary condition of
the form $\partial u/\partial n + h u = 0$ on $\partial \Omega^\ve$,
where $h$ is a function on $\partial
\Omega^\ve$ with a positive lower bound.  In both cases, the method in
his proof relies on the fact that the lowest eigenvalue of the channel
tends to infinity as the channel narrows.  However, it is no longer
true for Neumann boundary condition.  In this case, with some
restrictions on the shape of the channel, Beale proved that the
scattering frequencies converge not only to the eigenfrequencies of
$\Omega_1$ and scattering frequencies of $\Omega_2$ but also to the
longitudinal frequencies of the channel.  Similar results can be
obtained in domains of space dimension other than $3$.

There are other problems for partial differential equations in
dumbbell domains which undergo a singular perturbation
\cite{Hale88,Henry,Mazya00}.  For instance, in a series of articles,
Arrieta {\it et al.} studied the behavior of the asymptotic nonlinear
dynamics of a reaction-diffusion equation with Neumann boundary
condition \cite{Arrieta06,Arrieta09,Arrieta09b}.  In this context,
dumbbell domains appear naturally as the counterpart of convex domains
for which the stable stationary solutions to a reaction-diffusion
equation are necessarily spatially constant \cite{Casten78}.  As
explained in \cite{Arrieta06}, one way to produce ``patterns'' (i.e.,
stable stationary solutions which are not spatially constant), is to
consider domains which make it difficult to diffuse from one part of
the domain to the other, making a constriction in the domain.  Kosugi
studied the convergence of the solution to a semilinear elliptic
equation in a thin network-shaped domain which degenerates into a
geometric graph when a certain parameter tends to zero \cite{Kosugi00}
(see also \cite{Kuchment03,Rubinstein01,Post05,Rubinstein06}).


\subsection{Localization in irregularly-shaped domains}
\label{sec:irregular}

\begin{figure}
\begin{center}
\includegraphics[width=120mm]{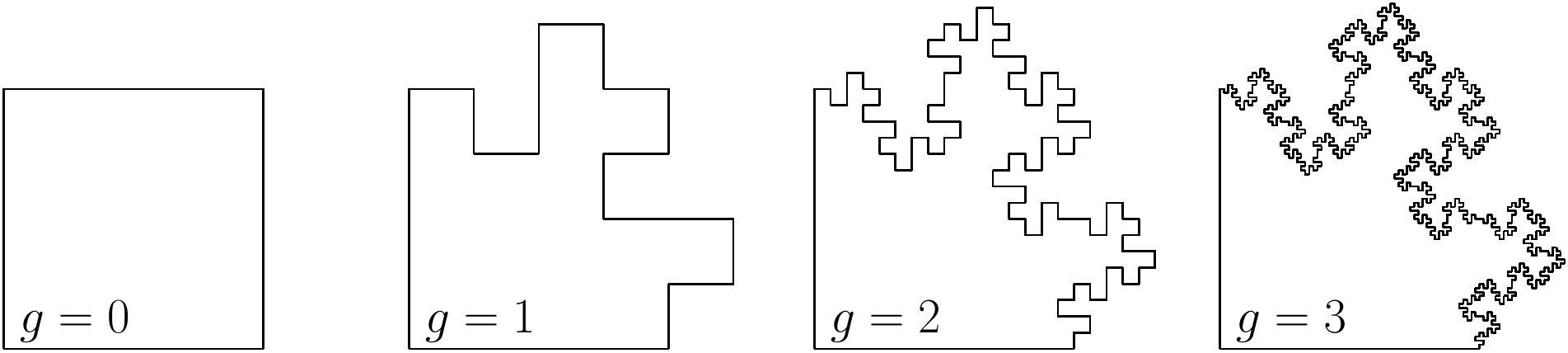}
\end{center}
\caption{
The unit square and three prefractal domains obtained iteratively one
from the other (two sides of these domains are finite generations of
the Von Koch curve of fractal dimension $3/2$).  These domains were
intensively studied, both numerically and experimentally, by Sapoval
and co-workers
\cite{Sapoval91,Sapoval93,Sapoval97,Russ97,Haeberle98,Hebert99,Even99,Russ02}.
}
\label{fig:fractal}
\end{figure}

\begin{figure}
\begin{center}
\includegraphics[width=125mm]{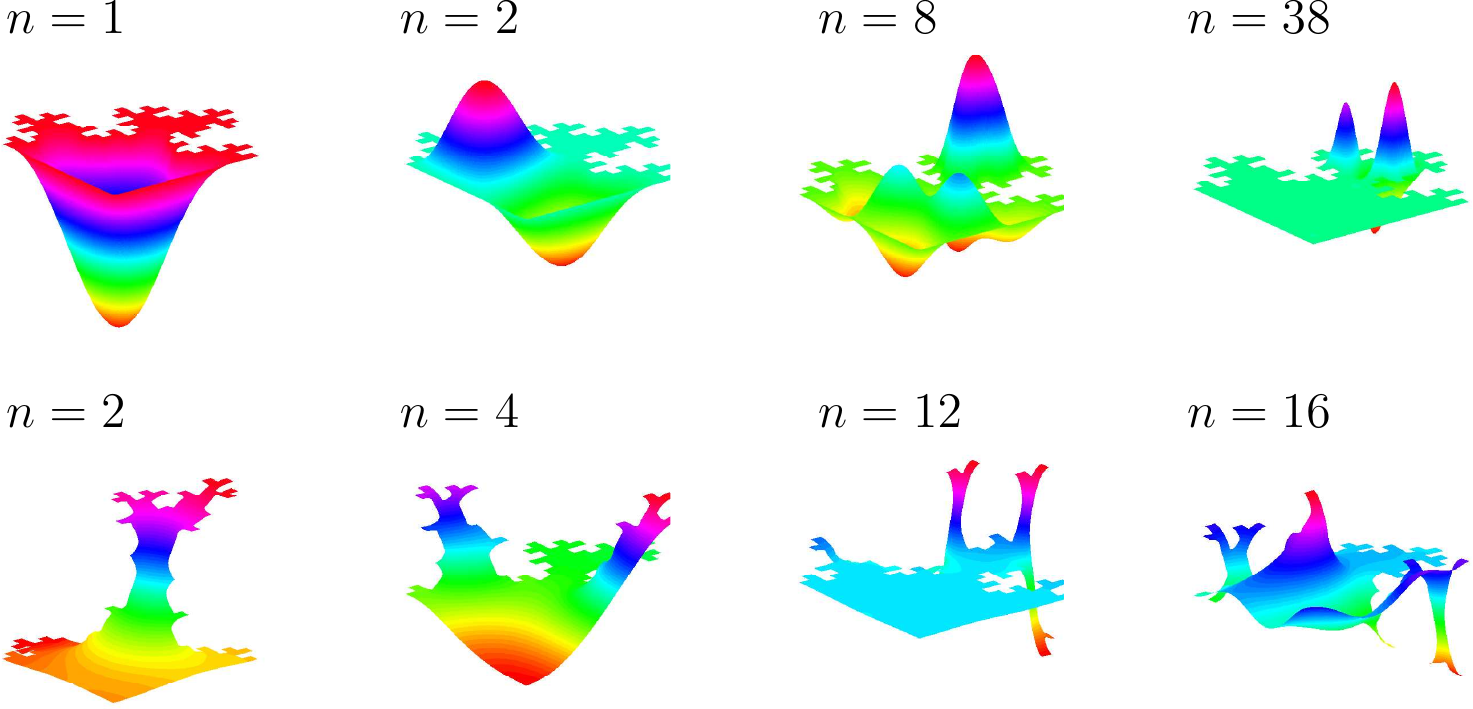}
\end{center}
\caption{
Several Dirichlet (top) and Neumann (bottom) eigenfunctions for the
third domain on Fig. \ref{fig:fractal} ($g = 2$).  The $38$th
Dirichlet and the $12$th Neumann eigenfunctions are localized in a
small subdomain (located in the upper right corner on
Fig. \ref{fig:fractal}), while the first/second Dirichlet and the
$4$th Neumann eigenfunctions are almost zero on this subdomain.
Finally, the $8$th Dirichlet and the second Neumann eigenfunctions are
examples of eigenfunctions extended over the whole domain. }
\label{fig:eigen_fractal}
\end{figure}

\begin{figure}
\begin{center}
\includegraphics[width=120mm]{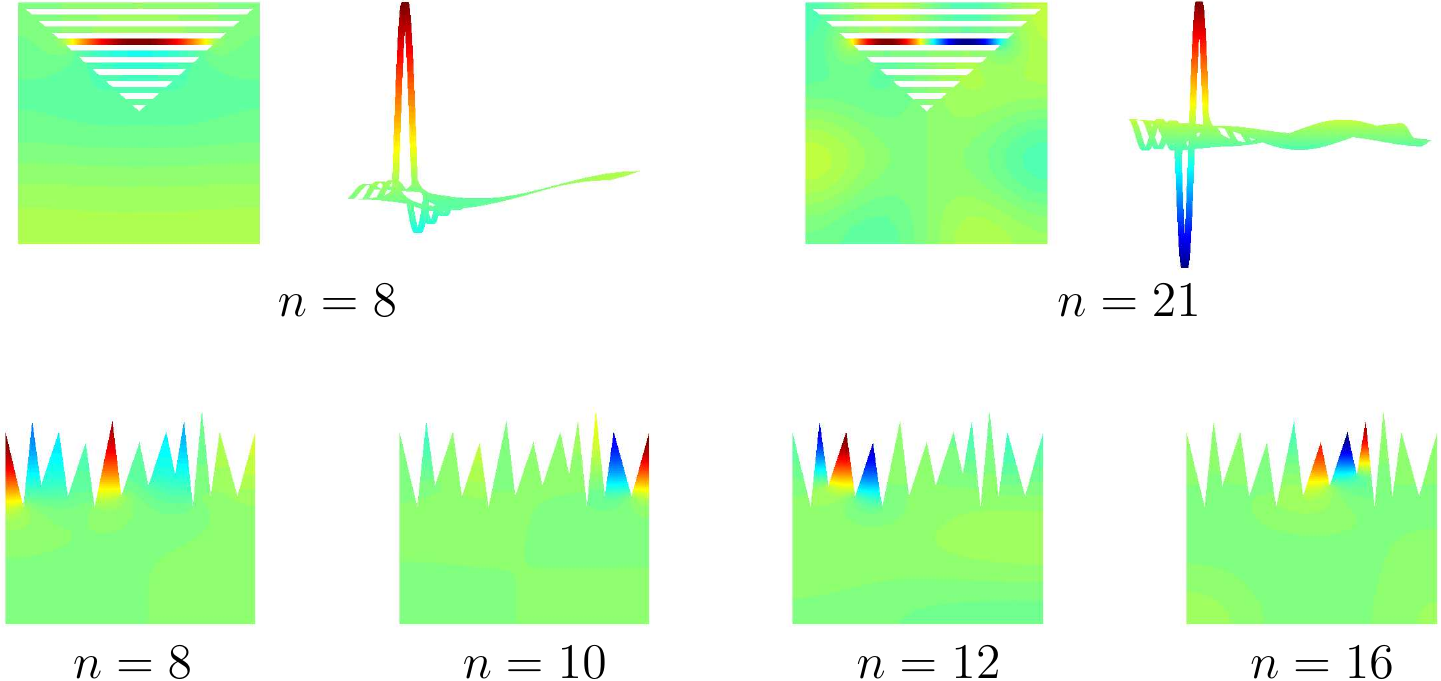}
\end{center}
\caption{
Examples of localized Neumann Laplacian eigenfunctions in two domains
adapted from \cite{Felix07}: square with many elongated holes (top)
and random sawteeth (bottom).  Colors represent the amplitude of
eigenfunctions, from the most negative value (dark blue), through zero
(green), to the largest positive value (dark red).  One can notice
that the eigenfunctions on the top are not negligible outside the
localization region.  This is yet another illustration for the
conventional character of localization in bounded domains. }
\label{fig:felix}
\end{figure}

As we have seen, a narrow connection between subdomains could lead to
localization.  How narrow should it be?  A rigorous answer to this
question is only known for several ``tractable'' cases such as
dumbbell-like or cylindrical domains (Sec. \ref{sec:dumbbell}).
%
%
Sapoval and co-workers have formulated and studied the problem of
localization in irregularly-shaped or fractal domains through
numerical simulations and experiments
\cite{Sapoval91,Sapoval93,Sapoval97,Russ97,Haeberle98,Hebert99,Even99,Russ02,Felix07}.
In the first publication, they monitored the vibrations of a
prefractal ``drum'' (i.e., a thin membrane with a fixed boundary)
which was excited at different frequencies \cite{Sapoval91}.  Tuning
the frequency allowed them to directly visualize different Dirichlet
Laplacian eigenfunctions in a (prefractal) quadratic von Koch
snowflake (an example is shown on Fig. \ref{fig:fractal}).  For this
and similar domains, certain eigenfunctions were found to be localized
in a small region of the domain, for both Dirichlet and Neumann
boundary conditions (Fig. \ref{fig:eigen_fractal}).  This effect was
first attributed to self-similar structure of the domain.  However,
similar effects were later observed through numerical simulations for
non-fractal domains \cite{Felix07,Sapoval08}, as illustrated by
Fig. \ref{fig:felix}.  In the study of sound attenuation by
noise-protective walls, F\'elix and co-workers have further extended
the analysis to the union of two domains with different refraction
indices which are separated by an irregular boundary
\cite{Felix07,Sapoval08,Felix09}.  Many eigenfunctions of the related
second order elliptic operator were shown to be localized on this
boundary (so-called ``astride localization'').  A rigorous
mathematical theory of these important phenomena is still missing.
Takeda {\it al.}  observed experimentally the electromagnetic field at
specific frequency to be confined in the central part of the third
stage of three-dimensional fractals called the Menger sponge
\cite{Takeda04}.  This localization was attributed to a singular
photon density of states realized in the fractal structure.

A number of mathematical studies were devoted to the theory of partial
differential equations on fractals in general and to localization of
Laplacian eigenfunctions in particular (see \cite{Kigami,Strichartz}
and references therein).  For instance, the spectral properties of the
Laplace operator on Sierpinski gasket and its extensions were
thoroughly investigated
\cite{Barlow88,Barlow89,Shima91,Fukushima92,Barlow95,Barlow99,Blasiak08}.
Barlow and Kigami studied the localized eigenfunctions of the
Laplacian on a more general class of self-similar sets (so-called post
critically finite self-similar sets, see \cite{Kigami93,Kigami93b} for
details).  They related the asymptotic behavior of the eigenvalue
counting function to the existence of localized eigenfunctions and
established a number of sufficient conditions for the existence of a
localized eigenfunction in terms of the symmetries of a set
\cite{Barlow97}.

Berry and co-workers developed a new method to approximate the Neumann
spectrum of the Laplacian on a planar fractal set $\Omega$ as a
renormalized limit of the Neumann spectra of the standard Laplacian on
a sequence of domains that approximate $\Omega$ from the outside
\cite{Berry09}.  They applied this method to compute the
Neumann-Laplacian eigenfunctions in several domains, including a
sawtooth domain, Sierpinski gasket and carpet, as well as nonsymmetric
and random carpets and the octagasket.  In particular, they gave a
numerical evidence for the localized eigenfunctions for a sawtooth
domain, in agreement with the earlier work by F\'elix {\it et al.}
\cite{Felix07}.

\begin{figure}
\begin{center}
\includegraphics[width=125mm]{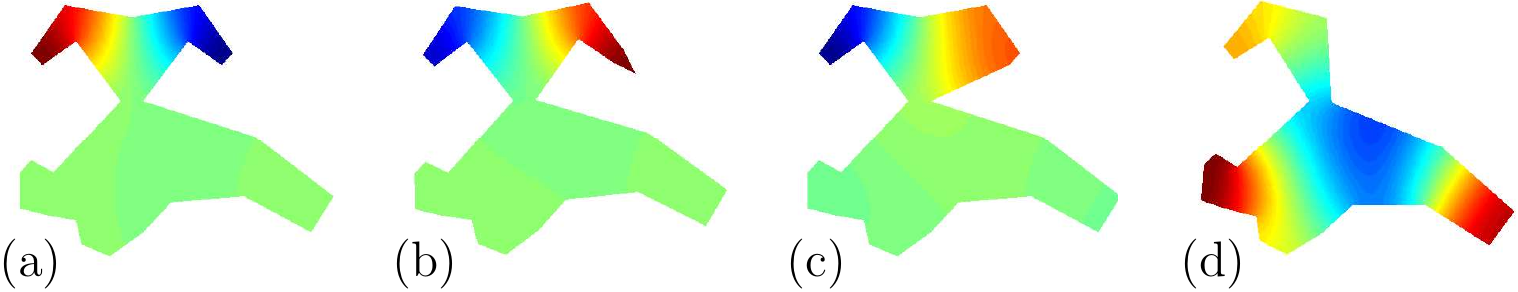}
\end{center}
\caption{
Neumann-Laplace eigenfunction $u_4$ in the original ``cow'' domain
from \cite{Heilman10} (a) and in three modified domains (b,c,d), in
which the reflection symmetry of the upper subdomain is broken.  The
fourth eigenfunction is localized for the first three domains (a,b,c),
while the last domain with the strongest modification shows no
localization (d).  Colors represent the amplitude of an eigenfunction,
from the most negative value (dark blue), through zero (green), to the
largest positive value (dark red). }
\label{fig:heilman}
\end{figure}

Heilman and Strichartz reported several numerical examples of
localized Neumann-Laplacian eigenfunctions in two domains
\cite{Heilman10}, one of them is illustrated on
Fig. \ref{fig:heilman}a.  Each of these domains consists of two
subdomains with a narrow but not too narrow connection.  This is a
kind of dumbbell shape with a connector of zero length.  Heilman and
Strichartz argued that one subdomain must possess an axis of symmetry
for getting localized eigenfunctions.  Since an anti-symmetric
eigenfunction vanishes on the axis of symmetry, it is necessarily
small near the bottle-neck that somehow ``prevents'' its extension to
the other domain.  Although the argument is plausible, we have to
stress that such a symmetry is neither sufficient, nor necessary for
localization.  It is obviously not sufficient because even for
symmetric domain, there exist plenty of extended eigenfunctions
(including the trivial example of the ground eigenmode which is a
constant over the whole domain).  In order to illustrate that the
reflection symmetry is not necessary, we plot on
Fig. \ref{fig:heilman}b,c examples of localized eigenfunctions for
modified domains for which the symmetry is broken.  Although rendering
the upper domain less and less symmetric gradually reduces or even
fully destroys localization (Fig. \ref{fig:heilman}d), its
``mechanism'' remains poorly understood.  We also note that methods of
Sec. \ref{sec:expon} are not applicable in this case because of
Neumann boundary condition.

Lapidus and Pang studied the boundary behavior of the Dirichlet
Laplacian eigenfunctions and their gradients on a class of planar
domains with fractal boundary, including the triangular and square von
Koch snowflakes and their polygonal approximations \cite{Lapidus95}.
A numerical evidence for the boundary behavior of eigenfunctions was
reported in \cite{Lapidus96b}, with numerous pictures of
eigenfunctions.  Later, Daubert and Lapidus considered more
specifically the localization character of eigenfunctions in von Koch
domains \cite{Daudert07}.  In particular, different ``measures'' of
localization were discussed.

Note also that Filoche and Mayboroda studied the problem of
localization for bi-Laplacian in rigid thin plates and discovered that
clamping just one point inside such a plate not only perturbs its
spectral properties, but essentially divides the plate into two
independently vibrating regions \cite{Filoche09}.

\subsection{High-frequency localization}
\label{sec:high-freq}

A hundred years ago, Lord Rayleigh documented an interesting
acoustical phenomenon in the whispering gallery under the dome of
Saint Paul's Cathedral in London \cite{Rayleigh1910} (see also
\cite{Raman1921,Raman1922}).  A whisper of one person propagated along
the curved wall to another person stood near the wall.  Keller and
Rubinow discussed the related ``whispering gallery modes'' and also
``bouncing ball modes'', and showed that these modes exist for a
two-dimensional domain with arbitrary smooth convex curve as its
boundary \cite{Keller60}.  A semiclassical approximation of Laplacian
eigenfunctions in convex domains was developed by Lazutkin
\cite{Babich68,Lazutkin68,Lazutkin73,Lazutkin81,Lazutkin93} (see also
\cite{Arnold72,Smith74,Ralston76,Ralston77}).  Chen and co-workers
analyzed Mathieu and modified Mathieu functions and reported another
type of localization named ``focusing modes'' \cite{Chen94}.  All
these eigenmodes become more and more localized in a small subdomain
when the associated eigenvalue increases.  This so-called
high-frequency or high-energy limit was intensively studied for
various domains, named as quantum billiards
\cite{Gutzwiller,Heller98b,Stockmann,Jakobson01,Sarnak11}.  In quantum
mechanics, this limit is known as semi-classical approximation
\cite{Berry72}.  In optics, it corresponds to ray approximation of
wave propagation, from which the properties of an optical, acoustical
or quantum system can often be reduced to the study of rays in
classical billiards.  Jakobson {\it et al.} gave an overview of many
results on geometric properties of the Laplacian eigenfunctions on
Riemannian manifolds, with a special emphasis on high-frequency limit
(weak star limits, the rate of growth of $L_p$ norms, relationships
between positive and negative parts of eigenfunctions, etc.)
\cite{Jakobson01} (see also \cite{Uhlenbeck76,Albert78}).
Bearing in mind this comprehensive review, we start by illustrating
the high-frequency localization and the related problems in simple
domains such as disk, ellipse and rectangle for which explicit
estimates are available.  After that, some results for quantum
billiards are summarized.

\subsubsection{Whispering gallery and focusing modes}

The disk is the simplest shape for illustrating the whispering gallery
and focusing modes.  The explicit form (\ref{eq:u_disk}) of
eigenfunctions allows one to get accurate estimates and bounds, as
shown below.  When the index $k$ is fixed and $n$ increases, the
Bessel functions $J_n(\alpha_{nk} r/R)$ become strongly attenuated
near the origin (as $J_n(z)\sim (z/2)^n/n!$ at small $z$) and
essentially localized near the boundary, yielding whispering gallery
modes.  In turn, when $n$ is fixed and $k$ increases, the Bessel
functions rapidly oscillate, the amplitude of oscillations decreasing
towards the boundary.  In that case, the eigenfunctions are mainly
localized at the origin, yielding focusing modes.

\begin{figure}
\begin{center}
\includegraphics[width=120mm]{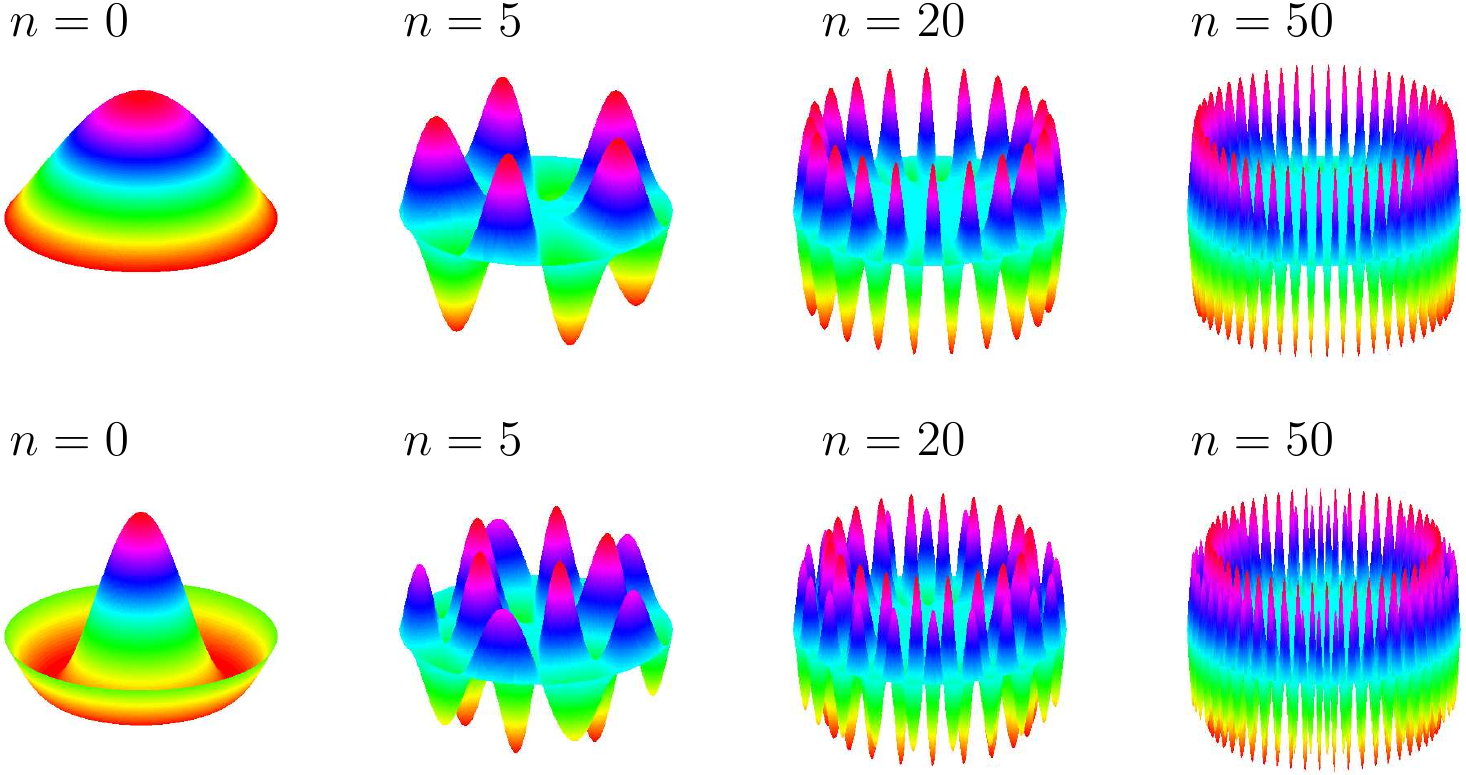}
\end{center}
\caption{
Formation of whispering gallery modes for the unit disk with Dirichlet
boundary condition: for a fixed $k$ ($k = 1$ for top figures and $k =
2$ for bottom figures), an increase of the index $n$ leads to stronger
localization of eigenfunctions near the boundary. }
\label{fig:whispering}
\end{figure}

These qualitative arguments were rigorously formulated in
\cite{Nguyen12}.  For each eigenfunction $u_{nk}$ on the unit disk
$\Omega$, one introduces the subdomain $\Omega_{nk} = \{ x\in\R^2~:~
|x| < d_n/\alpha_{nk}\} \subset \Omega$, where $d_n = n - n^{2/3}$,
and $\alpha_{nk}$ are, depending on boundary conditions, the positive
zeros of either $J_n(z)$ (Dirichlet), or $J'_n(z)$ (Neumann) or
$J'_n(z) + h J_n(z)$ for some $h > 0$ (Robin), with $n = 0,1,2,...$
denoting the order of Bessel function $J_n(z)$ and $k=1,2,3,...$
counting zeros.  Then for any $p \geq 1$ (including $p = \infty$),
there exists a universal constant $c_p > 0$ such that for any $k =
1,2,3,...$ and any large enough $n$, the Laplacian eigenfunction
$u_{nk}$ for Dirichlet, Neumann or Robin boundary condition satisfies
\begin{equation} 
\label{eq:eigen_disk_L2}
\frac{\|u_{nk}\|_{L_p(\Omega_{nk})}}{\|u_{nk}\|_{L_p(\Omega)}} < c_p n^{\frac{1}{3} + \frac{2}{3p}} \exp(-n^{1/3}\ln(2)/3) .
\end{equation}
The definition of $\Omega_{nk}$ and the above estimate imply
\begin{equation}
\lim\limits_{n \to \infty} \frac{\|u_{nk}\|_{L_p(\Omega_{nk})}}{\|u_{nk}\|_{L_p(\Omega)}} = 0,  \hskip 10mm 
\lim\limits_{n \to \infty} {\frac{\mu_2(\Omega_{nk})}{\mu_2(\Omega)}} = 1.
\end{equation}
This theorem shows the existence of infinitely many Laplacian
eigenmodes which are $L_p$-localized in a thin layer near the boundary
$\partial \Omega$.  In fact, for any prescribed thresholds for both
ratios in (\ref{eq:def_loc}), there exists $n_0$ such that for all $n
> n_0$, the eigenfunctions $u_{nk}$ are $L_p$-localized.  These
eigenfunctions are called ``whispering gallery eigenmodes'' and
illustrated on Fig. \ref{fig:whispering}.  

We outline a peculiar relation between high-frequency and
low-frequency localization.  The explicit form (\ref{eq:u_disk}) of
Dirichlet Laplacian eigenfunctions $u_{nk}$ leads to their simple
nodal structure which is formed by $2n$ radial nodal lines and $k-1$
circular nodal lines.  The radial nodal lines split the disk into $2n$
circular sectors with Dirichlet boundary conditions.  As a
consequence, whispering gallery eigenmodes in the disk and the
underlying exponential estimate (\ref{eq:eigen_disk_L2}) turn out to
be related to the exponential decay of eigenfunctions in domains with
branches (Sec. \ref{sec:expon}), as illustrated on
Fig. \ref{fig:eigen_triangle} for elongated triangles.

A simple consequence of the above theorem is that for any $p\geq 1$
and any open subset $V$ compactly included in the unit disk $\Omega$
(i.e., $\bar{V}\cap \partial \Omega = \emptyset$), one has
\begin{equation}
\label{eq:disk_local}
\lim\limits_{n \to \infty} \frac{\|u_{nk}\|_{L_p(V)}}{\|u_{nk}\|_{L_p(\Omega)}} = 0,
\end{equation}
and
\begin{equation}
\label{eq:CV}
C_p(V) \equiv \inf\limits_{nk} \left\{ \frac{\|u_{nk}\|_{L_p(V)}}{\|u_{nk}\|_{L_p(\Omega)}} \right\} = 0 .
\end{equation}
Qualitatively, for any subset $V$, there exists a sequence of
eigenfunctions that progressively ``escape'' $V$.

\begin{figure}
\begin{center}
\includegraphics[width=120mm]{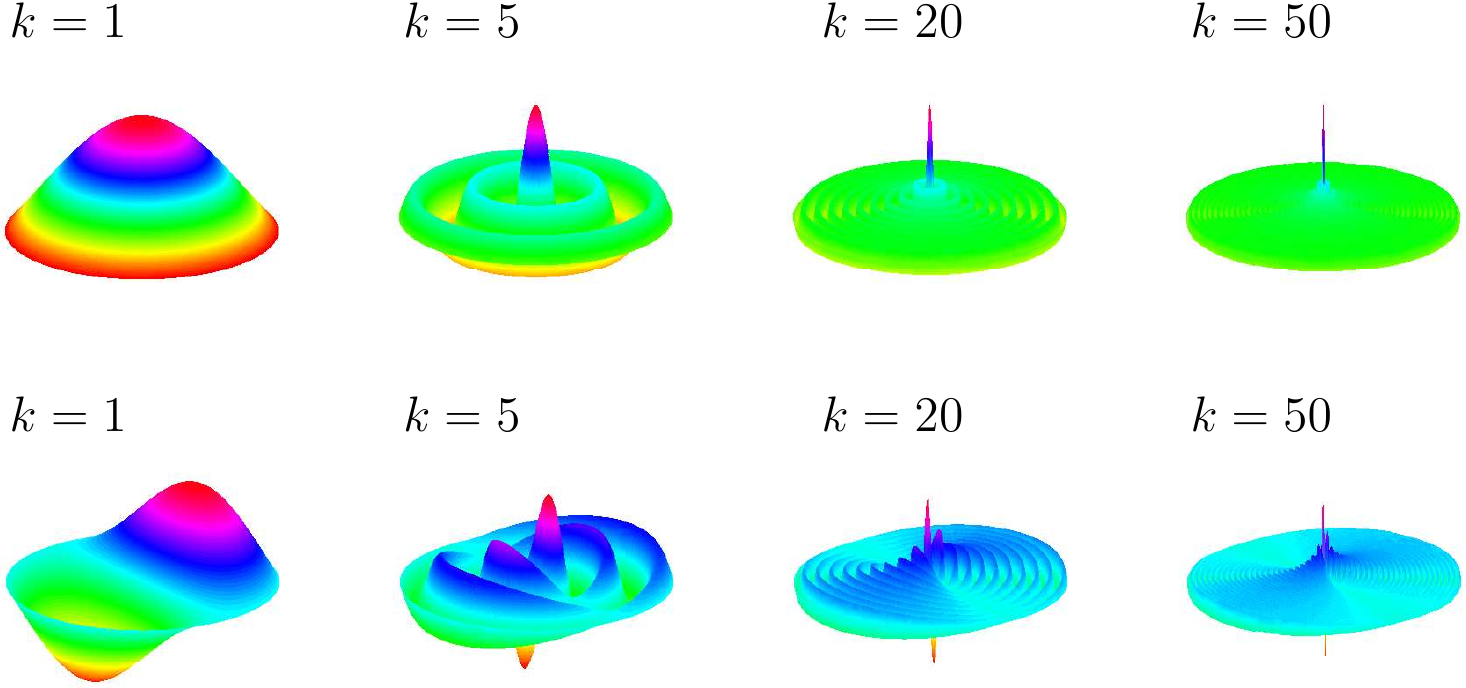}
\end{center}
\caption{
Formation of focusing modes for the unit disk with Dirichlet boundary
condition: for a fixed $n$ ($n = 0$ for top figures and $n = 1$ for
bottom figures), an increase of the index $k$ leads to stronger
localization of eigenfunctions at the origin. }
\label{fig:focusing}
\end{figure}

The localization of focusing modes at the origin is revealed in the
limit $k\to\infty$.  For each $R\in (0,1)$, one defines an annulus
$\Omega_R = \{ x\in\R^2 ~:~ R < |x| < 1 \} \subset \Omega$ of the unit
disk $\Omega$.  Then, for any $n = 0,1,2,...$, the Laplacian
eigenfunction $u_{nk}$ with Dirichlet, Neumann or Robin boundary
condition satisfies
\begin{equation}
\label{eq:focusingNorminfinity}
\lim\limits_{k\to \infty} \frac{\|u_{nk}\|_{L_p(\Omega_R)} }{ \|u_{nk}\|_{L_p(\Omega) }} = \begin{cases} (1 - R^{2-p/2})^{1/p} \quad (1\leq p < 4), \cr
\hskip 10mm  0 \hskip 22mm (p > 4). \end{cases} 
\end{equation}
When the index $k$ increases (with fixed $n$), the eigenfunctions
$u_{nk}$ become more and more $L_p$-localized near the origin when $p
> 4$ \cite{Nguyen12}.  These eigenfunctions are called ``focusing
eigenmodes'' and illustrated on Fig. \ref{fig:focusing}.  The theorem
shows that the definition of localization is sensitive to the norm:
the above focusing modes are $L_p$-localized for $p > 4$, but they are
not $L_p$-localized for $p < 4$.  Similar results for whispering
gallery and focusing modes hold for a ball in three dimensions
\cite{Nguyen12}.

\subsubsection{Bouncing ball modes}

\begin{figure}
\begin{center}
\includegraphics[width=120mm]{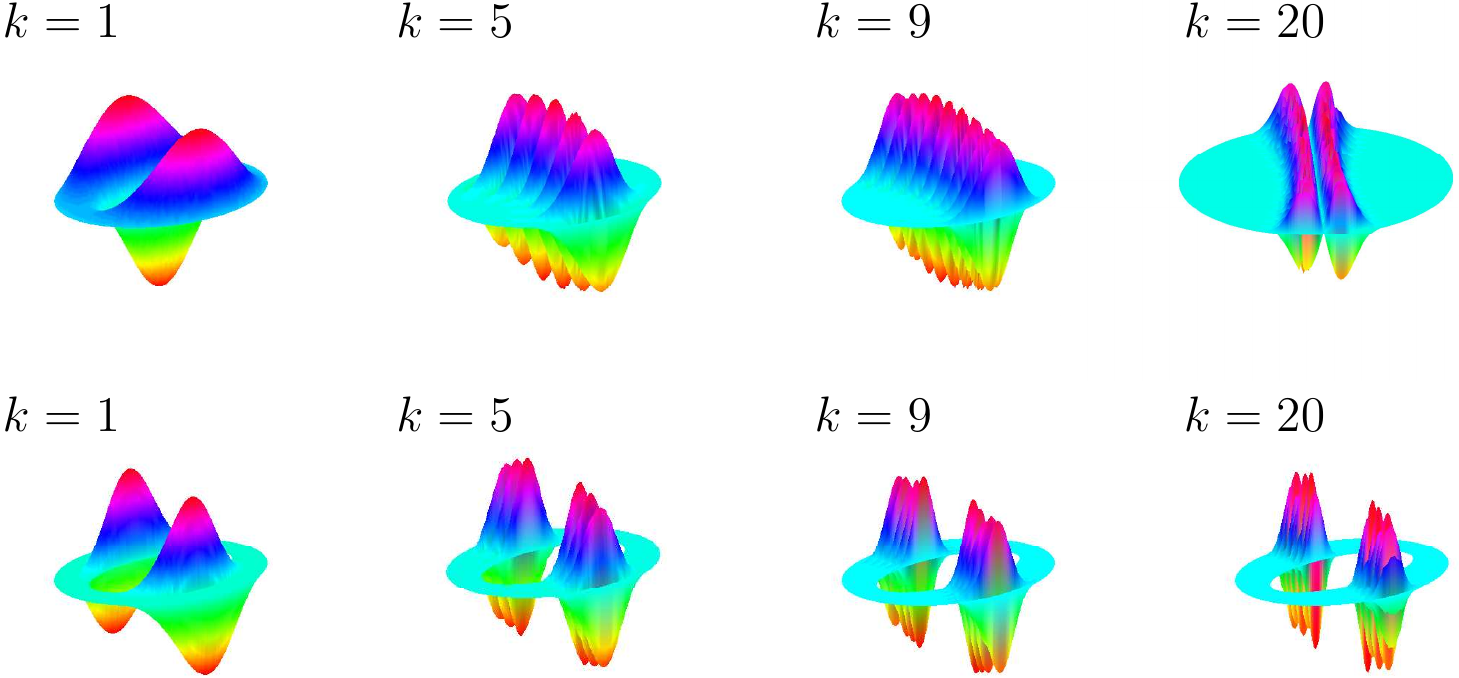}
\end{center}
\caption{
Formation of bouncing ball modes $u_{nkl}$ in a filled ellipse of
radius $R=1$ (top) and an elliptical annulus of radii $R_0 = 0.5$ and
$R = 1$ (bottom), with the focal distance $a = 1$, and Dirichlet
boundary condition.  For fixed $n = 1$ and $l = 1$, an increase of the
index $k$ leads to stronger localization of the eigenfunction near the
vertical semi-axis (from \cite{Nguyen12}).}
\label{fig:eigen_bouncing}
\end{figure}

Filled ellipses and elliptical annuli are simple domains for
illustrating bouncing ball modes.  For fixed foci (i.e., a given
parameter $a$ in the elliptic coordinates in Eq. (\ref{eq:elliptic})),
these domains are characterized by two radii, $R_0$ ($R_0 = 0$ for
filled ellipses) and $R$, as $\Omega = \{ (r,\theta)~:~ R_0 < r < R,~
0\leq \theta < 2\pi\}$, while the eigenfunctions $u_{nkl}$ were
defined in Sec. \ref{sec:ellipse}.  For each $\alpha \in
\left(0,\frac{\pi}{2}\right)$, we consider an elliptical sector
$\Omega_\alpha$ inside an elliptical domain $\Omega$
(Fig. \ref{fig:ellipse})
\begin{equation*}
\Omega_\alpha =  \left\{ (r,\theta) ~:~ R_0 < r < R,~ \theta \in (\alpha,\pi-\alpha) \cup (\pi+\alpha,2\pi-\alpha)\right\}.
\end{equation*}
For any $p \geq 1$, there exists $\Lambda_{\alpha,n} > 0$ such that
for any eigenvalue $\lambda_{nkl} > \Lambda_{\alpha,n}$, the
corresponding eigenfunction $u_{nkl}$ satisfies \cite{Nguyen12} (see
also \cite{Betcke11})
\begin{equation}
\label{eq:bound_elliptical}
\frac{\left\|u_{nkl}\right\|_{L_p(\Omega\setminus \Omega_\alpha)}}{\left\|u_{nkl}\right\|_{L_p(\Omega)}} < 
D_n \left(\frac{16\alpha}{\pi-\alpha/2}\right)^{1/p}
\exp\left(-a\sqrt{\lambda_{nkl}} \left[\sin\left(\frac{\pi}{4}+\frac{\alpha}{2}\right) - \sin\alpha\right]\right),
\end{equation}
where
\begin{equation*}
D_n = 3\sqrt{\frac{1+\sin\left({\frac{3\pi}{8}+\frac{\alpha}{4}}\right)}
{\bigl[\tan\left({\frac{\pi}{16}-\frac{\alpha}{8}}\right)\bigr]^n}} .
\end{equation*}
Given that $\lambda_{nkl}\to\infty$ as $k$ increases (for any fixed
$n$ and $l$), while the area of $\Omega_\alpha$ can be made
arbitrarily small by sending $\alpha\to \pi/2$, the estimate implies
that there are infinitely many eigenfunctions $u_{nkl}$ which are
$L_p$-localized in the elliptical sector $\Omega_\alpha$:
\begin{equation}
\lim\limits_{k \to \infty} \frac{\|u_{nkl}\|_{L_p(\Omega\setminus \Omega_\alpha)}}
{\|u_{nkl}\|_{L_p(\Omega)}} = 0. 
\end{equation}
These eigenfunctions, which are localized near the minor axis, are
called ``bouncing ball modes'' and illustrated on
Fig. \ref{fig:eigen_bouncing}.  The above estimate allows us to
illustrate bouncing ball modes which are known to emerge for any
convex planar domain with smooth boundary \cite{Keller60,Chen94}.  At
the same time, the estimate is as well applicable to elliptical
annuli, providing thus an example of bouncing ball modes for
non-convex domains.

\subsubsection{Domains without localization}

The analysis of geometrical properties of eigenfunctions in
rectangle-like domains $\Omega = (0,\ell_1)\times ... \times
(0,\ell_d) \subset \R^d$ (with sizes $\ell_1 > 0$, ..., $\ell_d > 0$)
may seem to be the simplest case because the eigenfunctions are
expressed through sines (Dirichlet) and cosines (Neumann), as
discussed in Sec. \ref{sec:u_rectangle}.  The situation is indeed
elementary when all eigenvalues are simple, i.e.,  $(\ell_i/\ell_j)^2$
are not rational numbers for all $i\ne j$.  For any $p\geq 1$ and any
open subset $V\subset \Omega$, one can prove that \cite{Nguyen12}
\begin{equation}
\label{eq:C2V}
C_p(V) \equiv \inf\limits_{n_1,...,n_d} \left\{ \frac{\|u_{n_1,...,n_d}\|_{L_p(V)}}{\|u_{n_1,...,n_d}\|_{L_p(\Omega)}} \right\} > 0 .
\end{equation}
This property is in sharp contrast to Eq. (\ref{eq:CV}) for
eigenfunctions in the unit disk (or ball).  The fact that $C_p(V) > 0$
for any open subset $V$ means that there is no eigenfunction that
could fully ``avoid'' any location inside the domain, i.e., there is
no $L_p$-localized eigenfunction.  Since the set of rational numbers
has zero Lebesgue measure, there is no $L_p$-localized eigenfunctions
in almost any randomly chosen rectangle-like domain.

When at least one ratio $(\ell_i/\ell_j)^2$ is rational, certain
eigenvalues are degenerate, and the associated eigenfunctions are
linear combinations of products of sines or cosines (see
Sec. \ref{sec:u_rectangle}).  Although the computation is still
elementary for each eigenfunction, it is unknown whether the infimum
$C_p(V)$ from Eq. (\ref{eq:C2V}) is strictly positive or not, for
arbitrary rectangle-like domain $\Omega$ and any open subset $V$.  For
instance, the most general known result for a rectangle $\Omega =
(0,\ell_1)\times (0,\ell_2)$ states that $C_2(V) > 0$ for any
$V\subset \Omega$ of the form $V = (0,\ell_1)\times \omega$, where
$\omega$ is any open subset of $(0,\ell_2)$ \cite{Burq05} (see also
\cite{Marklof12}).  Even for the unit square, the statement $C_p(V) >
0$ for any open subset $V$ appears as an open problem.  More
generally, one may wonder whether $C_p(V)$ is strictly positive or not
for any open subset $V$ in polygonal (convex) domains.

\subsubsection{Quantum billiards}
\label{sec:quantum_billiard}

\begin{figure}
\begin{center}
\includegraphics[width=100mm]{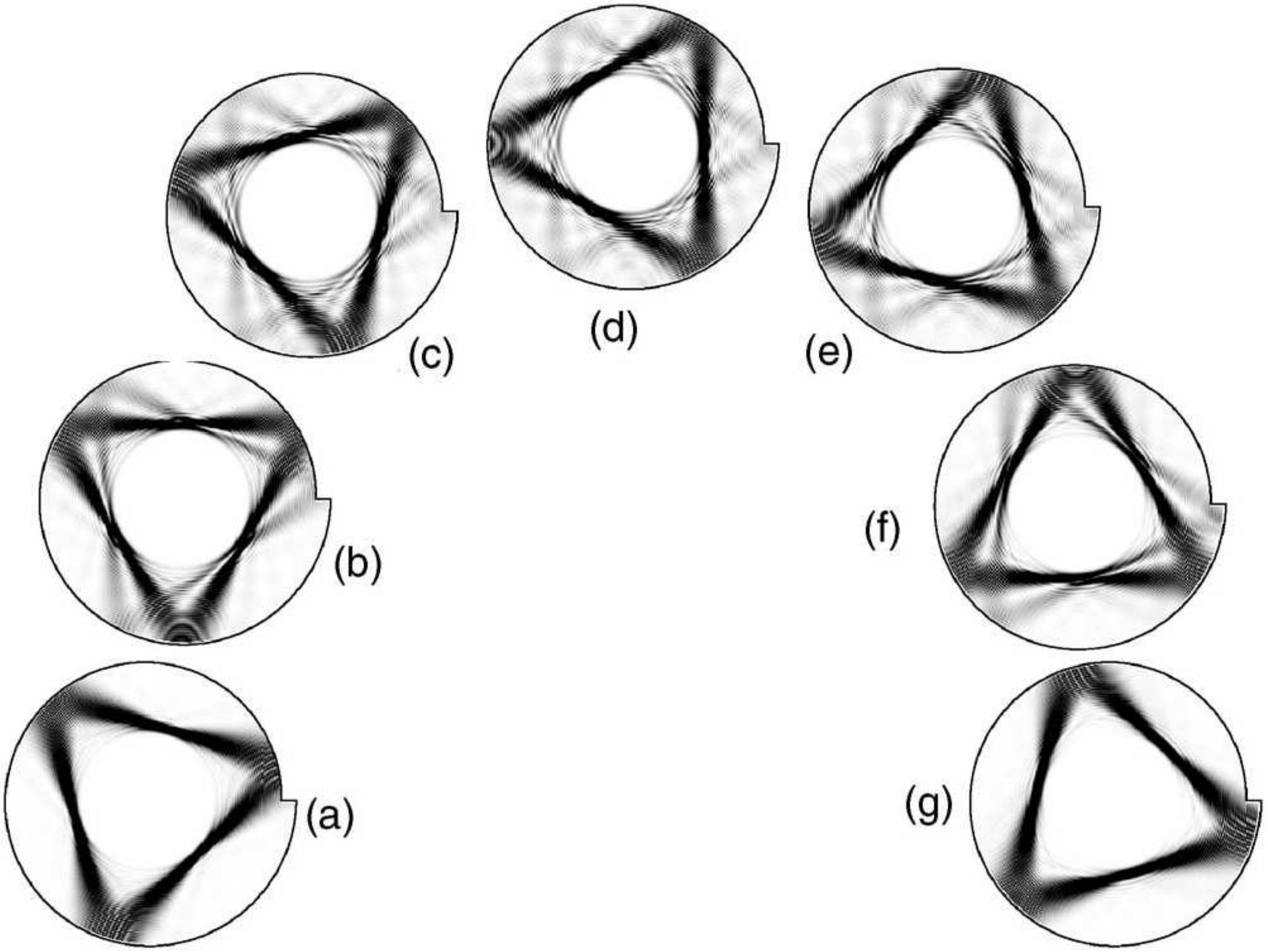}
\end{center}
\caption{
Several eigenstates with localization on periodic orbits for the
spiral-shaped billiard with $\epsilon = 0.1$, from \cite{Liu06}
(Reprinted Fig. 6 with permission from Liu {\it et al.}, Physical
Reviews E, 74, 046214, (2006).  Copyright (2006) by the American
Physical Society).}
\label{fig:liu}
\end{figure}

The above examples of whispering gallery or bouncing ball modes
illustrate that certain high-frequency eigenfunctions tend to be
localized in specific regions of circular and elliptical domains.  But
what is the structure of a high-frequency eigenfunction in a general
domain?  What are these specific regions on which a sequence of
eigenfunctions may be localized, and whether do localized
eigenfunctions exist for a given domain?  Answers to these and other
related questions can be found by relating the high-frequency behavior
of a quantum system (in our case, the structure of Laplacian
eigenfunctions) to the classical dynamics in a billiard of the same
shape \cite{Arnold,Lichtenberg,Gutkin86,Chernov,Katok}.  This relation
is also known as a semi-classical approximation in quantum mechanics
and a ray approximation of wave propagation in optics, while the
correspondence between classical and quantum systems can be shown by
the WKB method, Eikonal equation or Feynman path integrals
\cite{Morse,Feynman48,Feynman}.  For instance, the dynamics of a
particle in a classical billiard is translated into quantum mechanism
through the stationary Schr\"odinger equation $H u_n(x) = E_n u_n(x)$
with Dirichlet boundary condition, where the free Hamiltonian is $H =
p^2/(2m) = - \hbar^2 \Delta/(2m)$, and the energy $E_n$ is related to
the corresponding Laplacian eigenvalue $\lambda_n = 2m E_n/\hbar^2$.
Since $|u_n(x)|^2$ is the probability density for finding a quantum
particle at $x$, this density should resemble some classical
trajectory of that particle in the (high-frequency) semi-classical
limit ($\hbar\to 0$ or $m\to \infty$).
In particular, some orbits of a particle moving in a classical
billiard may appear as ``scars'' in the spatial structure of
eigenfunctions in the related quantum billiard
\cite{Heller84,Gutzwiller,Eckhardt88,Berry89,Aurich91,Steiner94,Sarnak97,Kaplan98,Heller98b,Kaplan99,Stockmann,Haake,Jakobson01,Sarnak11}.
This effect is illustrated on Fig. \ref{fig:liu} by Liu and co-workers
who investigated the localization of Dirichlet Laplacian
eigenfunctions on classical periodic orbits in a spiral-shaped billiard
\cite{Liu06} (see also \cite{Lee04}).

\begin{figure}
\begin{center}
\includegraphics[width=120mm]{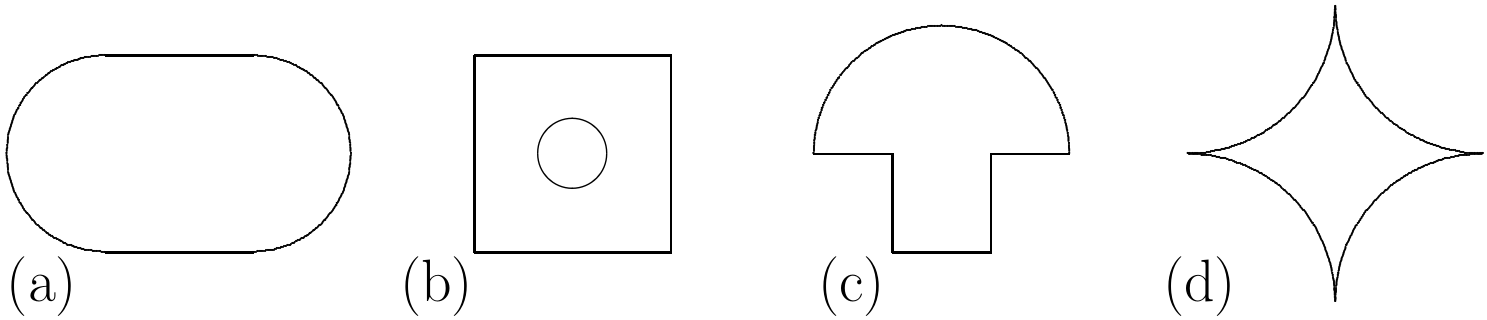}
\end{center}
\caption{
Examples of chaotic billiards: (a) Bunimovich stadium (union of a
square and two half-disks)
\cite{Bunimovich74,Bunimovich79,Heller84,Tomsovic93a,Tomsovic93b,Chinnery96,Bies01},
(b) Sinai's billiard \cite{Sinai63,Sinai70}, (c) mushroom billiard
\cite{Bunimovich01,Barnett07}, and (d) hyperbolic billiard
\cite{Agam94}.  Many other examples are given in
\cite{Bunimovich79}.}
\label{fig:billiards}
\end{figure}

In the classical dynamics, one may distinguish the domains with
regular, integrable and chaotic dynamics.  In particular, for a
bounded domain $\Omega$ with an ergodic billiard flow \cite{Sinai},
Shnirelman's theorem (also known as quantum ergodicity theorem
\cite{Colin85,Zelditch96,Zelditch96b}) states that among the set of
$L_2$-normalized Dirichlet (or Neumann) Laplacian eigenfunctions,
there is a sequence $u_{j_k}$ of density $1$ (i.e., $\lim\limits_{k\to
\infty} {j_k/k} = 1$), such that for any open subset $V \subset
\Omega$, one has \cite{Shnirelman74}
\begin{equation}
\lim\limits_{k\rightarrow \infty} \int\limits_V |u_{j_k}(x)|^2 dx = \frac{\mu_d(V)}{\mu_d(\Omega)} .
\end{equation} 
(this version of the theorem was formulated in \cite{Burq05}).
Marklof and Rudnick extended this theorem to rational polygons
(i.e., simple plane polygons whose interior is connected and simply
connected and all the vertex angles are rational multiplies of $\pi$)
\cite{Marklof12}.  Loosely speaking, $\{ u_{j_k} \}$ is a sequence of
non-localized eigenfunctions which become more and more uniformly
distributed over the domain (see
\cite{Jakobson01,Burq05,Gerard93} for further discussion and
references).  At the same time, this theorem does not prevent the
existence of localized eigenfunctions.  How large the excluded
subsequence of (localized) eigenfunctions may be?  In the special case
of arithmetic hyperbolic manifolds, Rudnick and Sarnak proved that
there is no such excluded subsequence \cite{Rudnick94}.  This
statement is known as the quantum unique ergodicity (QUE).  The
validity of this statement for other dynamical systems (in particular,
ergodic billiards) remains under investigation
\cite{Donnelly03,Barnett06,Hassell10}.  The related notion of weak
quantum ergodicity was discussed by Kaplan and Heller
\cite{Kaplan98b}.  A classification of eigenstates to regular and
irregular ones was thoroughly discussed (see \cite{Prosen93,Veble99}
and references therein).

There were numerous studies of Laplacian eigenfunctions in chaotic
domains such as, e.g., Bunimovich stadium
\cite{Bunimovich74,Bunimovich79,Heller84,OConnor88,Tomsovic93a,Tomsovic93b,Chinnery96,Bies01},
Sinai's billiard \cite{Sinai63,Sinai70}, mushroom billiard
\cite{Bunimovich01,Barnett07} or hyperbolic billiard \cite{Agam94},
illustrated on Fig. \ref{fig:billiards}.  The literature on quantum
billiards is vast, and we only mention selected works on the spatial
structure of high-frequency eigenfunctions.  McDonald and Kaufman
studied the Bunimovich stadium billiard and reported a random
structure of nodal lines of eigenfunctions and Wigner-type
distribution for eigenvalue spacings \cite{McDonald79,McDonald88}.
Bohigas and co-workers analyzed eigenvalue spacings for the Sinai's
billiard and also obtained the Wigner-type distribution
\cite{Bohigas84}.  It means that eigenvalue spacings for these
chaotic billiards obey the same distribution as that for random
matrices from the Gaussian Orthogonal Ensemble.  This is in a sharp
contrast to regular billiards for which eigenvalue spacings generally
follow a Poisson distribution.  The problem of circular-sector and
related billiards was studied (e.g., see \cite{Liboff94b}).

Polygon billiards have attracted a lot of attention, especially the
class of rational polygons for which all the vertex angles are
rational multiples of $\pi$ \cite{Richens81,Liboff94c,Liboff01}.  As
the dynamics in rational polygons is neither integrable nor ergodic
(except several classical integrable cases such as rectangles,
equilateral triangle, right triangles with an acute vertex angle
$\pi/3$ or $\pi/4$), they are often called {\it pseudo-integrable}
systems.  Bellomo and Uzer studied scarring states in a
pseudo-integrable triangular billiard and detected scars in regions
which contain no periodic orbits \cite{Bellomo94}.  Amar {\it et al.}
gave a complete characterization of the polygons for which a Dirichlet
eigenfunction can be found in terms of a finite superposition of plane
waves \cite{Amar91,Amar93} (see also \cite{Lauber94} for experimental
study).  Biswas and Jain investigated in detail the $\pi/3$-rhombus
billiard which presents an example of the simplest pseudo-integrable
system \cite{Biswas90}.  Hassell {\it et al.} proved for an arbitrary
polygonal billiard that eigenfunction mass cannot concentrate away
from the vertices \cite{Hassell09} (see also \cite{Marzuola07}).  The
level spacing properties of rational and irrational polygons were
studied numerically by Shimizu and Shudo \cite{Shimizu93}.  They also
analyzed the structure of the related eigenfunctions \cite{Shimizu95}.

\begin{figure}
\begin{center}
\includegraphics[width=120mm]{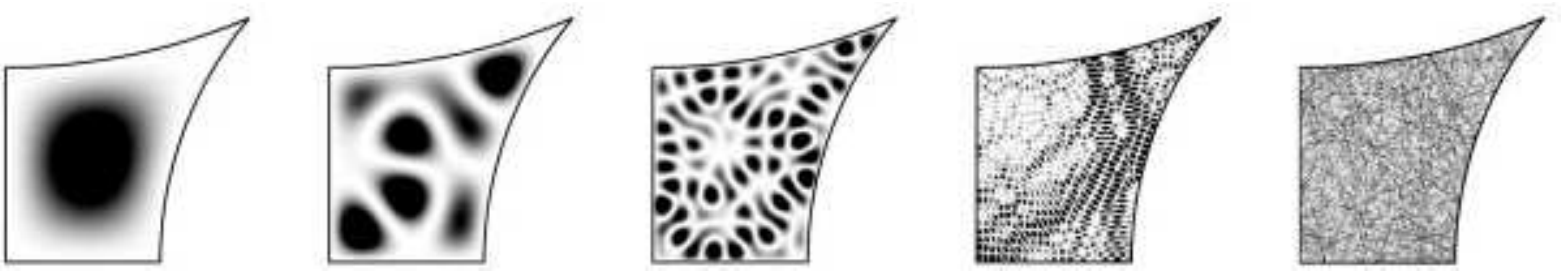}
\end{center}
\caption{
Illustration of spatial distribution of the Dirichlet Laplacian
eigenfunction $|u_m|^2$ (shown as density plots: larger values are
darker) with $m = 1, 10, 100, 1000$ and $m\approx 50 000$
\cite{Barnett06} (by A. Barnett, with permission).  }
\label{fig:barnett}
\end{figure}

\begin{figure}
\begin{center}
\includegraphics[width=62mm]{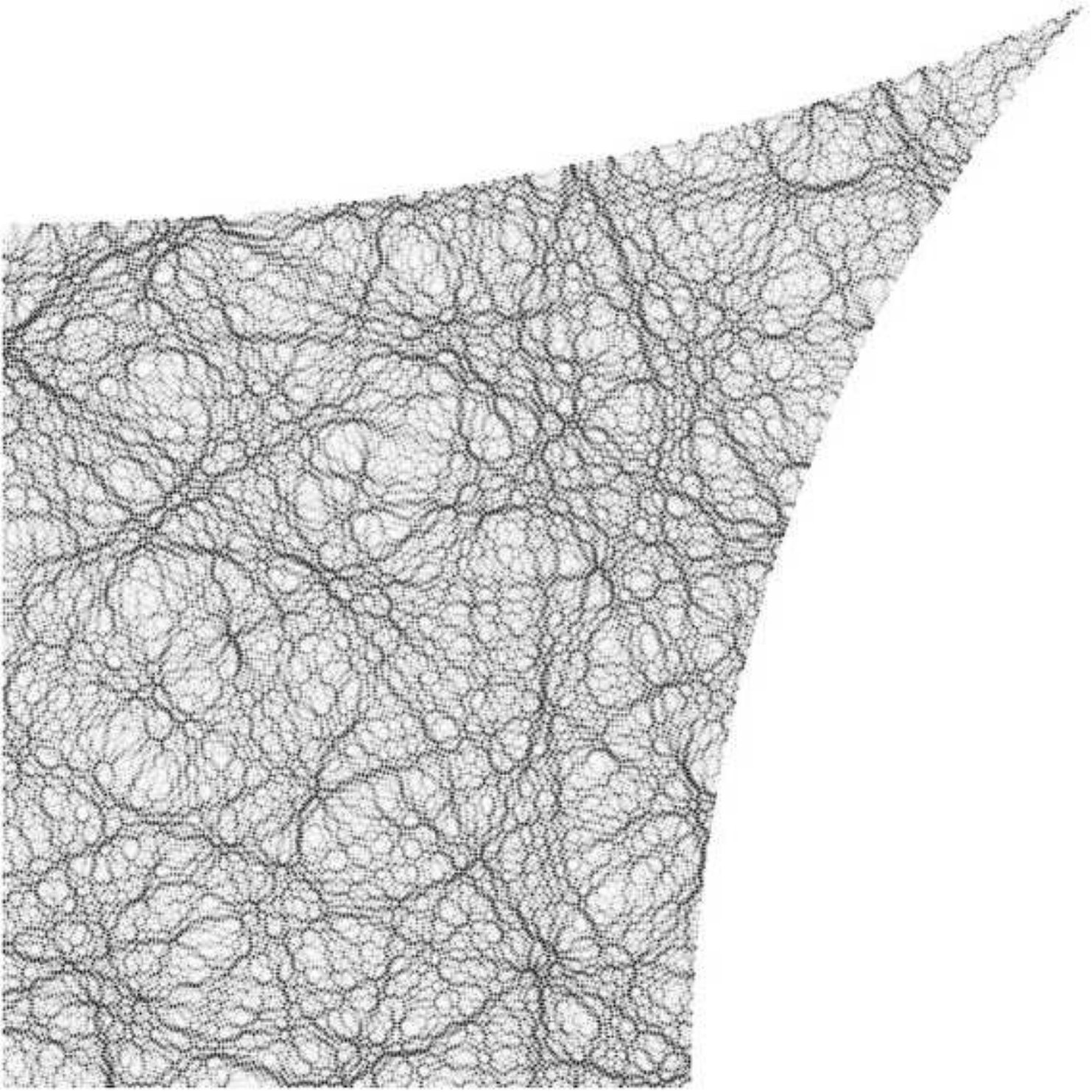} 
\includegraphics[width=62mm]{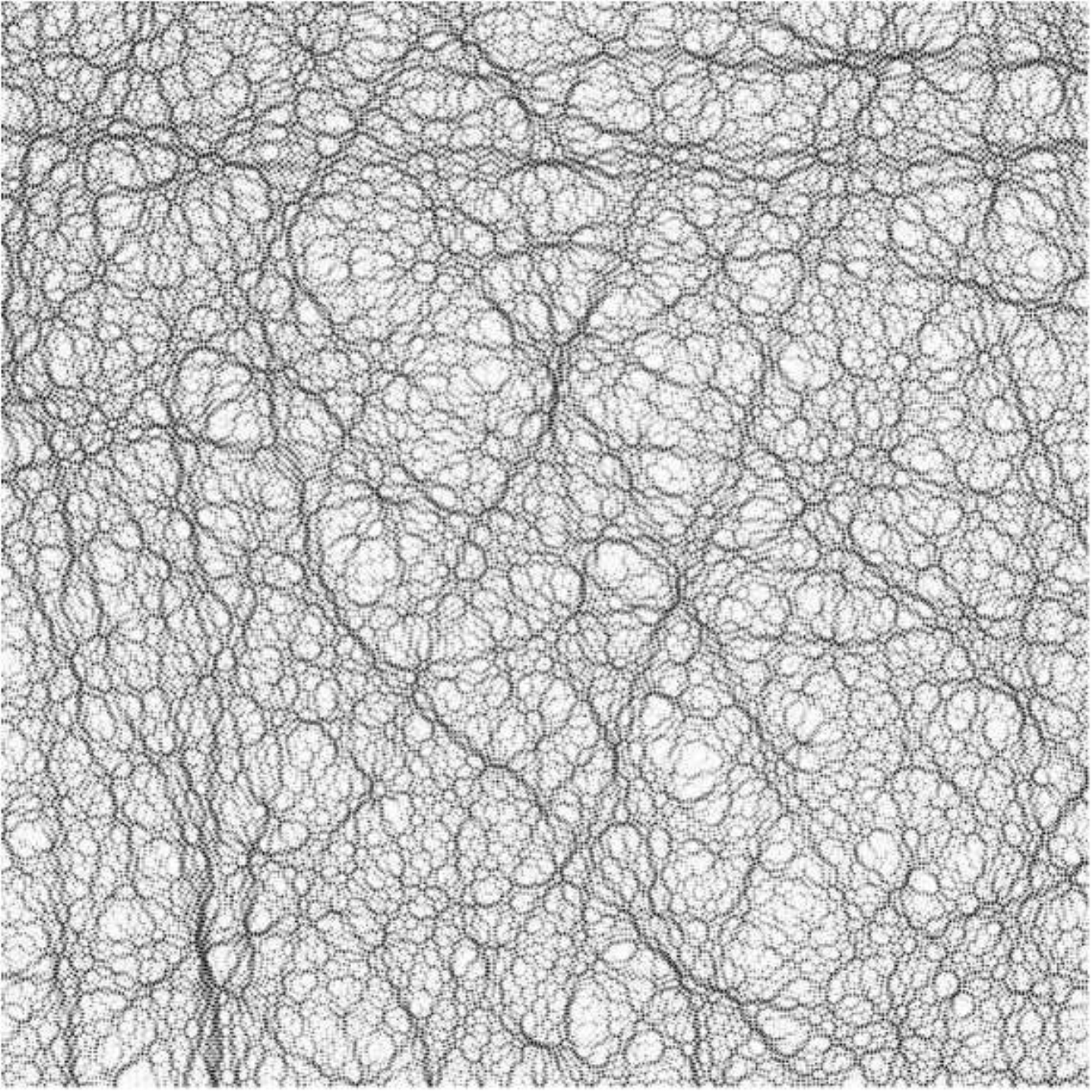}
\end{center}
\caption{
{\bf (Left)} Density plot of a Dirichlet-Laplacian eigenfunction
$|u_m|^2$ for $m\approx 50 000$ with the eigenvalue $\lambda_m \approx
10^6$.  There are about 225 wavelengths across the diagonal; {\bf
(Right)} Density plot of one sample from the ensemble of random plane
waves with the same wavenumber and mean intensity, shown in a square
region of space (with no boundary conditions) \cite{Barnett06} (by
A. Barnett, with permission).  }
\label{fig:barnett2}
\end{figure}

B\"acker and co-workers analyzed the number of bouncing ball modes in
a class of two-dimensional quantized billiards with two parallel walls
\cite{Backer97}.  Bunimovich introduced a family of simple billiards
(called ``mushrooms'') that demonstrate a continuous transition from a
completely chaotic system (stadium) to a completely integrable one
(circle) \cite{Bunimovich01}.  Barnett and Betcke reported the first
large-scale statistical study of very high-frequency eigenfunctions in
these billiards \cite{Barnett07}.  Using nonstandard numerical
techniques \cite{Barnett06}, Barnett also studied the rate of
equidistribution for a uniformly hyperbolic, Sinai-type, planar
Euclidean billiard with Dirichlet boundary condition, as illustrated
on Fig. \ref{fig:barnett}.  This study brought a strong numerical
evidence for the QUE in this system.  The spatial structure of
high-frequency eigenfunctions shown on Fig. \ref{fig:barnett} looks
somewhat random.  This observation goes back to Berry who conjectured
that high-frequency eigenfunctions in domains with ergodic flow should
look locally like a random superposition of plane waves with a fixed
wavenumber \cite{Berry77}.  This analogy is illustrated on
Fig. \ref{fig:barnett2} by Barnett \cite{Barnett06}.  O'Connor and
co-workers analyzed the random pattern of ridges in a random
superposition of plane waves \cite{OConnor87}.

Pseudo-integrable barrier billiards were intensively studied in a
series of theoretical, numerical and experimental works by Bogomolny
{\it et al.} \cite{Bogomolny04,Bogomolny06}.  They reported on the
emergence of scarring eigenstates which are related with families of
classical periodic orbits that do not disappear at large quantum
numbers in contrast to the case of chaotic systems.  These so-called
superscars were observed experimentally in a flat microwave billiard
with a barrier inside \cite{Bogomolny06}.  Wiersig performed an
extensive numerical study of nearest-neighbor spacing distributions,
next-to-nearest spacing distributions, number variances, spectral form
factors, and the level dynamics \cite{Wiersig02}.  Dietz and
co-workers analyzed the number of nodal domains in a barrier billiard
\cite{Dietz08}. 

Tomsovic and Heller verified a remarkable accuracy of the
semi-classical approximation that relates the classical and quantum
billiards \cite{Tomsovic93a,Tomsovic93b}.  In some cases,
eigenfunctions can therefore be constructed by purely semiclassical
calculations.  Li {\it et al.} studied the spatial distribution of
eigenstates of a rippled billiard with sinusoidal walls \cite{Li02}.
For this type of ripple billiards, a Hamiltonian matrix can be found
exactly in terms of elementary functions that greatly improves
computation efficiency.  They found both localized and extended
eigenfunctions, as well as peculiar hexagon and circle-like pattern
formations.  Frahm and Shepelyansky considered almost circular
billiards with a rough boundary which was realized as a random curve
with some finite correlation length. On a first glance it may seem
that such a rough boundary in a circular billiard would destroy the
conservation of angular momentum and lead to ergodic eigenstates and
the level statistics predicted by random matrix theory.  They showed,
however, that there is a region of roughness in which the classical
dynamics is chaotic but the eigenstates are localized \cite{Frahm97}.
Bogomolny {\it et al.}  presented the exact computation of the
nearest-neighbor spacing distribution for a rectangular billiard with
a pointlike scatterer inside \cite{Bogomolny02b}.

Prosen computed numerically very high-lying energy spectra for a
generic chaotic 3D quantum billiard (a smooth deformation of a unit
sphere) and analyzed Weyl's asymptotic formula and the nearest
neighbor level spacing distribution.  He found significant deviations
from the Gaussian Orthogonal Ensemble statistics that were explained
in terms of localization of eigenfunctions onto lower dimensional
classically invariant manifolds \cite{Prosen97a}. He also found that
the majority of eigenstates were more or less uniformly extended over
the entire energy surface, except for a fraction of strongly localized
scarred eigenstates \cite{Prosen97b}.  An extensive study of 3D
Sinai's billiard was reported by Primack and Smilansky
\cite{Primack00}.  Deviations from a semi-classical description were
discussed by Tanner \cite{Tanner97}.  Casati and co-workers
investigated how the interplay between quantum localization and the
rich structure of the classical phase space influences the quantum
dynamics, with applications to hydrogen atoms under microwave fields
\cite{Casati,Casati87,Casati95,Casati99} (see also references therein).

A large number of physical experiments were performed with classical
and quantum billiards.  For instance, Gr\"af and co-workers measured
more than thousand first eigenmodes in a quasi two-dimensional
superconducting microwave stadium billiard with chaotic dynamics
\cite{Graf92}.  Sridhar and co-workers performed a series of
experiments in microwave cavities in the shape of Sinai's billiard
\cite{Sridhar91,Sridhar92}.  In particular, they observed bouncing ball
modes and modes with quasi-rectangular or quasi-circular symmetry
which are associated with nonisolated periodic orbits (which avoid the
central disk).  Some scarring eigenstates, which are associated with
isolated periodic orbits (which hit the central disk, see
Fig. \ref{fig:billiards}b), were also observed.  Kudrolli {\it et al.}
investigated the signatures of classical chaos and the role of
periodic orbits in the eigenvalue spectra of two-dimensional billiards
through experiments in microwave cavities
\cite{Kudrolli94,Kudrolli95}.  The eigenvalue spectra were analyzed by
using the nearest neighbor spacing distribution for short-range
correlations and the spectral rigidity for longer-range correlations.
The density correlation function was used for studying the spatial
structure of eigenstates.  The role of disorder was also investigated.
Chinnery and Humphrey visualized experimentally acoustic resonances
within a stadium-shaped cavity \cite{Chinnery96}.  Bittner {\it et
al.} performed double-slit experiments with regular and chaotic
microwave billiards \cite{Bittner11}.  Chaotic resonators were also
employed for getting specific properties of lasers (e.g., high-power
directional emission or ``Fresnel filtering'') \cite{Gmachl98,Rex02}.

\section{Other points and concluding remarks}
\label{sec:conclusion}

This review was focused on the geometrical properties of Laplacian
eigenfunctions in Euclidean domains.  We started from the basic
properties of the Laplace operator and explicit representations of its
eigenfunctions in simple domains.  After that, the properties of
eigenvalues and their relation to the shape of a domain were briefly
summarized, including Weyl's asymptotic behavior, isoperimetric
inequalities, and Kac's inverse spectral problem.  The structure of
nodal domains and various estimates for the norms of eigenfunctions
were then presented.  The main Section \ref{sec:localization} was
devoted to the spatial structure of eigenfunctions, with a special
emphasis on their localization in small subsets of a domain.  One of
the major difficulties in the study of localization is that
localization is a property of an individual eigenfunction.  For the
same domain, two consecutive eigenfunctions with very close
eigenvalues may have drastically different geometrical structures
(e.g., one is localized and the other is extended).  One needs
therefore fine analytical tools which would differently operate with
localized and non-localized eigenfunctions.  In the review, we
distinguished two types of localization, for low-frequency and
high-frequency eigenfunctions.

In the former case, an eigenfunction remains localized in a subset
because of a geometric constraint that prohibits its extension to
other parts of the domain.  A standard example is a dumbbell (two
domains connected by a narrow channel), for which an eigenfunction may
be localized in one domain if its typical wavelength is larger than
the width of the channel (meaning that an eigenfunction cannot
``squeeze'' through the channel).  Such kind of ``expulsion'' from a
channel is quite generic, as the analysis is applicable to domains
with branches of variable cross-sectional profiles.  It is important
to note that a geometric constraint does not need to be strong (e.g.,
two domains may be separated by a cloud of point-like obstacles of
zero measure).  Another example is an elongated triangle, in which
there is no ``obstacles'' at all.  Low-frequency localization was
found numerically in many irregularly-shaped domains, for both
Dirichlet and Neumann boundary conditions.  From a practical point of
view, the low-frequency localization may find various applications,
e.g., it is important for the theory of quantum, optical and
acoustical waveguides and microelectronic devices, as well as for
analysis and engineering of highly reflecting or absorbing materials
(noise protective barriers, anti-radar coatings, etc).

The high-frequency localization manifests in quantum billiards when a
sequence of eigenfunctions tends to concentrate onto some orbits of
the associated classical billiard.  In this regime, the asymptotic
properties of eigenvalues and eigenfunctions are strongly related to
the underlying classical dynamics (e.g., regular, integrable or
chaotic).  For instance, the ergodic character of the classical system
is reflected in the spatial structure of eigenfunctions.  Working on
simple domains, we illustrated several kinds of localized
eigenfunctions which emerge for a large class of domains.  We also
provided examples of rectangular domains without localization.
Although a number of rigorous and numerical results were obtained
(e.g., quantum ergodicity theorem for ergodic billiards), many
questions about the spatial structure of high-frequency eigenfunctions
remain open, even for very simple domains (e.g., a square).

Although the review counts more than five hundred citations, it is far
from being complete.  As already mentioned, we focused on the Laplace
operator in bounded Euclidean domains and mostly omitted technical
details, in order to keep the review at a level accessible to
scientists from various fields.  Many other issues had to be omitted:

(i) Many important results for Laplacians on Riemannian manifolds or
weighted graphs could not be included.  In addition, we did not
discuss the spectral properties of domains with ``holes''
\cite{Samarski48,Ozawa81,Mazya84,Mazya00,Ward93,Flucher95,Ozawa96,Schatzman96,Kolokolnikov05,Cheviakov11,Lanza12},
as well as their consequences for diffusion in domains with static
traps \cite{Grassberger82,Kayser83,Kayser84,Torquato91,Nguyen10}.  The
behavior of the eigenvalues and eigenfunctions under deformations of a
domain was partly considered in Sec. \ref{sec:first_eigenfunction} and
Sec. \ref{sec:dumbbell}, while many significant results were not
included (see \cite{Rellich2,Kato,Simon91,Bruno01,Khelifi07} and
references therein).

(ii) There are important developments of numerical techniques for
computing the Laplacian eigenbasis.  In fact, standard finite
difference or finite element methods rely on a regular or adapted
discretization of a domain that reduces the continuous eigenvalue
problem to a finite set of linear equations
\cite{Ciarlet,Saad,Kuttler74,Cureton99,Heuveline03,Fried04,Grinfeld04}.
Since finding the eigenbasis of the resulting matrix is still an
expensive computational task, various hints and tricks are often
implemented.  For instance, for planar polygonal domains, one can
exploit the behavior of eigenfunctions at corners through radial basis
functions in polar coordinates and the integration of related
Fourier-Bessel functions on subdomains
\cite{Descloux83,Driscoll97,Platte04}.  Another ``trick'' is conformal
mapping of planar polygonal domains onto the unit disk, for which the
modified eigenvalue problem can be efficiently solved
\cite{Banjai03,Banjai07}.  Yet another approach, known as the method of
particular solutions was suggested by Fox and co-workers \cite{Fox67}
and later progressively improved \cite{Backer,Betcke05,Barnett09}.
The main idea is to consider various solutions of the eigenvalue
equation for a given value of $\lambda$ and to vary $\lambda$ until a
linear combination of such solutions would satisfy the boundary
condition at a number of sample points along the boundary.  One can
also mention a stochastic method by Lejay and Maire for computing the
principal eigenvalue \cite{Lejay07}.  The eigenvalue problem can also
be reformulated in terms of boundary integral equations that reduces
the dimensionality and allows for rapid computation of eigenvalues
\cite{Lu91}.  Kaufman and co-workers proposed a simple expansion
method in which wave functions inside a two-dimensional quantum
billiard are expressed in terms of an expansion of a complete set of
orthonormal functions defined in a surrounding rectangle for which
Dirichlet boundary conditions apply, while approximating the billiard
boundary by a potential energy step of a sufficiently large size
\cite{Kaufman99}.  Vergini and Saraceno proposed a scaling method for
computing high-frequency eigenmodes \cite{Vergini95}.  This method was
later improved by Barnett and Hassell \cite{Barnett11} (this reference
is also a good review of numerical methods for high-frequency
Dirichlet Laplacian eigenvalues).

(iii) We did not discuss various applications of Laplacian
eigenfunctions which nowadays range from pure and applied mathematics
to physics, chemistry, biology and computer sciences.  One can mention
manifold parameterizations by Laplacian eigenfunctions and heat
kernels \cite{Jones08}, the use of Laplacian spectra as a diagnostic
tool for network structure and dynamics \cite{McGraw08}, efficient
image recognition and analysis
\cite{Reuter06,Saito08,Saito09,Rhouma11}, shape optimization and
spectral partition problems
\cite{Pironneau84,Sokolowski92,Allaire02,Bucur05,Buttazzo11,Caffarelli07,vandenBerg11},
computation and analysis of diffusion-weighted NMR signals
\cite{Grebenkov07,Grebenkov08,Grebenkov09}, localization in
heterogeneous materials (e.g., photonic crystals) and the related
optimization problem
\cite{John87,John91,Figotin96,Figotin97,Figotin97b,Figotin98,Kohler96a,Kohler96b,Kuchment04,Dobson04,Sakoda,Kubytskyi12},
etc.

\section*{Acknowledgments}

The authors thank H. Obuse and A. Barnett for providing their images
of eigenfunctions and allowing us to reproduce them.  We acknowledge
fruitful discussions with Y. G. Sinai and A. L. Delitsyn.  We thank
M. Ashbaugh, A. Barnett, J. Clutterbuck, V. Nistor for their comments
on the original version and for reference updates.  We also appreciate
numerous remarks by the reviewers.  Finally, we are grateful to
B. Sapoval for many valuable discussions and for his passion to
localization that strongly motivated our work.

\end{document}